\documentclass[11pt]{amsart}
\usepackage{amsmath,latexsym,amssymb,amsthm,amscd}
\usepackage{hyperref}
\usepackage{lineno}
\usepackage{enumitem}
\usepackage{chngcntr}
\usepackage{pst-all}
\begin{document}
\newtheorem{theorem}{Theorem}[section]
\newtheorem{definition}[theorem]{Definition}
\newtheorem{proposition}[theorem]{Proposition}
\newtheorem{lemma}[theorem]{Lemma}
\newtheorem{remark}[theorem]{Remark}
\newtheorem{corollary}[theorem]{Corollary}
\newtheorem{question}{Question}
\newtheorem{example}{Examples}[section]
\newtheorem{notation}[theorem]{Notation}
\newtheorem{notation-lemma}[theorem]{Notation-Lemma}
\newtheorem{claim}{Claim}[section]
\newtheorem{fact}[theorem]{Fact}
\newcommand\cl{\begin{claim}}
\newcommand\ecl{\end{claim}}
\newcommand\rem{\begin{remark}\upshape}
\newcommand\erem{\end{remark}}
\newcommand\ex{\begin{example}\upshape}
\newcommand\eex{\end{example}}
\newcommand\nota{\begin{notation}\upshape}
\newcommand\enota{\end{notation}}
\newcommand\notalem{\begin{notation-lemma}\upshape}
\newcommand\enotalem{\end{notation-lemma}}
\newcommand\dfn{\begin{definition}\upshape}
\newcommand\edfn{\end{definition}}
\newcommand\cor{\begin{corollary}}
\newcommand\ecor{\end{corollary}}
\newcommand\thm{\begin{theorem}}
\newcommand\ethm{\end{theorem}}
\newcommand\prop{\begin{proposition}}
\newcommand\eprop{\end{proposition}}
\newcommand\lem{\begin{lemma}}
\newcommand\elem{\end{lemma}}
\newcommand\fct{\begin{fact}}
\newcommand\efct{\end{fact}}
\providecommand\qed{\hfill$\quad\Box$}
\newcommand\pr{{Proof:\;}}
\newcommand\prcl{\par\noindent{\em Proof of Claim: }}
\def\K{\mathcal{K}}
\def\cC{\mathcal C}
\def\cS{\mathcal S}
\newcommand\cP{\mathcal P}
\newcommand\U{\mathcal U}
\newcommand\M{\mathcal M}
\newcommand\X{\mathbf{X}}
\newcommand\y{\mathbf{y}}
\newcommand\x{\mathbf{x}}
\newcommand\ba{\mathbf{a}}
\newcommand\bb{\mathbf{b}}
\newcommand\cc{\mathbf{c}}
\newcommand\bu{\mathbf{u}}
\newcommand\bd{\bar{\delta}}
\newcommand\f{\mathfrak{f}}
\def\N{\mathbb N}
\newcommand\R{\mathbb R}
\newcommand\cR{\mathcal R}
\def\F{\mathbb F}
\def\A{\mathbb A}
\newcommand\Q{\mathbb Q}
\newcommand\Z{\mathbb Z}
\newcommand\cF{\mathcal F}
\newcommand\cL{\mathcal L}
\newcommand\cO{\mathcal O}
\newcommand\cA{\mathcal A}
\newcommand\cB{\mathcal B}
\newcommand\cK{\mathcal K}
\newcommand\cZ{\mathcal Z }
\newcommand\cG{\mathcal G}
\def\T{\mathcal T}
\newcommand\G{\Gamma}
\newcommand\gG{\mathfrak G}
\newcommand\cM{\mathcal M}
\newcommand\cN{\mathcal N}
\newcommand{\acl}{\mathrm{acl}}
\newcommand{\dcl}{\mathrm{dcl}}
\newcommand{\Cl}{\mathrm{cl}}
\newcommand{\graph}{\mathrm{graph}}
\newcommand{\dom}{\mathrm{dom}}
\newcommand{\im}{\mathrm{im}}
\newcommand{\I}{\mathrm{I}}
\newcommand{\DCF}{\mathrm{DCF}}
\newcommand{\RCVF}{\mathrm{RCVF}}
\newcommand{\ACVF}{\mathrm{ACVF}}
\newcommand{\RCF}{\mathrm{RCF}}
\newcommand{\ACF}{\mathrm{ACF}}
\newcommand{\Alg}{\mathrm{Alg}}
\newcommand{\ord}{\text{ord}}
\newcommand{\Ja}{\nabla}
\newcommand{\bF}{\mathbf{F}}
\newcommand{\DL}{\mathrm{DL}}
\newcommand{\0}{\emptyset}
\newcommand{\IFTA}{\mathrm{IFTA}}

\title[Definable groups]{On definable groups in dp-minimal topological fields equipped with a generic derivation}
\author{Fran\c{c}oise Point}
\address{Department of Mathematics (De Vinci)\\ UMons\\ 20, place du Parc 7000 Mons, Belgium}
\email{point@math.univ-paris-diderot.fr}

\begin{abstract}  
Let $T$ be a complete, model-complete,  geometric dp-minimal $\cL$-theory of topological fields of characteristic $0$ and let $T(\partial)$ be the theory of expansions of models of $T$ by a derivation $\partial$. We assume that $T(\partial)$ has a model-companion $T_{\partial}$. 
Let $\Gamma$ be a finite-dimensional $\cL_\partial$-definable group in a model of $T_\partial$. Then we show that $\Gamma$ densely and definably embeds in an $\cL$-definable group $G$.
 Further, using a $C^1$-cell decomposition result, we show that $\Gamma$ densely and definably embeds in a definable $D$-group, generalizing the classical construction of Buium of algebraic $D$-groups and extending for that class of fields, results obtained in \cite{PPP}, \cite{PPP2}. 
\end{abstract}
\maketitle

\section{Introduction}
\par In \cite{S}, M. Singer showed that the class of existentially closed ordered differential fields is elementary and gave a simple axiomatization. Later on, his work was put in a larger setting by M. Tressl who described a uniform companion for every theory of large, differential fields (in the sense of \cite{Pop}) of characteristic zero extending a model-complete theory of pure fields. Theories of fields that enter Tressl setting are, for instance, pseudo-finite fields (of characteristic $0$) and p-adically closed fields \cite[Corollary 8.2]{Tressl}. This was further extended to encompass the theory of algebraically closed valued fields of characteristic $0$ \cite{GP} and then using many-sorted structures the theory of any henselian field of characteristic $0$ \cite{K-P}. 
Then working in  o-minimal expansions of the field of real numbers where the derivation $\partial$ respects definable functions \cite{W} (see section \ref{sec:diff} for a more precise statement), these results were extended by Fornasiero and Kaplan \cite{FK}. Lately, work of 
 Fornasiero and Terzo encompasses the case of classes of non-large geometric fields with a model-complete theory.
 \par The common property of all these examples is that the classes of structures that one deals with, are models of a geometric theory, namely the model-theoretic algebraic closure $\acl$ satisfies the exchange property. 
 \par Let $T$ be a complete, model-complete  geometric $\cL$-theory of fields of characteristic $0$ and let $T(\partial)$ be the theory of expansions of models of $T$ by a derivation $\partial$. We will assume that $T(\partial)$ has a model-companion $T_{\partial}$. In \cite{PPP}, \cite{PPP2}, in all the above mentioned cases, we described definable groups in models of $T_\partial$ in terms of definable groups in models of $T$. The aim of the present work is to generalize some of these former results to the case where $T$ is in addition dp-minimal, not strongly minimal or  $T$ is an open theory of topological fields \cite{K-P}.
\par Given a finite-dimensional $\cL_{\partial}$-definable group $\G$ in a model of $T_{\partial}$, we embed $\G$ {\it densely} in an $\cL$-definable group where {\it dense} means that any generic type is realized by a tuple coming from the group $\G$ (and when the model is a topological field, {\it dense} will have its usual topological meaning). Let $\nabla_{N}(\G):=\{(a,\partial(a),\ldots,\partial^N(a))\colon a\in \G\}$, $N\geq 0$ (with the convention that $\partial^0(a)=a$), the hypothesis that $\G$ is finite-dimensional, means that for some $N\in \N$, $\nabla_{N}(\G)\subset \acl_{\mathcal L}(\nabla_{N-1}(\G))$. In a subsequent paper, the case where $\G$ is not finite-dimensional is being dealt with \cite{PPR}.
\medskip
\par We develop an axiomatic framework (see Definition \ref{AA}) for which we recover the main result in \cite{PPP} (see Theorem \ref{thm:def}).
\newline
{\bf Theorem} Let $(\U,\partial)\models T_{\partial}$ and assume that $\U$ is sufficiently saturated. Let $\G$ be an $\cL_{\partial}$-definable group in $(\U,\partial)$ in the field sort. Assume that $\G$ is finite-dimensional (as an $\cL_{\partial}$-definable set).
Then there is an $\cL$-definable group $G$ (possibly adding finitely parameters) and an $\cL_{\partial}$-definable group embedding of $\nabla_{N}(\G)$, for some natural number $N$, into $G$ with the property that any $\cL$-generic type in $G$ is realized by a tuple $\nabla_{N}(a)$ with $a\in \Gamma$.
\bigskip
\par Moreover when $T$ is in addition dp-minimal, not strongly minimal or an open theory of topological fields containing the theory of a complete nondiscrete valued field of rank $1$, we show that there is an $\cL$-definable $D$-group such that $\G$ is isomorphic to the sharp points of that $D$-group (see Definition \ref{D-groups}), following \cite[section 1, Theorem]{PPP2}. 
 \par The plan of the paper is the following. First we will assume that $T$ is a complete geometric theory of topological fields of characteristic $0$ which is in addition either dp-minimal, or an open theory of topological fields containing the theory of a complete nondiscrete valued field of rank $1$ and we will prove differentiability properties of definable functions in models of $T$.
 \par Then we will set up an axiomatic setting under which one can associate with a definable group in a generic differential expansion of a model of $T$, a definable group in models of $T$. 
 
For sake of completeness, we will go over the pre-group construction; one reason for doing so, is that we recover a definable group as in \cite{PPP} and not simply an interpretable one in the finite-dimensional case. 

Then we will follow the same strategy as in \cite{PPP} in order to get from an $\cL_{\partial}$-definable finite-dimensional group $\G$, an $\cL$-definable pre-group and then from that pre-group, a definable group into which $\G$ {\it densely} embeds. Finally we will adapt the $D$-group construction introduced in \cite{PPP2}, which generalized the classical algebraic $D$-groups construction due to Buium, to this abstract setting.
\medskip
\par {\bf Acknowledgments} I would like to thank Kobi Peterzil and Anand Pillay for their collaboration and numerous comments on a previous version. I also would like to thank an anonymous referee of our joint paper  \cite{PPP}.
\section{Preliminaries}
\par We fix a complete $\cL$-theory $T$, possibly many-sorted, expanding the theory of fields of characteristic $0$ (the language $\cL$ containing the language of rings $\cL_{rings}=\{+, -, \cdot, 0, 1\}$). We work in a sufficiently saturated model $\U$ of $T$ and $A, B, C$ will denote small sets of parameters. Since ultimately, we will be interested in differential expansions of $T$, we will concentrate on the field sort $K$
and we will make the following technical assumption: we assume that for every $K$-valued $\cL$-term $t(x,z)$ with $x$ a tuple of $K$-variables and $z$ a tuple of variables varying in the other sorts, there is an $K$-valued $\cL$-term $\tilde t(x)$ such that $T\models \forall z\forall x (t(x,z)=\tilde t(x))$ \cite{K-P}. 
\par When $T$ is in addition a geometric theory on the field sort  \cite{Hrushovski-Pillay-groupslocalfields}, the model-theoretic algebraic closure $\acl$ (on the field sort) satisfies the exchange property (1)  and the quantifier $\exists^{\infty}$ is eliminated (2)). Note that since we are in a field, the first property (1) implies that $T$ eliminates the quantifier $\exists^{\infty}$ (on the field sort) \cite[Lemma 3.47]{F} (there, the result is attributed to A. Dolich, C. Miller and C. Steinhorn \cite[1.18]{DMS}), see also \cite[Theorem 2.5]{JY}.
\par We first recall a few properties of the dimension function induced by $\acl$ on the field sort in models of $T$ \cite[Remarks 2.2]{Hrushovski-Pillay-groupslocalfields} (see \cite[section 2]{Hrushovski-Pillay-groupslocalfields}, \cite[section 2]{JY}).
\subsection{Geometric theories}
\par Let $a$ be a finite  tuple of elements in the field sort and let $C$ be a subset in the field sort. Then $\dim(a/C)$ is the cardinality of a maximal subtuple of $a$ which is $\acl$-independent over $C$. This well-defined since $\acl$ has the exchange property; we will use the following property:
\smallskip
\par Let $a, b$ be two finite  tuples in the field sort, then
\begin{equation}\label{dim}
\dim(a b/C)=\dim(a/C b)+\dim( b/C).
\end{equation}
If $\dim(a/Cb)=\dim(a/C)$, we say that $a$ is independent from $b$ over $C$.
\par One extends the dimension function $\dim$ on partial types $\pi(x)$ over $C$ in finite tuple $x$ of variables (varying in the field sort), closed under finite conjunctions, as follows:
\[
\dim_C(\pi(x)):=\max\{\dim(a/C)\colon a {\rm\;\; realizing}\;\; \pi\}
\]
If $D$ is a small set of parameters in the field sort containing $C$, then $\dim_C(\pi(x))=\dim_D(\pi(x)$ \cite[Lemma 2.2]{Hrushovski-Pillay-groupslocalfields}, so we will simply write $\dim(\pi(x))$.
\par Then, $\dim(\pi) = \min\{\dim(\varphi(x)): \varphi\in \pi\}$. When $X$ is the solution set of an $\cL(C)$-formula $\varphi(x)$ where $x$ varying in the field sort), then we define $\dim(X)$ as $\dim(\varphi)$.
\par When $p$ is a complete type over $C$ in finite tuple $x$ of variables (varying in the field sort), then $\dim(p) = \dim(a/C)$ for $a$ realizing $p$.  
\par A tuple $a$ in an $\cL(C)$-definable set $X$ is called {\it generic over $b$} if $\dim(X)=\dim(a/C b)$, where $b$ is some tuple of elements in the field sort (when $b$ is empty we simply say {\it generic}). (Sometimes, we will use the term $\cL(C b)$-generic, when we want to stress over which language we are working on.) If $p$ is the complete type (over $C$) of an element $a$ of $X$, we say that $p$ is generic if $a$ is generic.

 \subsection{dp-minimal fields}
\par When $T$ is dp-minimal and not strongly minimal, then W. Johnson constructed a Hausdorff definable $V$-topology on the field sort \cite[Theorem 1.3]{WJ} and P. Simon and E. Walsberg proved  a cell decomposition theorem 
 \cite[Proposition 4.1]{simon-walsberg2016}, where the building blocks of definable sets are correspondences (see Definition \ref{corres}).
\par Note that there are dp-minimal structures where the model-theoretic algebraic closure, does not have the exchange property \cite[Proposition 5.2]{simon-walsberg2016}. However, W. Johnson showed that in C-minimal fields, it does have the exchange property \cite{WJ2}.
\par When $T$ is in addition a geometric theory,  one can show that the correspondences can be chosen to be closed under projections onto an initial set of coordinates (see \cite[Corollary 1.5.5]{K-P}). 
\par For convenience of the reader, we will recall the definition of first uniform structure and then of correspondences, following  \cite{simon-walsberg2016}.
\par  Indeed, when dealing with correspondences (since they are multi-valued functions), it is convenient to put on $K$ a uniform structure, which is another way to describe a topology on $K$.  
$U\circ U:=\{(x,z): \exists y\;(x,y)\in U\& (y,z)\in U\}$. Then the elements of $\cB$ are called basic entourages if:
\begin{enumerate}
\item $\bigcap_{U\in \cB} U=\Delta$,
\item for all $U, V\in \cB$, there is $W\in \cB$ such that $W\subset U\cap V$,
\item for all $W\in \cB$, there is $U\in \cB$ such that $U\circ U\subset W$. 
\end{enumerate}
\par The uniform structure generated by $\cB$ in $K$ is the set of subsets of $K^2$ containing an element of $\cB$. The uniform structure is said to be definable if there is a definable set $D$ and a formula $\varphi(x,y,\bar z)$ such that $\cB=\{\varphi(K^2,\bar d)\colon \bar d\in D\}$.
\par This induces a topology on $K$ by defining a neighbourhood basis of $x\in K$ as the set of $W[x]:=\{y\in K: (x,y)\in W\}$, with $W\in \cB$.
\par If one starts with a field $K$ endowed with a definable topology and define $\cB$ as the collection of symmetric open subsets of $K^2$  containing $\Delta$ and satisfying $(1)$ up $(3)$ above, then the topology on $K$ induced by this uniform structure coincides with  the original topology on $K$.
\par On cartesian products $K^d$ of $K$, one puts the product topology, namely the product uniform structure generated by $\{U^d\colon U\in \cB\}$. By abuse of notation we will still denote the product uniform structure by $\cB$.
\medskip
\dfn \cite[Section 3.1]{simon-walsberg2016} \label{corres} A definable {\it correspondence} $f\colon E\rightrightarrows K^\ell$ consists of a definable set $E$ (with $E\subseteq K^n$) together with a definable subset $\mathrm{graph}(f)$ of $E\times K^\ell$ such that
\[
0< \vert\{y\in K^\ell : (x,y)\in \mathrm{graph}(f)\}\vert<\infty, \text{ for all } x\in E.
\]
The set $\{y\in K^\ell : (x,y)\in \mathrm{graph}(f)\}$ is also denoted by $f(x)$. For a positive integer $m$, we say $f$ is an \emph{$m$-correspondence} if $\vert f(x)\vert=m$ for all $x\in E$ (if $m=1$, $f$ is a function). By convention, if $n=0$, $\mathrm{graph}(f)$ is identified with a finite set and if $E$ is an open subset of $K^n$ and $\ell=0$, $\mathrm{graph}(f)$ is identified with an open set.
\par The correspondence $f$ is {\it continuous} at $x\in E$ if for every open set $V\in \cB$ there exists $U\in \cB$ such that for any $(x,x')\in U$, we have $(f(x),f(x'))\in V$.
\par We will say that $f$ is $\cL(A)$-definable if the above data is defined over $A$.
\edfn
 \par We use the following conventions. Let $X\subseteq K^n$,  $\rho\colon K^n\to K^d$ be a coordinate projection and $f\colon \rho(X) \rightrightarrows K^{n-d}$ be a correspondence. Let $\rho^\perp\colon K^n \to K^{n-d}$ be the complement projection of $\rho$ (i.e. to the complement set of coordinates). The \emph{graph of $f$ along $\rho$} is the set
\[
\{x\in K^n : \rho(x)\in \rho(X) \text{ and } \rho^\perp(x)\in f(\rho(x))\}.
\]
\par Simon and Walsberg showed that $\dim$, the dp-rank and the topological dimension of a definable set all coincide \cite[Proposition 2.4]{simon-walsberg2016}, without assuming that $\acl$ has the exchange property. 
\medskip
\dfn \label{dlarge} 
Let $Y\subset X\subset K^n$ be two $\cL(C)$-definable sets in $\U$, then $Y$ is large in $X$ if $\dim(X\setminus Y)<\dim(X)$.
\edfn
\par A property is said to hold {\it almost everywhere} on $X$, if it holds on a large definable subset of $X$.

\fct \label{continuous}\cite[Proposition 3.7]{simon-walsberg2016} Assume that $T$ is dp-minimal and not strongly minimal. Then, on the field sort, every definable correspondence $f$ defined on an open set 
 is continuous on an open large definable subset. 
 \efct
 Under the same hypothesis on $T$, we have the following further property of correspondences.
 \fct \label{function} \cite[Lemma 3.1]{simon-walsberg2016} Let $f:U\rightrightarrows K^{\ell}$ be a continuous definable $m$-correspondence defined on an open definable set $U$. Then every $a\in U$ has an open neighbourhood $V$ such that there are $m$ continuous definable functions $g_i:V\to K^{\ell}$,  $1\leq i\leq m$, with $\mathrm{graph}(f\restriction V)=\bigsqcup_{i=1}^m \mathrm{graph}(g_i)$. 
\efct

 \par In \cite[Proposition 4.1]{simon-walsberg2016}, Simon and Walsberg showed that any definable subset of a dp-minimal topological structure (satisfying two additional properties which holds in case the structure is a field) is a finite union of graphs of definable continuous $m$-correspondences defined on open subsets. In \cite{K-P}, we revisited their result for open theories of topological fields \cite[Definition 1.2.4]{K-P} (note that these are geometric theories) and to stress that we get the additional property that the correspondences are closed under coordinate projections on an initial set of coordinates, we called them {\it cells}. 
\par Cells are defined by induction on $n$ and from the definition it will follow that if $X\subset K^n$ is a cell and $\pi\colon K^n\to K^d$ is a coordinate projection onto the first $d$-variables, then $\pi(X)$ is also a cell.
\medskip

\dfn\label{def:cells} Let $A\subseteq K$ and $X$ a non-empty $A$-definable subset of $K^n$, then $(X,\rho_X)$ is an $A$-cell if  $\rho_X\colon K^n \to K^{\dim(X)}$ is a coordinate projection such that $\rho_X(X)$ is an open set (with the convention that if $\dim(X)=0$, then $\rho_X(X)$ is open) and $X$ is the graph of a continuous $\cL(A)$-definable $m$-correspondence $f:\rho_X(X)\rightrightarrows K^{n-d}$, along $\rho_X$ with the following additional properties:

\begin{enumerate}[leftmargin=*]
\item if $n=1$, namely $X\subseteq K$ is an $A$-definable cell if and only if $X$ is an $A$-definable non-empty open set and $\rho_X$ is the identity, or $X$ is an $A$-definable non-empty finite set and $\rho_X$ is the projection to $K^0$. 

\item Assume $A$-definable cells have been defined for all $k\leqslant n$ and let $\pi\colon K^{n+1}\to K^n$ be the projection onto the first $n$-coordinates. A pair $(X,\rho_X)$ with $X\subseteq K^{n+1}$ is an $A$-definable cell if and only if $X$ is the graph along $\rho_X$ of a continuous $A$-definable $m$-correspondence $f_X\colon \rho_X(X)\rightrightarrows K^{n+1-\dim(X)}$ (for some $m>0$), there is an $A$-definable cell $(C, \rho_C)$ such that $C=\pi(X)$ and one of the following three cases holds
\begin{enumerate}[resume]
\item $\rho_X$ is the identity (so $X$ is an $A$-definable open set); 
\item $\rho_X=\rho_C\circ\pi$; 
\item $\rho_X$ corresponds to 
\[
\rho_X(x_1,\ldots, x_{n+1}) = (\rho_C(x_1,\ldots, x_n), x_{n+1}),  
\]
the fiber $X_c$ is a non-empty open subset of $K$ for every $c\in C$, and for every $x\in X$
\[
f_X(\rho_X(x)) = f_C(\rho_C(\pi(x))). 
\]
\end{enumerate}
\end{enumerate}
By a cell, we mean a $K$-definable cell. We often omit the associated projection $\rho_X$ of a cell, and simply write $X$ for a cell. 
\edfn
\fct  \cite[Corollary 1.5.5]{K-P}\label{cd} Let $T$ be geometric, dp-minimal and not strongly minimal. Let $\mathcal U\models T$, let $Y$ be an $\cL(A)$-definable subset of the field sort $K$. Then $Y$ is a finite disjoint union of $A$-cells.
\efct 

\par In the next section we will improve the above cell decomposition result, using the existence of  partial derivatives for continuous $m$-correspondences, almost everywhere.
\subsection{Partial derivatives of correspondences}\label{cont}
\par Let $(K,\tau)$ be a topological field endowed with a $V$-topology, equivalently it is induced by either an ordering or a (nontrivial) valuation. Furthermore if $K$ is $\omega$-complete (namely the filter of neighbourhoods of $0$ is closed under countable intersection), we may always assume that $\tau$ is induced by a nontrivial valuation \cite[Theorem 3.1]{PZ}. 
\par We use the multiplicative notation for the valuation, denoted by $\vert\;\vert$ (the valuation ring of $K$ is the set of elements $x$ of $K$ with $\vert x\vert\leq 1$ and we extend the value group of $K$ by: $\vert 0\vert=0$). For $u\in K^n$, we denote by $\Vert u\Vert:=\max_{i}\vert u_{i}\vert$, the induced norm on cartesian products of $K$. Let $K^*:=K\setminus\{0\}$.

Let us define the notions of partial derivatives and the property of being $C^1$ for correspondences. 
\dfn \label{C1-corr} Let $U\subset K^n$ be an open set and let $f: U\rightrightarrows K$, be an $m$-correspondence. 
\par If for $a:=(a_{1},\ldots,a_{n})\in U$,  $\lim_{\substack{\epsilon\rightarrow 0\\ \epsilon\in K^*}}\frac{f(a_{1}+\epsilon,a_{2},\ldots,a_{n})-f(a_{1},\ldots,a_{n})}{\epsilon}$ exists, then we denote it by $\partial_{x_{1}} f(a)$. It $\partial_{x_{1}} f(a)$ exists for every $a\in U$, then we get an $m$-correspondence $\partial_{x_{1}} f: U\rightrightarrows K$. We define similarly $\partial_{x_{i}} f$, $1\leq i\leq n$. When $f$ is definable, so are its partial derivatives.
\medskip
\par Then we say that $f$ is $C^1$ at $u\in U$ if there are $m$ linear maps $T(u):=(T_1(u),\ldots,T_{m}(u))$, $T_{i}(u):K^n\to K$ a linear function, $1\leq i\leq m$
such that $\lim_{\substack{h\rightarrow 0\\ h\in K^n}}\frac{\Vert f(u+h)-f(u)-T(u) h\Vert}{\Vert h\Vert}=0$. 
The $m$-correspondence $f$ is $C^1$ on $U$ if it is $C^1$ at every $u\in O$ and if the $m$-correspondence $df:U\rightrightarrows K^{n}:u\mapsto T(u)$ is continuous. 
\par Then by induction one defines: $f$ is $C^k$ on $U$, $k\geq 1$, if $f$ is $C^1$ on $U$ and $df$ is $C^{k-1}$ on $U$
It is easy to see if $f$ is $C^1$ on $U$, then all partial derivatives exist and are continuous on $U$.
\edfn
\par Note that the property of being $C^1$ for a definable $m$-correspondence $f$ on a definable neighbourhood $U$, is a first-order property: we express that 
\[\forall u\in U\bigwedge_{i=1}^m \exists \alpha_{1}^i\ldots\exists \alpha_{n}^i\; \lim_{\substack{h\rightarrow 0\\ h\in K^n}}\frac{\Vert f(u+h)-f(u)-(\sum_{j=1}^n \alpha_{j}^1 h_{j},\ldots,\sum_{j=1}^n \alpha_{j}^1 h_{j})\Vert}{\Vert h\Vert}=0.
\]
\par Finally for correspondences $f$ with values in $K^\ell$, $f:U\rightrightarrows K^\ell$, $f=(f_{1},\ldots,f_{\ell})$, with $f_{j}:U\rightrightarrows K$, one says that $f$ is $C^1$ if each $f_{j}$ is $C^1$, $1\leq j\leq \ell$. Similarly we say that $\partial_{x_{i}} f:U\rightrightarrows K^{\ell}$ exists if each $\partial_{x_{i}} f_{j}$, $1\leq j\leq \ell$ exist, $1\leq i\leq n$.

\
 \medskip
 \par Recall that if $A$ is a definable open subset of $K^n$, where either $\U$ is a model of an open theory of topological fields or a dp-minimal topological field and $B$ is a definable subset of $A$ dense in $A$, then the interior $Int(B)$ of $B$ is dense in $A$ (and $B$ is large subset of $A$) \cite[Lemma 2.6]{simon-walsberg2016}, \cite[Lemma 1.4.4]{K-P}.
  \medskip
  \nota \par Let $V$ be an open subset of $K^n$ and let 
 $f: V\to K$ be a definable continuous function. For $h\in K^*$, we denote by $f_{h}(x)$, for $x:=(x_{1},\ldots,x_{n})\in V$, $f_{h}(x):=f(x_1+h,x_2,\ldots,x_n)-f(x_1,x_2,\ldots,x_n)$.
 \enota
 \par Following the same strategy as in \cite[Proposition 3.7]{simon-walsberg2016}, we get:
 \prop \label{partial} Let $T$ be dp-minimal and not strongly minimal. Let $V$ be an open definable subset of $K^n$ and let $f: V\rightrightarrows K$ be a definable continuous $m$-correspondence. Then for each $1\leq i\leq n$, $\partial_{x_{i}} f$ exists and is continuous almost everywhere.
 \eprop
 \pr Fix $i=1$ and let us show that $\partial_{x_{1}} f$ exists almost everywhere, namely that $\lim_{h\mapsto 0} \frac{f_{h}(x)}{h}$ exists for almost all $x\in V$. 
 By \cite[Lemma 2.6]{simon-walsberg2016}, it suffices to show that the sets of points where it exists is dense in $V$. So fixing an open subset $V'$ of $V$, it suffices to show that this limit exists at some point of $V'$. By the way of contradiction, assume this is not the case and replacing $V'$ by $V$, we suppose that for no points in $V$, this limit exists. 
 So for all $p:=(x_{1},\ldots,x_{n})\in V$ there exists $W\in \cB$ such that $(W,p)$ has the following property. For all $W'\in \cB$, there is $(x_{1}+h,x_{1}+h')\in W'$, with $h, h'\neq 0$, such that $(\frac{f_{h}(p)}{h},\frac{f_{h}(p')}{h})\notin W$. We consider the collection $B$ of all such $(W,p)$ and
 of the subsets of $V$ of the form $B_W:=\{p\in V: (W,p)\in B\}$; we have $V=\bigcup_W B_W$ and this set is directed: $B_{W_1}\cup B_{W_2}\subset B_{W_{1}\cap W_{2}}$. So by \cite[Lemma 3.5]{simon-walsberg2016}, we find $W_0\in \cB$ such that $B_{W_0}$ has non empty interior. By replacing $V$ by a smaller open subset we may assume that $V\subset B_{W_0}.$ So for any $p=(x_{1},\ldots, x_{n})\in V$, there are 
$h, h'\neq 0$ such that $(x_{1}+h,x_{1}+h')$ is as close as we wish of $(x_{1},x_{1})$ but $(\frac{f_{h}(p)}{h},\frac{f_{h'}(p)}{h'})\notin W_{0}$ $(\star)$.  Let $U\in \cB$ be such that $U\circ U\subset W_0$.
\par Let $D:=\{(p,h,\frac{f_h(p)}{h}): p\in V\}$; since $D$ is definable and infinite, $D$ is not discrete by \cite[Lemma 3.6]{simon-walsberg2016} (moreover we may assume that the projection of $D$ on $V$ is infinite and so the set of corresponding images by $\frac{f_h}{h}$ is infinite too. Let $(p_0, h_0, \frac{f_{h_0}(p_0)}{h_0})$ be an accumulation point of $D$.
Let $O$ be a neighbourhood of $(p_0, h_0, \frac{f_{h_0}(p_0)}{h_0})$ of the form $O_1\times O_2\times U[ \frac{f_{h_0}(p_0)}{h_0}]$, with $O_1$ an open neighbourhood of $p_0$, $O_2$ an open neighbourhood of $h_0$, so $O\cap D$ is infinite. So its projection $\pi$ onto either $O_1$ or $O_2$ is infinite and so has non-empty interior in $V$.

\par Since $p_0\in V$, there are $p_1, p_2\in V\cap \pi(O)$ such that $(\frac{f_{h_1}(p_1)}{h_1}, \frac{f_{h_2}(p_2)}{h_2})\notin W_0$. But $(\frac{f_{h_1}(p_1)}{h_1}, \frac{f_{h_0}(p_0)}{h_0})\in U$ and $(\frac{f_{h_2}(p_2)}{h_2},\frac{ f_{h_0}(p_0)}{h_0})\in U$. So by definition of $U$, $(\frac{f_{h_1}(p_1)}{h_1}, \frac{f_{h_2}(p_2)}{h_2})\in W_0$, a contradiction. 
\par The last part of the statement follows from \cite[Proposition 3.7]{simon-walsberg2016} (see Fact \ref{continuous}). \qed
 \nota \par Let $V$ be an open subset of $K^n$ and let 
 $f: V\to K$ be a definable continuous function and suppose that  for each $1\leq i\leq n$, $\partial_{x_{i}} f$ exists and is continuous on $V$. For $h=(h_1,\ldots,h_n)$ with $h_i\in K^*$, $1\leq i\leq n$, we denote by $f_{h}(x)$, for $x:=(x_{1},\ldots,x_{n})\in V$, $f_{h}(x):=f(x_1+h_1,x_2+h_2,\ldots,x_n+h_n)-f(x_1,x_2,\ldots,x_n)$, $T_h f(x):=\sum_{i=1}^n \partial_i f(x) h_i$, $D_h f(x):=\frac{\Vert f_h(x)-T_h f(x)\Vert}{\Vert h\Vert}$.
 \enota

 \prop \label{C1} Let $T$ be geometric, dp-minimal and not strongly minimal. Let $V$ be an open definable subset of $K^n$ and let $f: V\rightrightarrows K$ be a definable $m$-correspondence. Then for all $x\in V$, $df(x)$ exists and is continuous, almost everywhere on $V$.
 \eprop
\pr By assuming that $K$ is $\aleph_1$-saturated, we may also assume that the topology on $K$ is given by a non-archimedean absolute value $\vert\cdot\vert$. We first apply Fact \ref{continuous} and then Proposition \ref{partial} and so we may assume, by replacing $V$ by a dense open subset, that for each $1\leq i\leq n$, $\partial_{x_{i}} f$ exists and is continuous on $V$.
\par Let $x\in V$ and suppose that $x$ is generic. Let $\eta\in \vert K\vert$ and let $U\subset K^n$ be a neighbourhood of $0$ such that $x+U\subset V$.
Consider $U_\eta(x):=\{h\in U\colon  \Vert f_h(x)-T_h f(x)\Vert<\eta\Vert h\Vert\}$.
\cl For all $x\in V$, $\dim(U_\eta(x))=n$.
\ecl
\prcl 
For sake of simplicity assume that $n=2$ and let $h=(h_1,h_2)\in K^2$. Express 
\begin{equation}\label{eq_tri}
f_h(x)-T_h f(x)=
\end{equation}
\begin{align*}
&f(x_1+h_1,x_2+h_2)-f(x_1,x_2+h_2)-\partial_1 f(x_1,x_2+h_2) h_1+\\
&-f(x_1,x_2)+f(x_1,x_2+h_2)-\partial_2 f(x_1,x_2) h_2+\\
&\partial_1   f(x_1,x_2+h_2) h_1-\partial_1 f(x_1,x_2) h_1.
\end{align*}
Let $O\subset \vert K\vert$ be a neighbourhood of $0$ of the form $\{u\in \vert K\vert: u<\eta\}$ with $\eta\in \vert K\vert$.
\par Since $\vert f_h(x)-T_h f(x)\vert=\max\{\vert f(x_1+h_1,x_2+h_2)-f(x_1,x_2+h_2)-\partial_1 f(x_1,x_2+h_2) h_1\vert,\vert -f(x_1,x_2)+f(x_1,x_2+h_2)-\partial_2 f(x_1,x_2) h_2\vert, \vert \partial_1   f(x_1,x_2+h_2) h_1-\partial_1 f(x_1,x_2) h_1\vert\}$.

\par Since $\partial_2 f(x_1,x_2)$ exists for all $\eta>0$ there exists $U_0$ a neighbourhood of $0$ such that for any $h_2\in U_0\setminus\{0\}$, $\Vert \frac{f(x_1,x_2)-f(x_1,x_2+h_2)}{h_2}-\partial_2 f(x_1,x_2) \Vert<\eta$.
\par Since $\partial_1 f(x_1,x_2)$ is continuous in a neighbourhood of $(x_1,x_2)$, for all $\eta>0$ there exists $U_1\subset U_0$ a neighbourhood of $0$ such that for any $h_2\in U_0$, $\Vert \partial_1   f(x_1,x_2+h_2) -\partial_1 f(x_1,x_2) \Vert<\eta$.
\par Let $U_2\subset K$ be the set of $h_2$ such that $\vert \partial_1   f(x_1,x_2+h_2)-\partial_1 f(x_1,x_2)\vert<\eta$ and $\vert -f(x_1,x_2)+f(x_1,x_2+h_2)-\partial_2 f(x_1,x_2) h_2\vert< \eta\vert h_2\vert$. Then $U_2$ is an open neighbourhood of $0$. W.l.o.g. we may assume that $\pi_2(U)=U_2$. 
\par Then fixing $h_2\in U_2$, we let $U_2(h_2):=\{h_1\in K:\;\vert f(x_1+h_1,x_2+h_2)-f(x_1,x_2+h_2)-\partial_1 f(x_1,x_2+h_2)h_1 \vert\leq \eta\vert h_1\vert$. Note that $U_2(h_2)$ is an open neighbourhood of $0$. Indeed, since 
$(x_1,x_2+h_2)\in V$, $\partial_1 f(x_1, x_2+h_2)$ exists.
\par Then for $h_2\in U_2$ and $h_1\in U_2(h_2)$, we have that $\vert f_h(x)-T_h f(x)\vert<\eta \max\{\vert h_1\vert, \vert h_2\vert\}$. 
Define $\tilde U_\eta(x)\subset K^2$ as the set of $\{(h_1,h_2): h_2\in U_2\wedge h_1\in U_2(h_2)\}$. So for each $h\in \tilde U_\eta(x)$, $D_h f< \eta$.
\par By additivity of the dimension (we assume that $\acl$ had the exchange), we have that $\dim(\tilde U_\eta(x))=2$. So by \cite[Lemma 2.2]{simon-walsberg2016}, $\tilde U_\eta(x)$ has non-empty interior and so the same holds for $U_\eta(x)$. 
\qed
\bigskip
\par Now consider $U_2(0):=\{h_1\in K:\;\vert f(x_1+h_1,x_2)-f(x_1,x_2)-\partial_1 f(x_1,x_2)h_1 \vert\leq \eta\vert h_1\vert$. As we saw in the Claim above, $U_2(0)$ is an open neighbourhood of $0$ in $K$.
Consider $\{y_2\in \pi_2(V): \forall h_1\in U_2(0)\;\vert f(x_1+h_1,x_2+h_2)-f(x_1,x_2+h_2)-\partial_1 f(x_1,x_2+h_2)h_1 \vert\leq \eta\vert h_1\vert\}.$ This subset contains $x_2$ and since $x$ is generic in $V$ (and we may assume that $x_2$ is generic over $x_1$), it contains an open subset $U(x_2)$ containing $x_2$. 
\par Now by equation (\ref{eq_tri}), we have that for $(h_1,h_2)\in U_2(0)\times U_1$, $\Vert f_h(x)-T_h f(x)\Vert<\eta\Vert h\Vert.$ So $df(x)$ exists and it is continuous almost everywhere by Fact \ref{continuous}.
\qed
\medskip
\cor \label{partialopen} Let $T$ be an open theory of topological fields. Let $V$ be an open definable subset of $K^n$ and let $f: V\rightrightarrows K$ be a definable continuous $m$-correspondence. Then for each $1\leq i\leq n$, $\partial_{x_{i}} f$ exists almost everywhere.
 \ecor
\pr We use  \cite[Lemma 1.4.4]{K-P} instead of \cite[Lemma 2.6]{simon-walsberg2016},  \cite[Lemma 1.4.5]{K-P} instead of \cite[Lemma 3.5]{simon-walsberg2016} and \cite[Lemma 1.4.7]{K-P} instead of \cite[Lemma 3.6]{simon-walsberg2016}.
\qed
\subsection{Implicit function theorems} \label{sec:implicit}
\par Let $(K,\tau)$ be a topological field; recall that the field $(K,\tau)$ is $t$-henselian if for every $n\geq 1$, there exists $U\in \tau$ such for every polynomial $p(x)\in K[x]$ of the form $x^{n+1}+x^{n}+U[x]^{n-1}$ has a zero in $K$, where $U[x]^{n-1}$ is the set of polynomials of degree $\leq n-1$ with coefficients in $U$ \cite[Theorem 7.2]{PZ}.
\par A. Prestel and M. Ziegler showed that a $V$-topological field $(K,\tau)$ is $t$-henselian if and only if it satisfies the following implicit function theorem  \cite[Theorem 7.4]{PZ}, namely for any $p(x_0,\ldots,x_{n},y) \in K[x_0,\ldots,x_{n},y]$, $a_{0},\ldots,a_{n},b\in K$ such that 
$p(a_{0},\ldots,a_{n},b)=0$ and $\partial_{y} p(a_{0},\ldots,a_{n},b)\neq 0$, there exist $U, V\in \tau$ such that for any $a_{0}'\in a_{0}+U,\ldots, a_{n}'\in a_{n}+U$,
there exists a unique $b'\in b+V$ such that $p(a_{0}',\ldots,a_{n}',b')=0$ and the map from $a_{0}+U\times \ldots\times a_{n}+U$ to $K$ sending $(a_{0}',\ldots, a_{n}')$ to $b'$ is continuous.
\par We will need other versions of the implicit function theorem where we may infer that the implicitly defined functions are differentiable. 
\nota \label{Jac} Let $U\subset K^{\ell}, V\subset K^m$ be open sets, let $\vert x\vert=\ell$, $\vert y\vert=m$. Let $f(x,y):=(f_1(x,y),\ldots,f_m(x,y))$ be a tuple of functions from $U\times V\to K^m$.

Denote by $f_{a}(y):=(f_{1}(a,y),\ldots, f_{m}(a,y))$. Assume that for every $y\in V$,  the partial derivative with respect to $y_j$, $\partial_{j} f_i(a,y)$ exists, $1\leq i, j\leq m$. Then we denote by $J_{f_a}(y)$ the Jacobian matrix: $\left( \begin{array}{ccc}
     \partial_{1} f_1(a,y) & \cdots & \partial_{m} f_1(a,y) \\
      \vdots & \ddots & \vdots\\
      \partial_{1} f_m(a,y) & \cdots & \partial_{m} f_m(a,y)  
      \end{array} \right).$
      
\par As usual, we denote by $det(J_{f_a}(y))$ the determinant of the matrix $J_{f_a}(y)$. 
\enota

\

\par Let us state the following version of the implicit function theorem, which holds in real-closed fields \cite[Section4]{W}.
\dfn\label{imp}  Let $n=\ell+m$, $n>1$, $\ell, m>0$, let $(a,b)\in K^{\ell+m}$.
 Let $f_1(x,y),\ldots, f_m(x,y)$ be definable $C^{1}$-functions in $K$, $\vert x\vert=\ell$, $\vert y\vert=m.$
Then $K$ satisfies $(\IFTA)$ if the following holds.
Assume that $\bar f_{a}(b)=0$ and that $det(J_{f_{a}}(b))\neq 0$, with $f:=(f_{1},\ldots, f_{m})$.
Then  there are  neighbourhoods $O_{a}\subset K^\ell$ of $a$ (respectively $O_{b_i}\subset K$, $1\leq i\leq m$, of $b_i$) and continuous functions $g_{i}(x): O_{a} \to O_{b_i}$, $1\leq i\leq m$, whose partial derivatives exist, such that , setting $g:=(g_1,\ldots, g_m)$ and $\bar \ell:=(1,\ldots,\ell)$,
\begin{align} \label{eq9} & g(a)=b\;\wedge\\
\label{eq10} &\forall\; x\in O_{a}\;\;\big(f(x, g(x))=0\; \wedge J_{g}(x)=-(J_{f_{x}}(g(x)))^{-1}J_{f_{g(x)}}(x)\big)\;\wedge\\
\label{eq11} & \forall\; x\in O_{a}\;\forall\; y\in O_{b}\;(f(x,y)=0 \leftrightarrow y=g(x)).
\end{align}
We say that $K$ satisfies $(\IFTA)_{\cP}$ when we only let $f_1(x,y),\ldots, f_m(x,y)$ to vary in the class of polynomial functions.
\edfn
\par When $K$ is a complete non discrete valued field of rank $1$ (namely its value group is a subgroup of $\R$), the following version of the implicit function theorem holds for power series $K[[x,y]]$ \cite[Chapter 2, paragraph 10]{A}; it relies on Weierstrass preparation theorem. In particular it shows that $K$ satisfies $(\IFTA)_{\cP}$.
\thm\label{imp_p} \cite[10.8]{A} Let $n=\ell+m$, $n>1$, $\ell, m>0$, let $(a,b)\in K^{\ell+m}$.
 Let $f_1(x,y),\ldots, f_m(x,y)\in K[[x,y]]$, $\vert x\vert=\ell$, $\vert y\vert=m$ such that $f_{a}(b)=0$ and that $det(J_{f_{a}}(b))\neq 0$.
Then there exist unique elements $\xi_{1}(x),\cdots,\xi_{m}(x)\in K[[x]]$ such that $\xi_{1}(a)=b_{1},\ldots,\xi_{m}(a)=b_{m}$ and $\bigwedge_{i=1}^m f_{i}(x,\xi_{1}(x),\cdots,\xi_{m}(x))=0$.
Assume that the $f_{i}$'s are convergent power series.
Then  there are  neighbourhoods $O_{a}\subset K^\ell$ of $a$ (respectively $O_{b_i}\subset K$, $1\leq i\leq m$, of $b_i$) such that $f_{1},\ldots,f_{m}$ are convergent on $O_{a}\times O_{b_{1}}\times\ldots O_{b_{m}}$, $\xi_{i}(x)$ is convergent on $O_{a}$, $1\leq i\leq m$, and for all $u\in O_{a}$,
the tuple $(\xi_{1}(u),\cdots,\xi_{m}(u))$ is the unique solution of $\bigwedge_{i=1}^m f_{i}(u,z_{1},\ldots,z_{m})=0$ $(\star)$.
\ethm
So in case the $f_{i}$'s are convergent, for instance if there are polynomials in $x, y$, then we may require that (in a certain neighbourhood) there are unique elements $z_{1},\ldots,z_{m}$ of $K$ such that $(\star)$ holds. So we get definable functions $\xi_{1}(x),\cdots,\xi_{m}(x)$ in the corresponding neighbourhoods. So we get the following corollary.

\cor\label{imp_p-def} Let $K$ be a topological field endowed with a definable topology and assume that $K$ is elementary equivalent to a complete nondiscrete valued field of rank $1$.  Let $n=\ell+m$, $n>1$, $\ell, m>0$, let $(a,b)\in K^{\ell+m}$.
 Let $f_1(x,y),\ldots, f_m(x,y)\in K[x,y]$, $\vert x\vert=\ell$, $\vert y\vert=m$ such that $f_{a}(b)=0$ and that $det(J_{\bar f_{a}}(b))\neq 0$ (see Notation \ref{Jac}). 
Then there are  neighbourhoods $O_{a}\subset K^\ell$ of $a$ (respectively $O_{b_i}\subset K$), $1\leq i\leq m$, of $b_i$, there exist $C^1$-definable functions $g_{1}(x),\cdots, g_{m}(x)$ on $O_{a}$ (with values in $K$) such that $g_{1}(a)=b_{1},\ldots,g_{m}(a)=b_{m}$ and $\bigwedge_{i=1}^m f_{i}(x,g_{1}(x),\cdots,g_{m}(x))=0$. Furthermore $(g_{1}(x),\cdots,g_{m}(x))$ is the unique solution of $\bigwedge_{i=1}^m f_{i}(x,u_1,\cdots,u_{m})=0$. 
\par Finally, 
 $\forall\; x\in O_{a}\;\;\big(J_{g}(x)=-(J_{f_{x}}(g(x)))^{-1}J_{\bar f_{g(x)}}(x)\big)$ and $\forall\; x\in O_{a}\;\forall\; y\in O_{b}\;(f(x,y)=0 \leftrightarrow y=g(x)).$
\ecor
\pr We use Theorem \ref{imp_p} for complete nondiscrete valued fields of rank $1$. We first assume that $K$ has this property.
So, there is an open neighbourhood $O_a$ of $a$ and an open neighbourhood $O_b$ of $b$ with $O_{a}\times O_{b}\subset \tilde O_{j}$, such that there is a unique power series solution $\xi(x)\in K[[x]]$ convergent on $O_a$ such that for any $x\in O_a$, $P_j(x,\xi(x))=0$. In particular $\xi$ is a $C^1$-map on $O_a$. there is an open neighbourhood $O_a$ of $a$ and an open neighbourhood $O_b$ of $b$ with $O_{a}\times O_{b}\subset \tilde O_{j}$, such that there is a unique power series solution $\xi(x)\in K[[x]]$ convergent on $O_a$ such that for any $x\in O_a$, $P_j(x,\xi(x))=0$. In particular $\xi$ is a $C^1$-map on $O_a$. These properties are first-order expressed and so it can be transferred to any field elementary equivalent to $K$.
\qed
\medskip
\par Recall that an ordered abelian group $G$ is elementarily equivalent to an archimedean group if and only if it is regularly ordered, namely for every natural number $n\geq 2$, every infinite interval of $G$ contain an element divisible by $n$ \cite[Theorem 4.7]{RZ}. In Appendix 3, we will show that any henselian valued field $(K,v)$ of characteristic $0$ with regularly ordered value group is elementary equivalent to a complete valued field of rank $1$, provided the residue characteristic is equal to $0$, or in the finitely ramified case in mixed characteristic. 
\medskip
\par Recall that the value group $\Gamma$ of a dp-minimal valued field of characteristic $0$ is non singular, namely $\Gamma/p\Gamma$ is finite, for every prime number $p$ \cite[Proposition 5.1]{JSW}. Being non-singular characterizes dp-minimal ordered groups among regularly ordered groups.
\section{Differential expansions}\label{sec:diff}
\par From now on $T$ will denote a complete geometric theory of fields of characteristic $0$.
\par Now we consider the expansion of models $\M$ of $T$ by a derivation $\partial$, acting on the field sort.
Let $\cL_\partial$ be the language $\cL$ extended by a unary function symbol $\partial$. Denote by $(\M,\partial)$ the expansion of $\M$ to an $\cL_\partial$-structure.
We will extend $T$ to an $\cL_{\partial}$-theory containing the usual axioms of a derivation (on the field sort), namely,  
\[
\begin{cases}
\forall x\forall y(\partial(x+y)=\partial(x)+\partial(y)),\\
\forall x\forall y(\partial(xy)=\partial(x)y+x\partial(y)).
\end{cases}
\] 
When on the field sort, $\cL$ contains additional function symbols besides the field ones,
 we will add a scheme of axioms expressing that $\partial$ respects $\cL$-definable functions on the field sort (see Definitions \ref{def:comp}, \ref{T++}) and the above axiom will express that $\partial$ respects the ring operations.
\medskip
\par We denote by $(\U,\partial)$ a sufficiently saturated model of $T(\partial)$ and as before denote by $K$ the field sort of $(\U,\partial)$. 
When we want to stress that we are in a differential expansion, we use $(K,\partial)$ to indicate that the field sort of $\U$ has been expanded by $\partial$.

\nota \label{nota:der}  
Let $X$ be a subset of $K$; we denote by $\langle X\rangle$ the differential subfield of $K$ generated by $X$.
For~$m\geqslant 1$ and~$a\in K$, we define 
\[
\partial^m(a):=\partial(\partial^{m-1}(a)), \text{ with $\partial^0(a):=a$,}
\]
and~$\bar{\partial}^m(a)$ as the finite sequence~$(\partial^0(a),\partial(a),\ldots,\partial^m(a))\in K^{m+1}$; we will also use the notation $\nabla_{m}(a)$ (and replace $\nabla_{1}(a)$ by $\nabla(a)$). 

Similarly, given an element~$a=(a_1,\ldots,a_n)\in K^n$ and $\bar m\in \N^n$, we will write~$\partial^{\bar m}(a)$ (or $\nabla_{\bar m}(a)$) to denote the element~$(\partial^{m_{1}}(a_1),\ldots,\partial^{m_{n}}(a_n))\in K^{n}$. When all the $m_{i}$ are equal, we denote $\partial^{\bar m}(a)$ by $\partial^m(a)$.
We will also use the notation $\nabla_{\bar m}(a)$ or $\nabla_m(a)$ when the tuple $\bar m$ is constant.
For $A\subseteq K^n$, let $\nabla_{\bar m}(A):=\{(\bar{\partial}^{\bar m}(a):\;a\in A\}$
and accordingly $\nabla_m(A)$ when the tuple $\bar m$ is constant. 

\par We also consider the infinite sequence $\nabla_{\infty}(a):=(a,\partial(a),\ldots)$, $a\in K^n$. 
\par The order of a tuple $x:=(x_{1},\ldots,x_{n})$ in an $\cL_{\partial}$-formula of the form $\theta((\bar{\partial}^{m_{1}}(x_1),\ldots,\bar{\partial}^{m_{n}}(x_n))$ where $\theta$ is an $\cL$-formula, is less than $\bar m$. When we don't want to specify the length of a differential tuple, we allow ourselves to use $\varphi(\nabla_\infty(a))$, where $\varphi$ is an $\cL$-formula.
\par We denote the length of a tuple $x$ by $\vert x\vert$.
\enota

\par Let us recall the notion of $\partial$ compatible with (or respecting) definable functions \cite{W05}, \cite{FK}.
  \dfn \label{def:comp} Let $A$ be a subset of the field sort and let $f: O\rightrightarrows K$ be a  continuous-$\cL(A)$-definable $m$-correspondence, where $O$ is an open $\cL(A)$-definable subset of $K^n$, with the property that for each $1\leq i\leq n$, the partial derivatives $\partial_{i} f$ with respect to $x_{i}$ exist (on $O$) and are continuous. So we get a continuous-$\cL(A)$-definable $m$-correspondence sending $u=(u_1,\ldots, u_n)\in O$ to $\sum_{i=1}^n \partial_{i} f(\bar u)\partial(u_{i})$.
  Then $\partial$ is compatible with $f$ if  we have:
\par  in case $A=\emptyset$, for each $u\in O$,
\[
\partial(f(u))=\sum_{i=1}^n \partial_{i} f(u)\partial(u_{i}),
\]
\par in case $A\neq\emptyset$, there is a $C^0$-$\cL(A)$-definable $m$-correspondence $f^{\partial}$ on $O$ such that for each $u\in O$, 
\[
\partial(f(u))=\sum_{i=1}^n \partial_{i} f(u)\partial(u_{i})+f^{\partial}(u),
\]
moreover if $\partial(A)=0$, then $f^{\partial}=0$.
\edfn
\par As usual we extend this definition when $f$ takes its values in some cartesian product of $K$, composing $f$ with projections. 
\medskip

\rem
Let $(u,v)\in K^n\times K^\ell$ and let $O=O_1\times O_2$ with $O_1$ (respectively $O_2$) an open definable subset of $K^n$ (respectively $K^\ell$). Let $F(u,v)$ be an $\cL$-definable continuous $m$-correspondence on $O$. Let $c\in O_2$ and consider the correspondence $f(u):=F(u, c)$ on $O_1$. Suppose that all partial derivatives of $F$ exist and are continuous on $O$ and that $\partial$ is compatible with $F$ on $O$. Then using Notation \ref{Jac}, $\partial(F(u,v))=J_{F}(u,v) (\partial(u),\partial(v))=\sum_{i=1}^n\partial_{i}F(u,v) \partial(u_{i})+\sum_{j=1}^\ell \partial_jF(u,v) \partial(v_{j})$. 
Set $f^{\partial}(u)=\sum_{j=1}^\ell \partial_jF(u, c) \partial(c_{j})$. Then $f^{\partial}$ is continuous on $O_1$ and $\partial(f(u))=\sum_{i=1}^n \partial_{i} f(u)\partial(u_{i})+f^{\partial}(u)$. So $\partial$ is compatible with $f$.
\par In particular $\partial$ is compatible with any polynomial map with coefficients in $K$. 
\erem
\dfn\label{T++} The expansion $T(\partial)$ of $T$ consists of $T$ together the following scheme of axioms: the derivation $\partial$ is compatible with every $\cL$-definable continuous $m$-correspondence $f$ on some definable open set $O\subset K^n$, with the property that for each $1\leq i\leq n$, the partial derivative $\partial_{i} f$ with respect to $x_{i}$ exists (on $O$) and is continuous. 
\edfn

\rem The set of $C^1$-correspondences compatible with $\partial$ is closed under composition by \cite[Lemma 2.6]{FK} and  by the implicit function theorem (it was already stated in \cite[Lemma 2.7]{FK} but there, one implicitly assumed that everything was defined over $\0$). For convenience of the reader we will show both properties in  Appendix 3. 
\erem
\lem Let $T$ be an open theory of topological fields containing the theory of a complete nondiscrete valued field of rank $1$. Let $\U\models T$ and let $\partial$ be a derivation on the field sort $K$ of $\U$. Let $f: K^n\rightrightarrows K$ be an $\cL(A)$-definable continuous $m$-correspondence on some definable open set $O\subset K^n$.
Then $\partial$ is compatible with $f$.
\elem
\pr For $\cA$ a finite subset of $K[x,y]\setminus\{0\}$, $S\in K[x,y]\setminus\{0\}$, where $\vert x\vert= n$, $\vert y\vert=1$, denote by $Z_{\cA}^{S}(x,y)$ the $\cL_{ring}$-formula $\bigwedge _{P\in \cA} P(x,y)=0\wedge S(x,y)\neq 0$. Set $\cA^y:=\{P\in \cA: deg_y(P)>0\}$. By \cite[Corollary 1.3.3]{K-P}, the graph of $f$, as an $\cL(A)$-definable subset of $K^n\times K$ can be expressed as
 a finite union of definable subsets of the following form:
\[
 Z_{\cA_\ell}^{S_\ell}(x,y) \wedge \theta_\ell(x,y)
\]
where $S_\ell\in K[x,y]$, $\cA_\ell$ a finite subset of $K[x,y]\setminus\{0\}$, 
$\theta_\ell$ an $\cL$-formula that defines an open subset of $K^{n+1}$,  $\cA_\ell^y=\{P_\ell\}$ and $\partial_{ y} P_\ell$ divides $S_\ell$.
\par Further we may assume that for $\ell\neq \ell'$, $Z_{\cA_\ell}^{S_{\ell}}\cap  Z_{\cA_\ell'}^{S_{\ell'}}=\emptyset$ (otherwise we perform Euclidean algorithm to decompose an intersection of the form $Z_{\cA_j}^{S_j}(x,y)\cap Z_{\cA_i}^{S_i}(x,y)$ as a finite union of sets of the same form.
Let $(a,b)\in graph(f)$ with $\vert a\vert=n$ and suppose that $P_\ell(a,b)=0$, $\partial_y P_\ell(a,b)\neq 0$ and $\theta_\ell(a,b)$ holds. Since $P_\ell(x,y)$ is a polynomial, $\partial$ is compatible with the polynomial function from $K^n\times K\to K: (a,b)\mapsto P_\ell(a,b)$ (which is a $C^1$-map). So we have $\partial(P_\ell(a,b))=\sum_{i=1}^n \partial_i P_\ell(a,b) \partial(a_i)+ \partial_y P_\ell(a,b) \partial(b)+P_\ell^{\partial}(a,b)$. Let $\pi$ be the projection on the first $n$ coordinates and let $O_\ell\subset \pi(\theta_\ell(K))\cap O$.
By the implicit function theorem (see Corollary \ref{imp_p-def}), for each $a\in O_\ell$, we can find a neighbourhood $V\subset O_\ell$ of $a$ such that there are $m$ $C^1$-functions $(g_1,\ldots, g_m)$ on $V$ such that $(a,g_1(a),\ldots,g_m(a))$ is the unique solution on $P_\ell(x,y)=0$ in $(V\times K)\cap \theta_\ell(K)$. So, 
$f\restriction V=(g_1,\ldots,g_m)$.
Replacing $b$ by one of $g_j(a)$ with $1\leq j\leq n$, we get $\sum_{i=1}^n \partial_i P_\ell(a,g_j(a)) \partial(a_i)+ \partial_y P(a,g_j(a)) \partial(g_j(a))+P_\ell^{\partial}(a,g_j(a))=0$. Then we use Lemma \ref{implicit} in the Appendix. 
\par We set $g_j^{\partial}:=-\frac{P_\ell^{\partial}(a,g_j(a))}{\partial_y P_\ell(a,g_j(a))}$ on $V$. (This is well-defined since $S_\ell(a, b)\neq 0$ holds and $\partial_{ y} P_\ell$ divides $S_\ell$.)
Then define for $a\in V$, the $m$-correspondence  that we denote by $f^{\partial}$,  by 
\[
f^{\partial}(a):=-\frac{P_\ell^{\partial}(a,f(a))}{\partial_y P_\ell(a,f(a))},
\]
it is continuous and $\cL(A)$-definable. 
\par  By equation (\ref{eq_implicit}) in Lemma \ref{implicit},  we get that $\partial$ is compatible with $f\restriction O_\ell$ (note that since $\partial_i f$, $1\leq i\leq n$, exist and are continuous, we get that $\partial_i f=(\partial_i g_j)_{j=1}^m$ on $V$. Since $\bigcup_\ell O_\ell=O$, we get the result.
\qed
\medskip
\par  Therefore, when $T$ is an open theory of topological fields, then $T(\partial)$ is simply $T$ together with the axiom stating that $\partial$ is a derivation. 
\medskip
\lem \label{dp-comp} Let $T$ be dp-minimal, not strongly minimal. Assume that $(\U,\partial)\models T(\partial)$. Let $k$ be a differential subfield of $K$ ($k$ sufficiently saturated) and let $f: K^n\rightrightarrows K$ be an $\cL(k)$-definable $C^1$ $m$-correspondence on some definable open set $O\subset K^n$
Then $\partial$ is compatible with $f$.
\elem
\pr Since $f(x)$ is an $\cL(k)$-definable $m$-correspondence on $O$, there is a $\cL$-formula $\varphi(x,y,z)$ such that for all $x\in O$, $y\in f(x)$ if and only if $\varphi(x,y,d) \wedge \forall x\in O \exists^{=m} y\,\varphi(x,y,d)$ holds, where $d\subset k$ and $\vert d\vert=\vert z\vert$.  Consider $Z:=\{z\in K^s: \forall x\in O \exists^{=m} y\,\varphi(x,y,z)\}$. By the cell decomposition theorem (see Fact \ref{cd}), $Z$ is a finite union of cells and we may assume that $d=(d_1,d_2)$ with $d_1$ generic in an open definable subset $V_1$ and $d_2\in h(d_1)$ where $h$ is a continuous correspondence on $V_1$. Moreover $dh$ exists and is continuous on a dense open subset of $V_1$ (see Proposition \ref{C1}).
Since $d_1$ is generic in $V_1$, it belongs to that dense open subset of $V_1$, that we rename $V_1$.

Further by Fact \ref{function}, there are finitely many $C^1$-functions $(\tilde h_1,\ldots,\tilde h_t)$ defined on a neighbourhood of $d_1$ such that $h(d_1)=(\tilde h_1(d_1),\ldots,\tilde h_t(d_1))$ and $d_{2})=\tilde h_1(d_1)$. W.l.o.g. we assume that this neighbourhood of $d_1$ is again $V_1$, with $\vert d_1\vert=\vert z_1\vert=:s_1$. Also since being $C^1$ is a first order property, we also assume that for each $1\leq i\leq t$, $f(x,d_1,\tilde h_i(d_1))$ is a $C^1$-correspondence in a neighbourhood of $a$. 
We consider the correspondence $G(x,z_1)$ sending $(x,z_1)$ to first $(x,z_1,h(z_1))$ and then each $(x,z_1,\tilde h_i(z_1))$ to $f(x,z_1,\tilde h_i(z_1))$. We will first show that if $G(x,d_1)$ is $C^1$ at $a$, then $G(x,z_1)$ is $C^1$ at $(a,d_1)$.
For doing that first step there is no harm in assuming that $h$ is a $1$-correspondence.
Let $a\in O$ and decompose $a$ as $(a_1,a_2)$ with $(a_1,d_1)$ $\acl$-independent and $a_2\subset \acl(a_1,d_1)$. Write $a_2:=(a_{21},\ldots, a_{2\ell})$ with $\ell= \vert a_2\vert$. Again by the cell decomposition theorem (see Fact \ref{cd}), we may assume that $a_1$ generic in an open definable subset $O_1\subset K^{n-\ell}$.
So for each $1\leq i\leq \ell$, there is a $\cL$-definable correspondence $\alpha_i$, with $a_{2i}\in \alpha_i(a_1,d_1)$; furthermore we may assume that $\alpha_i$ is $C^1$ on a neighbourhood $W\subset O_1\times V_1$ of $(a_1,d_1)$ and so $\partial$ is compatible with each $\alpha_i$. Set $\alpha(x_1,z_1):=(\alpha_1(x_1,z_1),\ldots, \alpha_\ell(x_1, z_1))$.
After a permutation of the variables, rewrite $x=(x_2,x_1)$.
We proceed as in \cite[Proposition 5.1]{PPP2}, in order to show that $G(x_2,x_1,z_1)$ is $C^1$ in a neighbourhood of $(a,d_1)$. Set $y:=(x_1,z_1)$.
Consider the map sending $(x_2,y)\mapsto (x_2-\alpha(y),y)$ and let $\tilde G(x_2,y):=G(x_2-\alpha(y), y)$. Note that $\tilde G(0, a_1,d_1)=G(a_2,a_1,d_1)$. We will show that $\tilde G$ is $C^1$ at $(0,b)$.
\par One uses the following Claim (see \cite[Claim 5.5]{PPP2}), with $b=(a_1,d_1)$. Let $X\subset K^{n+s_1}$, recall that $Int(X)$ denote the interior of $X$ and for $b\in K^{n-\ell+s_1}$, let $X^b:=\{u\in K^{\ell}: (u,b)\in X\}$.

\cl Let $X$ be $\emptyset$-definable and let $(a_2,b)\in X$ with $b$ $\acl$-independent over $\emptyset$. If $a_2\in Int(X^b)$, then $(a_2,b)\in Int(X)$.
\ecl
\par The claim was proven using \cite[Fact 5.2(1)]{PPP2}, which can be restated as follows.
\cl There are some parameters $c$ and open subset $O_{c}$ of $O_{a_2}$ containing $a_2$ defined over $c$ such that $b$ is generic over $c$. 
\ecl

\prcl 
Recall that by taking $\U$ sufficiently saturated, we may assume that the topology on $K$ is given by a valuation and using the multiplicative notation we denoted it by $\vert\cdot\vert$. Suppose that $O_{a}$ is given by $B(a,\gamma)$, where $\gamma$ is in the value group of $K$, namely $B(a,\gamma)=\{u\in K: \vert a-u\vert <\gamma\}$. Then choose $e\in K$ such that $\gamma=\vert e\vert$. So the former formula translates to: $\vert a-u\vert<\vert e\vert$, (equivalently in the one-sorted case by $div(e,a-u)$). First choose $e'$ such that $\gamma'=\vert e'\vert$ is such that $B(a, \gamma')\subset B(a,\gamma)$ and $d_1$ is still generic over $e'$. Then choose $a'\in B(a,\gamma')$ with $d_1$ generic over $e', a'$. Since $B(a',\gamma')=B(a,\gamma)$, we get the required open set $U_{c}$. \qed

\medskip
\par Then by replacing $\tilde G$ by $\tilde G(x,y)-\tilde G(0,a_1,d_1)$, we reduce to showing that $G$ is $C^1$ in a neighbourhood of $(0,a_1, d_1)$ with now $\tilde G(0,a_1,d_1)=0$. Then we follow the proof of \cite[Proposition 5.1]{PPP2}. So we replace $\tilde G$ with $\tilde G(x_2,y)-\tilde G(0,y)$ in order to get that $\tilde G(0,y)$ is zero in a neighbourhood of $(0,a_1,d_1)$ and so are all its partial derivatives with respect to $y$.
\par We have to show that $\vert \tilde G(x_2,y)-\sum \partial_{x_{2i}} \tilde G (0,a_1,d_1) (x_2-0,y-(a_1,d_1))\vert<\eta \vert (x,y-(a_1,d_1))\vert$.
\par By assumption, $\tilde G(x_2,a_1,d_1)$ is $C^1$ at $0$ and so $0\in Int(\{x_2\in K^{.}: \vert \tilde G(x_2,a_1,d_1)- -\sum \partial_{x_{2i}} \tilde G (0,a_1,d_1) x_2\vert<\eta \vert x_2\vert\}$.
\par Since $(a_1,d_1)$ is generic, $(0,a_1,d_1)\in Int\{(x_2,y): \vert \tilde G(x_2,a_1,d_1)- -\sum \partial_{x_{2i}} \tilde G (0,a_1,d_1) x_2\vert<\eta \vert x_2\vert\}$.
So there exists $\eta_1\in \vert K\vert$ such that $\vert \tilde G(x_2,a_1,d_1) -\sum \partial_{x_{2i}} \tilde G (0,a_1,d_1) x_2\vert<\eta \vert x_2\vert.$
\par Set $H(y):=\sum \partial_{x_{2i}} \tilde G (0,y)$. Then there exists $\eta_2$ such that if $\vert y-(a_1,d_1)\vert<\eta_2$, then $\vert H(y)-H(a_1,d_1)\vert \epsilon$.

\par So we get $\vert \tilde G(x_2,y)-\sum \partial_{x_{2i}} \tilde G (0,a_1,d_1)x_2\vert \leq \max\{\vert \tilde G(x_2,y)-\sum \partial_{x_{2i}} \tilde G (0,y) x_2\vert, \vert \sum \partial_{x_{2i}} \tilde G (0,y) x_2-\sum \partial_{x_{2i}} \tilde G (0,a_1,d_1) x_2\vert<\epsilon \vert x_2\vert \leq \epsilon (x_2,y-(a_1,d_1))\vert.$
\medskip

\par So $G$ is $C^1$ in a neighbourhood of $(a_2,a_1,d_1)$ and in order to deduce that $\partial$ is compatible with $G$, we no longer assume that $h$ is a $1$-correspondence but we consider the general case in order to use that $G$ is $\emptyset$-definable, we write:

\[
\partial(G(a,d))=\sum_{i=n-\ell+1}^{n} \partial_{x_{2i}} G(a_2,a_1,d_1) \partial(a_{2i})+\sum_{i=1}^{n-\ell} \partial_{x_{1i}} G(a_2, a_1,d_1) \partial(a_{1i})+\sum_{i=1}^{s_1} \partial_{z_{1i}} G(a_2,a_1,d_1) \partial(d_{1i}).
\]
Set $f^{\partial}(x):=\sum_{i=1}^{s_1} \partial_{z_{1i}} G(x_2,x_1,d_1) \partial(d_{1i})=\sum_{i=1}^{s_1} \partial_{z_{1i}} G(x_2,x_1,d_1) \partial(d_{1i})$. Then $f^{\partial}$ is continuous in a neighbourhood of $a$.

Using that  $\partial(f(a,d))=\partial(G(a,d))$, we get that 

\[
\partial(f(a,d))=\sum_{i=1}^{n} \partial_{x_{i}} f(a,d) \partial(a_{i})+f^{\partial}(a).
\]

\dfn \label{D-cell} Let $k$ be a differential subfield of $K$. Let $(X,\rho_X)\subset K^n$ be a $k$-cell, where $\rho_X(X)$ is an open $\cL(k)$-definable subset in some $K^d$ and $X$ is the graph of a continuous $\cL(k)$-definable $m$-correspondence $f$ (see Definition \ref{partial}). Then $(X,\rho_X)$ is a $\partial$-compatible $k$-cell, if for each $1\leq i\leq \ell$, the partial derivatives $\partial_{i} f$ with respect to $x_{i}$ exist and are continuous (on $\rho_X(X)$).
 We will say that $(X,\rho_X)$ is a $\partial$-compatible $C^1$-$k$-cell if $f$ is $C^1$ (on $\rho_X(X)$).
\medskip
\par Then a $D$-cell $(X,s)$ (over $k$) is a $\partial$-compatible $k$-cell $(X,\rho_X)$, equipped with a section $s:X\to X \times K^{n}$ defined as follows. By Fact \ref{function}, each $a\in \rho_X(X)$ has an open neighbourhood $V$ such that there are $m$ continuous definable functions $g_i: V\to K^{n-d}$, $1\leq i\leq m$, such that $\graph(f\restriction V)=\bigsqcup_{i=1}^m \graph(g_i)$.
 The current assumptions on $f$, implies that $\partial$ is compatible with each $g_i$, $1\leq i\leq m$, and for $1\leq j\leq d$, $\partial_j f\restriction V=(\partial_j g_1,\ldots, \partial_j g_m)$ and $f^{\partial}\restriction V=(g_1^{\partial},\ldots, g_m^{\partial})$.
Then $s$ is an $\cL(k)$-definable function sending $a\in X$ to $(a,b)$ with for some $1\leq i\leq m$, $a=(a_1, g_i(a_1))$, $b=(b_1,b_2)$, $b_1\in K^d$ and $b_2=J_{g_i}(a_1) b_{1}+ g_i^{\partial}(a_1)$.
\edfn
\par Note that 
if $b_1=\partial(a_1)$, then $\partial(g_i(a_1))=J_{g_i}(a_1) \partial(a_{1})+ g_i^{\partial}(a_1))$ since $\partial$ is compatible with $g_i$, $1\leq i\leq m$.
\medskip
\lem \label{alg++} Let $T$ be either dp-minimal and not strongly minimal or an open theory of topological fields. Let $(\M,\partial)\models T(\partial)$. Let $k$ be a differential subfield of $K$ and let $\bar a\in K$. Let $u\in \acl(\bar a,\ldots, \partial^{m-1}(\bar a), k)$.
 Then $\partial(u)\in \dcl_{\cL}(\bar a,\ldots, \partial^{m}(\bar a),u,k)$
\elem
\pr We extract from $\bar a,\ldots, \partial^{m-1}(\bar a)$ an $\acl$-independent (over $k$) subtuple that we denote $\bar d$. Now let $\bar e$ an $\acl$-independent subtuple in $k$ such that $u\in \acl(\bar d,\bar e)$.
Let $O\subset K^{\ell}$ be a definable open set containing $(\bar d,\bar e)$.  By Fact \ref{cd}, there is a continuous $\cL$-definable $m$-correspondence $f$ on an open large subset $\tilde O$ of $O$ containing $(\bar d,\bar e)$ and such that $u\in f(\bar d,\bar e)$. Moreover by Proposition \ref{partial} and since $(\bar d,\bar e)$ is $\acl$-independent, we may assume that the partial derivatives of $f$ exist on $\tilde O$ and are continuous. Since $\partial$ is compatible with $f$, 
we have $\partial( f(\bar d,\bar e))=\sum_{i=1}^n \partial_i f(\bar d,\bar e)\partial(d_i)+\sum_{j=1}^\ell \partial_j f(\bar d,\bar e) \partial(e_j)$, with $\bar d=(d_1,\ldots,d_n)\subset \bar a,\ldots, \partial^{m-1}(\bar a)$, $\bar e=(e_1,\ldots,e_\ell)\in k$.
\par By Fact \ref{function}, there is a continuous $\cL(u)$-definable function $g$ on a definable neighbourhood $V$ of $(\bar d, \bar e)$
such that $u=g(\bar d,\bar e)$.
 The function $g$ is defined as follows: one takes an open neighbourhood $W$ of $u$ such that $W\cap \graph(f)$ is the graph of a function that we denote by $g$. 
Furthermore we may again assume that the partial derivatives of $g$ exist and are continuous on $V$ by Proposition \ref{partial} and Corollary \ref{partialopen}.
\par So $\partial(u)=\partial(g(\bar d,\bar e))=\sum_{i=1}^n \partial_{i} g(\bar d,\bar e) \partial(d_i)+\sum_{j=1}^\ell \partial_j g(\bar d,\bar e) \partial(e_j)\in \dcl_{\cL}(\bar a,\ldots, \partial^{m}(\bar a),u,k)$.
\qed
\rem The proof of the above lemma is modelled on \cite[Lemma 2.3]{FK} in the o-minimal case. When $T$ is an open theory of topological fields, an equivalent form of the above lemma was proven in \cite[Lemma 3.1.9]{K-P}.
\erem
\medskip
\medskip
\par We  gather below other compatibility properties between  the derivation $\partial$ and the language $\cL$ (besides Lemma \ref{alg++}) that will play a role later in the description of $\cL_\partial$-definable groups.
\par Below we only assume that $T$ is a complete geometric theory of fields of characteristic $0$. Let $(\U,\partial) \models T(\partial)$ and let us single out the following properties (C1) up to (C3) that can hold in models of $T(\partial)$ and which are well-known in (pure) differential fields.
\begin{enumerate}[label={(C\arabic*)}] 
\item In models of $T(\partial)$, $\partial\restriction{\acl(\emptyset)}=0$.
\item Let  $\M\models T$ be an $\cL$-substructure of $\U$. Let $A$ be a subset of the field sort of $\M$ and suppose that $\partial\restriction A: A\to A$. Then for any finite subset  $B$ of the field sort of $\M$ consisting of $\acl$-independent elements over $A$, we can define a map $\partial^*$ on $B$ by sending $B$ to any non-empty subset of $\acl(A,B)\cap M$ and expand $\M$ to a model of $T(\partial)$ where the derivation agrees with $\partial$ on A and with $\partial^*$ on $B$. Furthermore the derivation $\partial^*$ on $B$ extends in a unique way on $\acl(B)$.
\item  Let $A\subset K$ with $\partial(A)=A$. Let $a$ be a finite tuple in $K$ and $u\in K$. Then if  $u\in \acl_{\cL}(a, A)$, then
 $\partial(u)\in \dcl_{\cL}(\nabla(a), u, A)$.
\end{enumerate}
\par As a consequence of the above conditions, letting $a$ be a $m$-tuple of $\acl$-independent elements of $K$, $b\in K^m$ and $u\in \acl(a)$, we have that in $tp_\cL(a,u,b)$, an 
$\cL$-formula $\varphi(x,y,z,w)$ such that $\exists^{=1} w\; \varphi(a,b,u,w)$ and further if $\partial$ is a derivation sending $a$ to $b$ then $w=\partial(z)$. Indeed,  any derivation sending $a$ to $b$ extends uniquely on $u\in \acl(a)$. So $\partial(u)$ is $\cL$-definable over $a, b,u$. Let us show that $\varphi$ only depends on $tp(a,u,b)$.
\lem\label{unicity} Assume that $T(\partial)$ satisfies condition (C2) and let  $(\M,\partial)\models T(\partial)$. Let $(k,\partial)$ be a substructure of $(\M,\partial)$ and suppose that $\cM$ is $\vert k\vert^+$-saturated and satisfies condition (C2). Let $a$ be a uple of $\acl$-independent elements of the field sort of $\M$ and $d\in \acl(a/k)$. 
Let $\varphi$ be an $\cL(k)$-formula  such that $(\M,\partial)\models \varphi(\nabla(a), d,\partial(d))\wedge (\exists^{=1} z \varphi(\nabla(a), d,z))$. Let $b=\partial(a)$, then 
\[
T(\partial)\cup tp_{\cL}(a, d, b/k)\cup Diag(k,\partial)\models (\exists^{=1} z \varphi(x,y, u, z))\wedge \varphi(x,\partial(x),u,\partial(u)),
 \]
with $c$ (respectively $c', u$) a tuple of new constants with $\vert c\vert=\vert a\vert$ (respectively $\vert c'\vert=\vert b\vert$, $\vert u\vert=1$).
\elem
\pr Let $\cN\models T$ containing $(k,\partial)$. Let $\tilde a$ be a tuple of $\acl$-independent elements in the field sort of $\cN$ with $\vert \tilde a\vert=\vert a \vert$, let $\tilde b$ be a tuple of elements of the field sort of $\cN$ with $\vert \tilde b\vert=\vert b \vert$ and let $\tilde d$ in the field sort of $\cN$ and in addition suppose that $tp_{\cL}(a,d, b)=tp_{\cL}(\tilde a, \tilde d, \tilde b)$.
 Let $\cM'$ be an $\cL$-elementary extension of $\M$ such that there is an elementary embedding $f$ of $\cN$ in $\cM'$ sending $\tilde a$ to $a$, $\tilde d$ to $d$, $\tilde b$ to $b$ and being the identity on $k$. 
 On $\cN$, define a derivation $\partial^*$ as $f^{-1} \partial f$. Since in $\M$ we have $\varphi(a,\partial(a),d,\partial d)$, in $\cN$ we have $\varphi(\tilde a,\partial^*(\tilde a), u, \partial^*(u))$. But $\partial^*(\tilde a)=\tilde b$, and by (C2), any other derivation on $\cN$ extending $\partial$ on $k$ and sending $\tilde a$ to $\tilde b$ will send $\tilde d$ to $\partial^*(\tilde d)$.\qed
 \medskip

\rem \label{c2} In case $T$ is an $\cL_{rings}$-theory of geometric fields of characteristic $0$, condition (C2) is a restatement of the well-known extension theorem for derivations \cite[Chapter 10, section 7, Theorem 7]{Lang} (see also \cite[Theorem 1.1]{PP}). Indeed, let $\M\models T$ and let $A$ be a subfield of $M$ equipped with a derivation $\partial$. Then for any finite subset  $B$ of the field sort of $\M$ consisting of algebraically independent elements over $A$, we can extend $\partial$ on the subfield $A(B)$ by sending $B$ to any non-empty subset of $A(B)$ and then to the relative algebraic closure of $A(B)$ inside $M$. Then we iterate the construction of $\partial$ inside $M$.
\par Since the language $\cL$ in case $T$ is an open theory of topological fields is a relational expansion of the language of fields (together possibly with additional constants) and $\acl$ is the field algebraic closure \cite[Proposition 1.3.4]{K-P}, it is easy to see that (C2) holds as well in that setting.
\erem 
\lem \label{condition2} Let $T$ be dp-minimal , not strongly minimal. Let $(\U,\partial)\models T(\partial)$ (with field sort $K$) and assume that $\M\models T$ be an $\cL$-substructure of $\U$. Let  $A\subset K(M)\subset K$ with $\partial\restriction A: A\to  A$ and $\acl(A)=A$. Then for any finite subset  $B\subset K(M)\setminus \acl(A)$ of $\acl$-independent elements over $A$, we can extend $\partial$ on $B$, sending $B$ to any non-empty subset of $\acl(A,B)\cap M$ and expand $M$ to a model of $T(\partial)$ where the derivation extends the one already defined on $A\cup B$. In particular $T(\partial)$ satisfies the compatibility condition (C2).
\elem
\pr We proceed as in \cite[Lemma 2.13]{FK}. 

\par Expand $B$ to an $\acl$-independent subset $\tilde B$ of elements of $K(M)$ over $A$ such that $K(M)=\acl(\tilde B, A)$. Let $c\in K(M)\setminus A$, so there is a finite subset $\{b_{1},\ldots, b_{n}\}$ of $\tilde B$, with $c\in \acl(b_{1},\ldots, b_{n}, \bar a)$, for some parameters $\bar a\in A$. W.l.o.g. we may assume that $\bar a$ consists of $\acl$-independent elements (over $\emptyset$).
This is witnessed by an $\cL$-formula $\varphi(x,\bar y,\bar z)$ with $\vert x\vert=1$, $\vert \bar y\vert=n$, $\vert \bar z\vert=\vert \bar a\vert=\ell$ and $\M\models \exists^{\leq m} \varphi(x,\bar b,\bar a)$, $\bar b:=(b_{1},\ldots, b_{n})$. By \cite[Proposition 4.1]{simon-walsberg2016}, there is an $\cL$-definable continuous $m'$-correspondence $f$ defined on an open definable set $O$, $m'\leq m$, such that $c\in f(\bar b, \bar a)$. Since $(\bar b,\bar a)$ are $\acl$-independent, by Proposition \ref{partial}, we may assume that the partial derivatives of $f$ exist and are continuous on an open neighbourhood $O'$ of $(\bar b,\bar a)$.
Further, by assumption on $(\U,\partial)$, $\partial$ is compatible with $f$. So $\partial( f(\bar b,\bar a))=\sum_i \partial_{i=1}^n f(\bar b,\bar a)\partial(b_i)+\sum_{j=1}^\ell \partial_j f(\bar b,\bar a) \partial(a_j)$.
\medskip
\par Now let $\bar d$ be any $n$-tuple in $\acl(A,B)\cap K(M)$ and define $\partial^*(f(\bar b, \bar a)):=\sum_{i=1}^n \partial_i f(\bar b,\bar a) d_i+\sum_{j=1}^\ell \partial_j f(\bar b,\bar a) \partial(a_j)$ and $f^{\partial^*}(\bar b,\bar a):=\sum_{j=1}^\ell \partial_j f(\bar b,\bar a) \partial(a_j)$.
\medskip
\par First let us show this is well-defined (w.r. to the chosen subset $\tilde B$ of $\acl$-independent elements). Let $g$ be another $\cL$-definable continuous $m''$-correspondence defined on an open definable set $O$, $m''\leq m$, such that $c\in g(\bar b, \bar a')$, with $(\bar b,\bar a'')$ $\acl$-independent. Suppose that $1\leq \vert g(\bar b,\bar a')\cap f(\bar b,\bar a)\vert=\tilde m$. Again, by Proposition \ref{partial}, we may assume that the partial derivatives of $g$ exist and are continuous on an open neighbourhood $O''$ of $(\bar b,\bar a')$.
By Fact\ref{function} there is an open neighbourhood $W_c$ of $c$, an open neighbourhood $V$ of $\bar b$, $V_1$ an open neighbourhood of $\bar a$ (respectively $V_2$ of $\bar a'$) such that $\graph(f) \cap (V\times V_1)\times W_c$ (respectively $\graph(g) \cap (V\times V_2)\times W_c$ ) is the graph of a function that we denote by $f_1$ (respectively $g_1$).

Let us show that 
\[
\sum_i \partial_{i=1}^n f_1(\bar b,\bar a) d_i+\sum_{j=1}^\ell \partial_j f_1(\bar b,\bar a) \partial(a_j)=\sum_i \partial_{i=1}^n g_1(\bar b,\bar a') d_i+\sum_{j=1}^\ell \partial_j g_1(\bar b,\bar a') \partial(a'_j).
\]

\par Let $\pi_1$ be the coordinate projection on the first coordinate. Consider $X:=\{x\in \pi_1(V): \vert g(x,b_{2},\ldots, b_{n},\bar a'))\cap f(x,b_{2},\ldots, b_{n},\bar a)\vert=\tilde m\}$, since $b_{1}$ is not an isolated point in $X$, by  \cite[Proposition 4.1]{simon-walsberg2016}, $X$ contains an open subset containing $b_{1}$. We iterate the reasoning for $b_{2}, \ldots, b_{n}$ and we find an open neighbourhood $\tilde V$ of $\bar b$ where $f(\bar y,\bar a)$ and $g(\bar y,\bar a')$ coincide on $\tilde m$ values containing $c$. In particular $f_1(\bar y,\bar a)$ and $g_1(\bar y,\bar a')$ coincide on $\tilde V$. So $J_{f_1}(\bar y,\bar a)$ and $J_{g_1}(\bar y,\bar a')$ coincide on that open neighbourhood. Since in $(\U,\partial)$, $\partial$ is compatible with $f$ and $g$, we have  that $\partial( f(\bar y,\bar a))=\sum_i \partial_{i=1}^n f(\bar y,\bar a)\partial(y_i)+\sum_{j=1}^\ell \partial_j f(\bar y,\bar a) \partial(a_j)$ and an analogous statement holds for $g$. In particular it also holds for $f_1$ and $g_1$. We can express with an $\cL$-formula in $\bar a, \bar a'$ and $c$ that in $\U$ the following holds: $\forall \bar y\in \tilde V$
\[
f_1(\bar y,\bar a)=g_1(\bar y,\bar a')\rightarrow \bigwedge_{1\leq i\leq n} \partial_{i} f_1(\bar y,\bar a)=\partial_{i} g_1(\bar y,\bar a')) \wedge (\sum_{j=1}^\ell \partial_j f_1(\bar y,\bar a) \partial(a_j)=\sum_{j=1}^\ell \partial_j g_1(\bar y,\bar a) \partial(a'_j)).
\]
Since $\M\models T$ and $T$ is model-complete, the above formula also holds in $\M$. So $\partial^*(c)$ is well-defined.

\par Now let us show that the derivation $\partial^*$ is compatible with every $\cL$-definable continuous $m$-correspondence $g$ on some definable open set $O\subset K^n$, with the property that for each $1\leq i\leq n$, the partial derivatives $\partial_{i} g$ with respect to $x_{i}$ exist (on $O$) and are continuous. 
Let $\bar u\in O$. Since $\bar u\in \acl(\tilde B, A)$, for each $u_i$, $1\leq i\leq n$, there is  a finite tuple $\bar b\subset \tilde B$ and $\bar a\in A$ such that $u_i\in f_i(\bar b,\bar a)$ and $f_i$ is  an $\cL$-definable $m_i$-correspondence and as before we may assume that $\bar b, \bar a$ are $\acl$-independent, that there is an open set $V\subset K^m$ containing $(\bar b,\bar a)$ and for each $1\leq i\leq m$, the partial derivatives $\partial_{j} f_i$ with respect to $y_{j}$ exist and are continuous on $V$, with $\vert \bar y\vert=\vert \bar b\vert$. 

Further, we may assume that $f_i(V)\subset O$. Let $f:=(f_1,\ldots,f_m)$.
\par We consider $h:=g\circ f$, it is  $\cL$-definable continuous correspondence on $V$ and its partial derivatives $\partial_{j} h$ with respect to $y_{j}$ exist and are continuous on $V$. 

Then by the way $\partial^*$ was defined, $\partial^*(h(\bar b))=J_{h}(\bar b) \partial^*(\bar b)+h^{\partial}(\bar b)=\partial^*(g(\bar u))$.

 By the above, we have $\partial^*(\bar u)=(\partial^*(u_1),\ldots,\partial^*(u_n))\in (\partial^*(f_1(\bar b),\ldots,\partial^*(f_n(\bar b)))=\partial^*(f(\bar b))$, where $\partial^*(f(\bar b))=J_{f}(\bar b) \partial^*(\bar b)+f^{\partial}(\bar b).$  

\par It remains to show that $\partial^*(g(\bar u))=J_{g}(\bar u) \partial^*(\bar u)$. 

This follows from the rule of composition of partial derivatives: $J_{h}(\bar x)=J_{g}(f(\bar x)) J_{f}(\bar x)$ and the statement (in $\cL$) which holds in $\U$ that: $h^{\partial}(\bar x)=J_{g}(f(\bar x)) f^{\partial}(\bar x)$.

We apply this to $\bar b$ and we get $\partial^*(g(\bar u))=J_{g}(\bar u) (J_{f}(\bar b) \partial^*(\bar b)+f^{\partial}(\bar b))=J_{g}(\bar u)\partial^*(\bar u)$.
\qed

\subsection{Cell decomposition}
\dfn\label{delta_compa_funct}  Let $T$ be either dp-minimal, not strongly minimal or $T$ an open theory of topological fields. Let $(\U,\partial)\models T(\partial)$.
Let $f$ be an $\cL(A)$-definable correspondence $f:K^n\to K$ with $A$ a differential subfield of $K$. We will say that $f$ is $\partial$-compatible if there is a finite partition of the domain of $f$ into $\cL(A)$-definable $\partial$-compatible cells $C_{i}=(\rho_{C_{i}}(C_{i}), h_{i})$ (see Definition \ref{D-cell}), such that the partial derivatives of $f\restriction C_{i}\circ h_{i}$
exist and are continuous on $\rho_{C_{i}}(C_{i})$ and $\partial$ is compatible with $f\restriction C_{i}\circ h_{i}$. 
\par We will say that $f$ is a $\partial$-compatible $C^1$-correspondence if there is a finite partition of the domain of $f$ into $\cL(A)$-definable $\partial$-compatible $C^1$-cells $C_{i}=(\rho_{C_{i}}(C_{i}), h_{i})$, such that the correspondence sending $x_1\mapsto f(x_1,h_i(x_1)$ that we denote by $f\restriction C_{i}\circ h_{i}$, is $C^1$ and $\partial$ is compatible with $f\restriction C_{i}\circ h_{i}$.
\par Let $\pi_{(j)}$ be the projection from $K^{m}$ to $K$, sending $(x_{1},\ldots,x_{m})$ to $x_{j}$, $1\leq j\leq m$ and $f:K^n\to K^m$. Let $f_{j}:K^n\to K$ be equal to $\pi_{(j)}\circ f$, $1\leq j\leq m$. Then we say that $f$ is $\partial$-compatible (respectively a $\partial$-compatible $C^1$-correspondence) if each $f_{j}:K^n\to K$ is $\partial$-compatible (respectively $\partial$-compatible $C^1$-correspondence). 
\edfn
\prop \label{C1-decomp-dp} Let $\U$ be a dp-minimal not strongly minimal field with a geometric theory.
Consider a differential expansion $(\U,\partial)\models T(\partial)$ and let $A$ be a differential subfield of $\U$. 
Let $X$ be a $\cL(A)$-definable subset of $K^n$.
Then $X$  can be partitioned into a finite union of $\partial$-compatible $C^1$-$A$-cells. Moreover, if $f$ is a definable $m$-correspondence, then $f$ is a  $\partial$-compatible $C^1$-correspondence.
\eprop
\pr We first apply \cite[Theorem 1.3]{WJ}. Then we follow the proof of \cite[Theorem 1.5.3, Corollary 1.5.5]{K-P}, which by induction first on $n$ and then on $\dim(X)$. We replace the property that a correspondence on open set $U$ is continuous on an open dense definable subset of $U$ by the assumption that $f$  
is $C^1$ and $\partial$-compatible 
on an open definable dense subset of $U$. The other ingredient is that the property of being a $\partial$-compatible $C^1$-cell is closed under restriction to an open subset of $\rho_{C}(C)$, where $C$ is a $\partial$-compatible $C^1$-cell \cite[Lemma 1.5.2]{K-P}. 
\par Now let $f$ be a definable $m$-correspondence in $K^n$. By the first part $\dom(f)$ can be partitioned into a finite union of $\partial$-compatible $C^1$-cells $(\rho_{C_i}(C_i),h_i)$, where $\rho_{C_i}(C_i)$ is an open subset of some $K^{d_i}$ and $h_i$ a $C^1$-correspondence from $K^{d_i}$ to $K^{n-d_i}$, $0\leq d_i\leq n$, $i\in I$ with $I$ finite. Let $(a_1,a_2)\in \dom(f)\cap C_i$ with $a_1\in \rho_{C_i}(C_i)$ and consider
the correspondence $f\circ h_i$ sending $a_1$ to $f(a_1,h_i(a_1))$. By Proposition \ref{C1}, this correspondence is $C^1$ on an open dense subset of $\rho_{C_i}(C_i)$. Then we re-apply the first part of the proposition to the subset of $u\in \rho_{C_i}(C_i)$ on which $f\circ h_i$ is not $C^1$ on a neighbourhood of $u$. This definable set can be expressed as a finite union of $\partial$-compatible $C^1$-cells of smaller dimension than $\rho_{C_i}(C_i)$. Iterating, we arrive to the desired decomposition of $\dom(f)$.
\qed
\medskip
\par From now on we assume that $T$ is geometric and model-complete and that $T(\partial)$ has a model-companion $T_{\partial}$.  Let $(\U,\partial)$ be a sufficiently saturated model of $T_\partial$.
\medskip

\dfn\label{AA}
\begin{enumerate}[label={(A\arabic*)}]
\item any $\cL_\partial$-formula $\varphi(x)$ is equivalent in $T_\partial$ to a formula of the form $\psi(\bar \partial^m(x))$, for some $\cL$-formula $\psi(x_0, x_1,\ldots, x_m)$, with $x_0=x$, $\vert x_i\vert=\vert x\vert$, $1\leq i\leq m$ and $m\geq 0$. 
\item Let $X\subset K$ and let $a\in K$ with $a\in \dcl_{\partial}(\langle X\rangle)$, then there is $m\geq 0$ such that $a\in \dcl(\nabla_{m}(\langle X\rangle))$.
\end{enumerate}
\edfn
\rem\label{openAA}
Note that when $T$ is an open theory of large topological fields of characteristic $0$ and $T$ is model-complete, then $T_\partial$ exists, is axiomatizable and model-complete and properties (A1) and (A2) hold \cite[Theorem 2.4.2, Corollary 2.4.9, Lemma 3.1.9]{K-P}.
Moreover, if $T$ has a model which is elementary equivalent to a complete non-trivial valued field of rank $1$, then all models of $T$ are large fields.
\erem

\lem Suppose in addition that $T$ is dp-minimal, not strongly minimal. Then $T_\partial$ has properties (A1) till (A2).
\elem
\pr 
Using Proposition \ref{C1-decomp-dp} and that property (C2) holds in models of $T(\partial)$ (see Lemma \ref{condition2}), we will show that (A1) holds in models of $T_\partial$.
\par First (A1) holds for quantifier-free $\cL_\partial$-formulas since $\partial$ is a compatible derivation. Indeed, suppose we have an atomic $\cL_\partial$-formula with a subterm of the form $\partial(t)$ where $t$ is an $\cL$-term.

Suppose $\vert x\vert=n$ and consider the $\cL_\partial$-function associated with $t(x)$ from $K^n\to K$. Then, by Proposition \ref{C1-decomp-dp}, $t(x)$ is a $\partial$-compatible $C^1$-function, namely there is a partition of $\dom(t)$ into finitely many cells $(\rho_{C_i}(C_i), h_i)$ with $\rho_{C_i}(C_i)$ an open set in some $K^{d_i}$ and $\tilde t:=t\restriction C_{i}\circ h_i$ a $\partial$-compatible $C^1$-correspondence into $K^{n-d_i+1}$. So for $x\in C_i$ with $x=(x_1,h_i(x_1))$, $x_1\in \rho_{C_i}(C_i)$, we have $t(x)=t(x_1,h_i(x_1))=\tilde t_i(x_1)$ and $\partial(\tilde t_i(x_1))=\sum_ {j=1}^{d_i}\partial_j \tilde t_i(x_1) \partial(x_{1j})$, where $1\leq \vert x_1\vert=d_i\leq n$.

So we replace $\partial(t(x))=y$ with $\vert y\vert=n$ by $\bigvee_i (x\in C_i \wedge x=(x_1, h_i(x_1)\wedge y=\sum_ {j=1}^{d_i}\partial_j \tilde t_i(x_1) \partial(x_{1j}))$ which can be expressed by an $\cL$-formula in $(\partial(x_{1j}))_{j=1,\ldots, d_i}$.

\par Since $T_\partial$ is model-complete, any $\cL_\partial$-formula $\chi(x)$ is equivalent to an existential formula $\exists y\;\varphi(y, x)$ with $\varphi$ a quantifier-free $\cL_\partial$-formula. So $\varphi(y, x)$ is equivalent to $\psi(\nabla_m(y),\nabla_\ell(x))$, with $\psi(\bar u, \bar z)$ an $\cL$-formula and $\vert \bar u\vert=(m+1)\vert y\vert$, $\vert \bar z\vert=(\ell+1)\vert x\vert$. 
\par By induction on the number of the existential quantifiers occurring in the formula $\exists y\;\varphi(y, x)$, we show that such formula is equivalent  to $\tilde \psi(\nabla_\infty(x))$, for some $\cL$-formula $\psi$.
\par It suffices to consider the case where $\vert y\vert=1$ and suppose that the quantifier-free$\cL_\partial$-formula $\varphi(y, x)$ is equivalent to the formula $\psi(\nabla_m(y),\nabla_\ell(x))$, where $\psi(y_0,\ldots,y_m,\bar z)$ is an $\cL$-formula, with $\vert \bar z\vert=\vert \nabla_\ell(x)\vert$. We will show that there is an $\cL$-formula $\chi$ such that the formula $\exists y\;\psi(y,\partial(y),\ldots,\partial^m(y),\bar z)$ is equivalent to $\exists y\;\chi(y,\nabla_\infty(\bar z)$.
\par Let $(K,\partial)\models T_\partial$. Consider the tuple $(y_0,\ldots, y_m)$ and consider the projections $\pi_\ell$ of $\psi(K)$ onto its $\vert \bar z\vert+\ell$ coordinates, $0\leq \ell\leq m$, (with for $\ell\geq 1$, $y_{\ell-1}$ as its ${\vert \bar z\vert+\ell}^{th}$ coordinates). Let $\ell\geq 1$, by Proposition \ref{C1-decomp-dp}, $\pi_\ell(\psi(K))$ is a finite disjoint union of $\partial$-compatible $C^1$-cells $X_i$, $i\in I$. 
\par Recall that a cell in $K^{\vert \bar z\vert+\ell+1}$ is either a finite set or an open set or it is of the form $(X,\rho_X)$ where $X$ is the graph along $\rho_X$ of a $\partial$-compatible $C^1$-correspondence and there is a cell $(C,\rho_C)$ such that $C=\pi(X)$ where $\pi$ is the coordinate projection from $K^{\vert \bar z\vert+\ell+1}$ to $K^{\vert \bar z\vert+\ell}$ and 

\par (i) $\rho_X=\rho_C\circ \pi$ or
\par (ii) $\rho_X(\bar z, y_0,\ldots, y_{\ell})=(\rho_C(\bar z, y_0,\ldots, y_{\ell-1}), y_\ell)$ and for every $c\in C$, the fiber $X_c$ is an open subset of $K$ and $f_X(\rho_X(\bar z, y_0,\ldots,y_\ell))=f_C(\rho_C(\bar z,y_0,\ldots,y_{\ell-1}))$. 

\par For each cell $X_i$, $i\in I$, we let $\ell$ be minimum such that case (i) occurs and denote by $C_i$ the cell previously denoted by $C$. Note that if there is no such $\ell$ it implies that in the cell $X_i$ that the tuple $(y_0,\ldots,y_m)$ varies in an open set. Even though $\ell$ depends on $i$, we simply denote it by $\ell$.
\par So, we have that $y_\ell\in f_{X_i}(\rho_{C_i}(\bar z,y_0,\ldots, y_{\ell-1}))$ and since $\partial$ is compatible with $f_{X_i}$, $\partial(y_\ell)$ is $\cL$-definable from
 $y_\ell,\ldots, y_0, \nabla_\infty(\bar z)$. 
 \par We replace the formula $\psi$ by a finite disjunction of formulas $\chi_i(y_0,\ldots,y_\ell,\bar z)$ and $\chi_i$ is obtained for $(y_0,\ldots,y_\ell,\bar z)\in X_i$ by replacing $y_\ell$ by $f_{X_i}(y_0,\ldots, y_{\ell-1},\bar z)$, by setting $y_{\ell+i}$ by $\partial^i(f_{X_i})$ with $1\leq i\leq m-\ell$ in case there is a minimum $\ell$ and requiring that $\psi(y_0,\ldots,y_m,\bar z)$ holds.
\par Let $\bar a\in K$ such that $(K,\partial)\models \exists y\;\varphi(y,\bar a)$. Let $A$ be the differential subfield of $K$ generated by $\bar a$. Let $(b_0,\ldots, b_{\ell},\bar a)$ be a tuple of elements in some $X_i\cap \pi_{\ell}(\psi(K))$. By choice of $\ell$, we may choose $(b_0,\ldots,b_{\ell-1})$ $\acl$-independent over $A$. By (C2), there is a derivation $\partial^*$ sending $(b_0,\ldots,b_{\ell-1})$ to $(b_1,\ldots, b_{\ell})$ and extending $\partial$ on $A$.  Furthermore, we can expand $\partial^*$ on $K$ such that $(K,\partial^*)$ is a model of $T(\partial)$. Since $T_\partial$ is the model-companion of $T(\partial)$, we can embed $(K,\partial^*)$ in a model of $T_\partial$, which we may assume to be an elementary extension of $(K,\partial)$. So if $(K,\partial)\models \exists y_0\ldots\exists y_m \bigvee_{i\in I} \chi_i(y_0,\ldots,y_\ell,\nabla_\infty(\bar a)$, then  $(K,\partial)\models \exists y\;\varphi(y,\bar a)$.

\medskip
\par It remains to show (A2). Let $u\in \dcl_{\partial}(\langle X\rangle)$. Let $\varphi$ be an $\cL_\partial$-formula with $\exists^{=1}y\; \varphi(y,\bar z)$ and such that for some $\bar d\in \langle X\rangle$, $\varphi(u,\bar d)$ holds. Suppose that $\varphi$ has order $m$. Let $\bigvee_{i\in I} \chi_i$ be the disjunction of $\cL$-formulas constructed above.  Then we have that for some $i$, $\chi_i(u_0,\ldots,u_\ell,\nabla_\infty(\bar d))$ holds, with $0\leq \ell\leq m$. By construction if $\ell>0$, then we have infinitely many differential solutions of $\chi_i$ and if $\ell=0$ it means that $u_0=f_{X_i}(\bar d)$. Moreover if $f_{X_i}$ is not a $1$-correspondence, then we would have more than one solution to $\varphi$.\qed  

\medskip
\subsection{Open topological fields} 
Let us refine the former cell decomposition theorem for open theories of topological fields (see Fact \ref{cd} and \cite[Theorem 1.5.3]{K-P}).

\prop \label{C1-decomp} Let $T$ be an open theory of topological fields of characteristic $0$ and let $\U$ be a model of $T$ and further assume that the field sort $K$ is elementary equivalent to a complete non-trivial valued field of rank $1$. Let $(X,\rho_{X})$ be a $k$-cell in the field sort $K$, $X\subset K^{n}$, $k$ a differential subfield of $K$. Then $X$  can be further partitioned into a finite union of $\partial$-compatible $C^1$-$k$-cells.
\eprop
\pr  
We proceed by induction on $n$ and use the description of cells given in Definition \ref{def:cells}. The case $n=1$ is clear. So let us assume the proposition holds for $X\subset K^n$ and let us show it holds for $n+1$. Let us assume that $0<\dim(X)<n+1$. In particular $\dim(\pi(X))\leq n$, where $\pi$ is the coordinate projection on the first $n$-variables.
\par By induction hypothesis, $C=\pi(X)$ can be partitioned into a finite union of $k$-cells $C_i$ which are graphs of $\cL(k)$-definable $m_i$-correspondences $h_i$ which are $C^1$ on $\rho(C_i)$.
\par In case (c) we have that $\rho_X(x,y)=(\rho_C(x),y)$ and $f_X(\rho_X(x,y))=f_C(\rho_C(x))$ with $\dim(\rho_C(C))<n$. 
Applying the induction hypothesis, we have that for $x$ such that $\pi(x)\in C_i$, $(\rho_C(x),f_C(x))=(\rho_{C_i}(x),f_{C_i}(x))$ where $f_{C_i}$ is a $\cL(k)$-definable $m_i$-correspondence which is $C^1$ on $\rho_{C_i}(C_i)$. So for $x$ such that $\pi(x)\in C_i$, we have that $f_X(\rho_X(x,y))=f_{C_i}(\rho_{C_i}(x))$. Since $\partial$ is compatible with $f_{C_i}$, so is $f_X$ since the value of $f_X$ does not depend on the last variable.

\par In case (b), $\rho_X(X)=\rho_C\circ\pi(X)$ and $f_X$ is a $\cL(k)$-definable $m$-correspondence, which is continuous on $\rho_C\circ\pi(X)$.
By \cite[Corollary 1.3.3]{K-P}, $X$ is a finite union of $\cL(k)$-definable subsets of the following form:
\[
 Z_{\cA_\ell}^{S_\ell}(x,y) \wedge \theta_\ell(x,y)
\]
where $\vert x\vert= n$, $\vert y\vert=1$, and $\theta_\ell$ is an $\cL(k)$-formula that defines an open subset of $K^{n+1}$, $S_\ell\in k[x,y]$, $\cA_\ell\subseteq k[x,y]\setminus\{0\}$, $\cA_\ell^y=\{P_\ell\}$ and $\partial_{ y} P_\ell$ divides $S_\ell$.
\par Further we may assume that for $\ell\neq \ell'$, $Z_{\cA_\ell}^{S_{\ell}}\cap  Z_{\cA_\ell'}^{S_{\ell'}}=\emptyset$ (otherwise we refine the given covering by performing the Euclidean algorithm). Let $(a,b)\in X$ with $\vert a\vert=n$ and suppose that for some $\ell$, $P_\ell(a,b)=0$, $\partial_y P_\ell(a,b)\neq 0$ and $\theta_\ell(a,b)$ holds. In particular $(a,b)$ belongs to the locally Zariski closed set $Z_{P_\ell}^{\partial_y P_\ell}$.
We apply Corollary \ref{imp_p-def}. So there are neighbourhoods $O_a\subset K^n$ of $a$, $O_b\subset K$ of $b$ and a definable $C^1$ function $g: O_a\to O_b$ such that $g(a)=b$ and $P_\ell(a,g(a))=0$. Further $g(a)$ is the unique solution of $P(a,y)=0$ in $O_{b}$.
Moreover by Lemma \ref{implicit} (in the Appendix), $\partial$ is compatible with $g$. From now on we will denote such function by $g_{\ell(a,b)}$.

\par Since $\pi(X)=\bigsqcup_i C_i$, there is a unique $i$ such that $a\in C_i$, so up to a permutation of coordinates $a=(a_1,a_2)$ with $\rho_{C_i}(a)=a_1$ and $a_2\in h_i(a_1)$ with $h_i$ a $C^1$-$m_i$-correspondence on $\rho_{C_i}(C_i)$, with $\partial$ compatible with $h_i$, by induction. By Fact \ref{function}, on a neighbourhood $V$ of $a_1$, $h_i=(h_{i1},\ldots,h_{im_i})$, where $h_{ij}$ are $C_1$-definable functions and since $\partial$ is compatible with $h_i$ on a neighbourhood of $a_1$ which is sent by $h_{ij}$ into $O_a$. 
On a neighbourhood of $a_1$, consider the maps sending $a_1\mapsto (a_1,h_{ij}(a_1), g_{\ell(a,b)}(a_1,h_{ij}(a_1)))$, $1\leq j\leq m_i$; it is $C^1$ on $\rho_{C_i}(C_i)$. 
Note that $a=(a_1,h_{ij}(a_1))$ and $b=g_{\ell(a,b)}(a_1,h_{ij}(a_1))$ and since $(a,b)\in X$, given $a_1$ there are exactly $m$ distinct tuples of the form $(h_{ij}(a_1), g_{\ell(a,b)}(a))$. 
\par The required partition of $X$ is obtained by defining on $\rho_{C_i}(C_i)$ the correspondence sending $a_1$ to $(h_{ij}(a_1), g_{\ell(a,b)}(a_1, h_{ij}(a_1)))$ with $a=(a_1,h_{ij}(a_1))$ and $(a,b)\in X$. As noted above there are only $m$ possibilities and this correspondence is definable and $C^1$ on $\rho_{C_i}(C_i).$
\qed

\prop \label{prop:comp1} Let $T$ be an open theory of topological fields and suppose that $T$ has a model which is a complete nondiscrete valued field of rank $1$.
Let $\K$ be a model of $T(\partial)$. Let $f$ be an $\cL(A)$-definable correspondence $f:K^n\rightrightarrows K$ with $A$ a differential subfield of $K$. 
Then there is a finite partition of the domain of $f$ into $\cL(A)$-definable cells $C_{i}=(\rho_{C_{i}}(C_{i}), h_{i})$, where $h_{i}$ is a $C^{1}$-definable, $\partial$-compatible correspondence and finitely many open disjoint subsets $W_{\ell}$ of $K^n$, finitely many $C^1$-definable correspondences $g_{\ell}$ on $W_{\ell}$ with $\partial$ compatible with $g_\ell$ such that for $x\in W_{\ell}\cap C_{i}$, $f(x)=g_{\ell}(\rho_{C_{i}}(x), h_{i}(\rho_{C_{i}}(x)))$. 
\par If $f$ is defined on an open subset $U\subset K^n$.
Then on a large open subset of $U$, $f$ is a $C^1$-$\cL(A)$-definable correspondence and $\partial$ is compatible with $f$.
\eprop
\pr The proof of this proposition is similar to the proof of Case (2) of the previous proposition; sometimes, one says that $f$ is {\it prepared} on the cells $C_{i}$. Let $\pi$ be the projection from $K^{n+1}$ to $K^n$ sending $(x_{1},\ldots,x_{n+1})$ to $(x_{1},\ldots,x_{n})$.
\par By Proposition \ref{C1-decomp}, the domain of $f$ is partitioned into a finite union of cells which are graphs of $C^1$-correspondences with whom $\partial$ is compatible.
 
\par By \cite[Corollary 1.3.3]{K-P}, the graph of $f$, as a definable subset of $K^n\times K$ can be expressed as
 a finite union of definable subsets of the following form:
\[
 Z_{\cA_\ell}^{S_\ell}(x,y) \wedge \theta_\ell(x,y)
\]
where $\vert x\vert= n$, $\vert y\vert=1$, and $\theta_\ell$ is an $\cL$-formula that defines an open subset of $K^{n+1}$, $S_\ell\in K[x,y]$, $\cA_\ell\subseteq K[x,y]\setminus\{0\}$, $\cA_\ell^y=\{P_\ell\}$ and $\partial_{ y} P_\ell$ divides $S_\ell$.
\par Further we may assume that for $\ell\neq \ell'$, $Z_{\cA_\ell}^{S_{\ell}}\cap  Z_{\cA_\ell'}^{S_{\ell'}}=\emptyset$ (otherwise we perform Euclidean algorithm to decompose an intersection of the form $Z_{\cA_j}^{S_j}(x,y)\cap Z_{\cA_i}^{S_i}(x,y)$ as a finite union of sets where $\vert \cA^y\vert= 1.$)). 
Also since this finite union is the graph of a function, we may assume that for $\ell\neq \ell'$, $\pi(Z_{\cA_\ell}^{S_\ell}\cap O_{\ell})\cap \pi(Z_{\cA_{\ell'}}^{S_{j'}}\cap O_{\ell'})=\emptyset$.

Let $(a,b)\in graph(f)$ with $\vert a\vert=n$ and suppose that $P_\ell(a,b)=0$, $\partial_y P_\ell(a,b)\neq 0$ and $\theta_\ell(a,b)$. 
We apply Corollary \ref{imp_p-def}. So there neighbourhoods $O_a\subset K^n$ of $a$, $O_b\subset K$ of $b$ and definable $C^1$ function $g:O_{a}\to O_b$ such that $g(a)=b$ and $P_\ell(a,g(a))=0$.
On the locally Zariski closed set $Z_{P_\ell}^{\partial_y P_\ell}$, we define the function $g_\ell$ such that for $(a',b')\in Z_{P_\ell}^{\partial y P_\ell}$, we have $P_\ell(a',g_{\ell}(a'))=0$ and there is a neigbourhood $O_{a'}\times O_{b'}$ of $(a',b')$ on which $g_\ell(x)$ is the unique solution of $P_{\ell}(x,y)=0$. Moreover $g_\ell$ is $C^1$ and by Lemma \ref{implicit} (in the Appendix), $\partial$ is compatible with $g_\ell$.
Let $W_{\ell}:=\{x\in K^n: \vert\{y: (x,y)\in Z_{P_\ell}^{\partial_y P_\ell}\}\vert=1\}$.
\par Since $\dom(f)=\bigsqcup_i C_i$, there is a unique $i$ such that $a\in C_i$, so up to a permutation of coordinates $a=(a_1,a_2)$ with $\rho_{C_i}(a)=a_1$ and $a_2\in h_i(a_1)$ with $h_i$ a $C^1$-$m_i$-correspondence on $\rho_{C_i}(C_i)$, which is $\partial$-compatible. By Fact \ref{function}, on a neighbourhood $V$ of $a_1$, $h_i=(h_{i1},\ldots,h_{im_i})$, where $h_{ij}$ are $C_1$-definable functions and since $h_i$ is $\partial$ compatible, $\partial$ is compatible with each $h_{ij}$, $1\leq j\leq m_i$. Let $a_2=h_{ij}(a_1)$ for some $1\leq j\leq m_i$. Consider the map sending $a_1\mapsto (a_1,h_{ij}(a_1), g_\ell(a_1,h_{ij}(a_1)))$; it is $C^1$ and $\partial$ is compatible with it on a neighbourhood of $a_1$ included in $\rho_{C_i}(C_i)$ (as a composition of $C^1$-functions with whom $\partial$ is compatible). 
\par Then we define on $(\rho_{C_i}(C_i))$ the correspondence sending $x_1$ to $(g_{\ell}(x_{1},h_{ij}(x_1)))_{1\leq j\leq m_{i}}$. We have that for some $\ell$, $f(x_{1},h_{ij}(x_{1}))=g_{\ell}(x_{1},h_{ij}(x_{1}))$, $1\leq j\leq m_{i}$. 
\par Moreover if $\dom(f)=U$ an open subset of $K^n$, then $I_{0}:=\{i: \dim(C_{i})=n\}$, the subset $\bigcup_{i\in I_{0}} C_{i}$ is a large open subset of $U$, on which $\partial$ is compatible with $f$. Indeed, if $a\in C_{i}$, $i\in I_{0}$, (in particular $C_{i}$ is open), then for some $\ell$, $f(x)$ coincide with $g_{\ell}(x)$ on a neighbourhood of $a$ and so we may define $f^{\partial}$ as $g_\ell^{\partial}$ and it has the required properties.  
We have  $\partial(f(a))=\partial(g_{\ell}(a))=\sum_{i=1}^{n} \frac{\partial}{\partial x_i} g_{\ell}(a) \partial(a_i)+g_{\ell}^{\partial}(a)$. Since both $f$ and $g_{\ell}$ are $C^1$ on $C_{i}\cap W_{\ell}$ and coincide there, we have  $\{\frac{\partial}{\partial x_i} g_{\ell}(a)\}= \frac{\partial}{\partial x_i} f(a)$. So $\partial$-compatible is compatible with $f$ on $C_{i}$. Therefore, $f$ is $C^1$ and $\partial$ is compatible with $f$ on $\bigcup_{i\in I_{0}} C_{i}$. \qed

\medskip
\subsection{Topological groups}
In this subsection, we will place ourselves into two settings: either $T$ is a dp-minimal not strongly minimal geometric theory, or $T$ is an open  theory of topological fields with the property that $T$ has a model which is a complete nondiscrete valued field of rank $1$. 
\par We will use the $C^1$-cell decomposition (see Propositions \ref{C1-decomp-dp}, \ref{C1-decomp}), in order to show that define a topology on a definable group in such a way the group operations are $C^1$-maps, following \cite{P88}. 
\prop\label{prop:topo} 
Let $(G, \times, F_{-1})$ be group definable in $\U$.
Then there is a definable topology $\tau_{g}$ on $G$
such that $(G,\times, F_{-1})$ is an $\cL$-definable topological group where the group operations $\times$ and $F_{-1}$
 are $C^1$-$\partial$-compatible maps. Moreover the topology $\tau_{g}$ on $G$ is related to the topology on $K$ as follows. There is a large definable subset $V$ of $G$ such that $V$ is the finite disjoint union of $\partial$-compatible cells $C_{i}$ which are relatively open in $G$ and open in $\tau_{g}$.
\eprop
\pr Assume that $G\subset K^n$. It suffices to check that we can apply \cite[Proposition 2.5]{P88}, choosing large subsets of $G$ (respectively of $G\times G$) on which the group operations are $C^1$-$\partial$-compatible maps and satisfy the pre-group properties (see \cite[section 4]{PPP2}). 
\par We decompose $G$ as a finite disjoint union of $\partial$-compatible cells $C_{i}$ and we only keep the cells of maximal dimension: $C_{1},\ldots,C_{r}$. Let $V:=\bigsqcup_{i=1}^r C_i$. 
Note that if $x$ is generic and $x\in C_j$, $1\leq j \leq r$, there is $1\leq i\leq r$ such that $F_{-1}(x)\in C_{i}$. Let $C_j^i$ be the finite disjoint union of $\partial$-compatible relatively open cells included in $C_j$ such that for some $x\in C_j$, $x$ generic, $F_{-1}(x)\in C_i.$ We have that $\bigsqcup_{i=1}^r C_j^i$ is a large relatively open subset of $C_j$ and $V_{0}=\bigsqcup_{i, j=1}^r C_j^i$ is large and open in $G$. 
By Propositions \ref{C1-decomp-dp}, \ref{prop:comp1}, we may assume that on $C_j^i$,  the inverse map $F_{-1}$ is a $C^1$-$\partial$-compatible map.
Let $V_{1}:=\{a\in V_{0}\colon F_{-1}(a)\in V_{0}\}$; it is a large subset of $G$, $V:=V_{1}\cap F_{-1}(V_{1})$ is invariant under $F_{-1}$ and $V$ is open in $G$.
\par Similarly we define $Y_{i,j}^k$, $1\leq i, j, k\leq r$, to be the finite disjoint union of $\partial$-compatible relatively open cells included in $C_i\times C_j$ such that for some $(a,b)\in C_{i}\times C_{j}$, with $a, b$ are independent generic, we have $a\times  b\in C_k$ (recall that $a\times b$ is again generic). 
By Propositions  \ref{C1-decomp-dp}, \ref{prop:comp1}, we may assume that on $Y_{i,j}^k$, the map $\times$ is a $C^1$-$\partial$-compatible map. 
Then $\bigsqcup_{k=1}^r Y_{i,j}^k$ is large, relatively open in $C_i\times C_j$. Let $Y_{0}:=\bigsqcup_{1\leq i ,j \leq r}(\bigsqcup_{k=1}^r Y_{i,j}^k)$, it is a large subset of $G\times G$, on which $\times$ is $C^1$-$\partial$-compatible. Let $Y:=\{(a,b)\in Y_{0}\colon a\times b\in V\}\cap (V\times V)$; it is open and large in $Y_{0}$.
\par Note that for any $a\in V$ if $b$ is generic in $V$ over $a$, then $(b,a)\in Y$ and $(F_{-1}(b),b\times a)\in Y.$
\par Then one replaces the lemma in \cite[Proposition 2.5]{P88},  by the following:
\begin{enumerate}[label=\alph*)]
\item for any $a, b\in G$, the set $\{x\in V: ((a\times x)\times b)\in V\}$  is open in $V$ and the map $x\mapsto ((a\times x)\times b)$ and its inverse are $C^1$-$\partial$-compatible maps.
\item for any $a, b\in G$, the set $\{(x,y)\in V\times V: (((a\times x)\times b))\times y)\in V\}$ is open in $V\times V$ and the corresponding map is $C^1$-$\partial$-compatible.
\end{enumerate}

\par Therefore, we may  use \cite[Proposition 2.5]{P88} and define a topology $\tau_{g}$ on $Z$ as follows. 
\par For $g\in G$ and $U\subset Z$, denote by $g\times U:=\{u\in G: \exists a\in U\;g\times a=u\}$. Then take as open sets, the subsets $U\subset G$ such that for all $g\in G$, $g\times U\cap V$ is open. The key property is that if $U\subset V$ and $a\in G$, then $a\times U$ is open in $\tau_{g}$ iff $U$ is open in $V$. 
Endowed with this topology, $G$ is a topological group where the group operations are $\partial$-compatible $C^1$-maps.

\qed
\nota\label{cor:topo} 
A $\cL$-definable group $G$ which can be endowed with a definable topology for which the group operations become $\partial$-compatible $C^1$-maps will be called  a $C^1$-group.
\enota

\section{$\cL_\partial$-Definable sets in models of $T_\partial$}\label{sec:defsets}
In this section, we don't necessarily place ourselves in the class of (non discrete) topological fields. We suppose that 
\begin{enumerate}
\item $T$ is geometric, model-complete,
\item $T(\partial)$ has a model-companion $T_\partial$ and $T(\partial)$ satisfies the compatibility conditions (C1) up to (C3),
\item $T_\partial$ satisfies the compatibility conditions (A1) till (A2). 
\end{enumerate}
Let $(\U,\partial)$ denote a sufficiently saturated model of $T_{\partial}$, $K$ the field sort and $k$ now denotes a small differential subfield of $K$. Let us first define the analog of a prolongation $\tau(V)$ of an irreducible algebraic set $V\subset K^n$ (see \cite[section 1]{M}). Recall that if $V$ is smooth, then $\tau(V)$ is again irreducible and that $\tau(V)$ is the Zariski closure of $\nabla(V)$.
In the following we will replace $V$ by an arbitrary $\cL(k)$-definable set $X$ and so to $X$ we will associate a finite union of $\cL(K)$-definable sets in which $\nabla(X)$ is dense. In this setting we will use assumption (C3) on $T(\partial)$.
\medskip
\par First let us introduce some notation. Let $a=(a_1,\ldots, a_n)\in K^n$ and suppose that $m=\dim(a/k)$. Let $a_{j_1},\ldots, a_{j_m}$ be $\acl$-independent (over $k$). Set $[n-m]:=\{1,\ldots, n\}\setminus\{j_1,\ldots,j_m\}$; note that when $n=m$, $[0]=\emptyset$. For $i\in [n-m]$, we have $a_{i}\in \acl(a_{j_1},\ldots, a_{j_m}, k)$.
\notalem \label{notsection} Let $X$ be an infinite $\cL(k)$-definable set in $(\U,\partial)$, included in $K^n$ for some $n\in \N$ and assume $0<\dim(X)<n$.
Let $a\in X$ with $\dim(a)=m$. Let $a[m]:=(x_{j_1},\ldots, x_{j_m})$ be a subtuple of $a$ of $\acl$-independent elements and let $b\in K^m$. By Lemma \ref{unicity}, we have for each $a_{i}\in \acl(a_{j_1},\ldots, a_{j_m}, k)$, $i\in [n-m]$, 
an $\cL(k)$-formula $\varphi_i$ such that
\[
T(\partial)\cup tp_{\cL}(a[m], a_i, b/k)\cup Diag(k,\partial)\models \exists^{=1} z_i \varphi_i(x,u_i, y, z_i))\wedge \varphi_i(x,u_i,\partial(x),\partial(u)),
 \]
with $\vert x\vert=\vert y\vert=m$, $\vert u_i\vert=\vert z_i\vert=1$.
So the $\cL$-type $p:=tp(a, b/k)$ contains an $\cL(k)$-formula $\chi_i$, $i\in [n-m]$ such that $\chi_i(a[m],a_i,b)$ holds and such that
\begin{equation}\label{eqdcl}
T(\partial)\cup Diag(k,\partial) \cup\chi_i(x,u,y)\models \exists^{=1} z \varphi_i(x,u,y, z))\wedge \varphi_i(x,u,\partial(x),\partial(u)),
\end{equation} 
Let $\chi(x,u,y):=\bigwedge_{i\in [n-m]} \chi_i(x,u_i,y)$ with $\vert u\vert=n-m$ and $\varphi(x,u,y,z):=\bigwedge_{i\in [n-m]} \varphi_i(x,u_i,y,z_i)$ with $\vert z\vert=n-m$.
By construction the projection onto the first $n$-coordinates of $\chi(K)$ has dimension $m$.  
For $i\in [n-m]$, let $f_{p,i}$ be the $\cL(k)$-definable function (partially defined) sending $a_{i}$ to the element $b_{i}$ such that $\varphi_{i}(a[m], a_i, b,b_i)$ holds. Denote by $f_p$ the function sending $(a, b)$ to $(f_{p,i}(a[m], b, a_{i}))_{i\in [n-m]}$.  Note that the function $f_p$ defined above now only depends on $\chi$; so from now on, we will change the subscript and denote it by $f_{\chi}$.
By compactness, varying the tuple $(a,b)\in X\times K^m$, we have finitely many $\cL$-formulas $\chi, \varphi$ with $\chi\in tp(a,b)$ 
such that the above implication (\ref{eqdcl}) holds in $T(\partial)\cup Diag(k,\partial)$; denote this set by $\cF(\nabla(X)).$

\enotalem

\dfn\label{defsection} Let $X$ be an $\cL(k)$-definable set in $(\U,\partial)$, included in $K^n$ for some $n\in \N$. 
\par First assume that $0<\dim(X)<n$, for each $(\chi,\varphi)\in \cF(\nabla(X))$,
define 
the prolongation $\tau(X)\subset X\times K^n$ as follows.
\[
\tau(X):=\bigcup_{(\chi,\varphi)\in \cF(\nabla(X))}\{(a,b)\in X\times K^{n}\colon (\chi\wedge \varphi)(a,b)\}.
\]

\par Then we define a family of maps $s$ from $X\to \tau(X)$ as follows. Let $a\in X$ with $a[m]$ a maximal subtuple of $\acl$-independent elements, then there is $b\in K^m$ such that $\chi(a,b)$ holds. Then let  $\varphi$ be an $\cL$-formula such that (\ref{eqdcl}) holds in $T(\partial)\cup Diag(k,\partial)$; so  $\varphi(a, b, f_{\chi}(a,b))$ holds.

Define a map $s:X\to \tau(X)$ by  sending $a$ to $(a, b, f_{\chi}(a,b))\in \tau(X)$. Note that if $b=\partial(a[m])$, then $\partial(a)=(b, f_{\chi}(a,b))$.
\par Second assume that $\dim(X)=n$, then define $\tau(X):=X\times K^n$ and for $a\in X$, define $s(a)=(a,b)$ by choosing for each $a\in X$ some $b\in K^n$.
\smallskip
\par Let $\pi$ be the projection from $\tau(X)$ to the first $n$-coordinates. Then $\pi\circ s$ is the identity on $X$; we call $s$ a section associated with $(\chi,\varphi)\in \cF(\nabla(X))$.
We denote this set of sections by $\cS(\nabla(X))$.
\edfn
\ex \label{ex:section}
\par (A) Let $\ACF_0$ be the theory of algebraically closed fields of characteristic $0$, $\cL_r$ the language of rings, $(K,\partial)$ an algebraically closed differential field of characteristic $0$ and $X$ an algebraic subset of some $K^n$. (Note that in this case $\acl$ is simply the field algebraic closure.)
\par Let $a$ be a
 tuple in $X$ and let $I(a)$ be the ideal in $K[X]$ of polynomials which vanish at $a$. Suppose that $I(a)$ is generated by $p_1,\ldots, p_{\ell}$, with $p_i\in k[x]$, $1\leq i\leq \ell$ with $k$ a differential subfield of $K$.  Then 
 \[\tau(X)=\bigcup_{a\in X, I(a)=\langle p_1,\ldots,p_{\ell}\rangle} \{(a,b)\in X\times K^n\colon \bigwedge_{i=1}^{\ell} (p_i(a)=0\;\wedge\;0=p_i^{\partial}(a)+\sum_{j=1}^n \partial_j p_i(a) b_j=0\},\] with $b=(b_1,\ldots,b_n)$.
 \par When $X$ is irreducible, then the above definition coincides with the classical one. Indeed, let $a$ a generic point in $X$ and  $I(a)=\langle p_1,\ldots,p_{\ell}\rangle$, then
 $$\tau(X)=\{(a,b)\in X\times K^n\colon \bigwedge_{i=1}^\ell (p_i(a)=0\;\wedge\;0=p_i^{\partial}(a)+\sum_{j=1}^n \partial_j p_i(a) b_j=0)\}.$$
  \par Let $X_{a}:=\{u\in X: I(u)=I(a)\}=V(I(a))$, then $a$ is a generic point of $X(a)$ and in particular a simple point. So the rank of the Jacobian matrix $(\partial_j p_i)_{1\leq i\leq \ell, 1\leq j\leq n}$ is equal to $n-\dim(V).$
 We can extend $\partial$ on $k[x]/I(a)$ and then on $Frac(k[x]/I(a))$ by sending a maximal subset (of cardinality $\dim(V)$) of algebraically independent elements of $(a_1,\ldots,a_n)$ arbitrarily and then the image of the $n-\dim(V)$ other ones is uniquely determined by the above equations.

\medskip
\par (B)  Now assume that $T$ is dp-minimal, not strongly minimal or that $T$ is an open theory of topological fields with the property that $T$ has a model which is a complete nondiscrete valued field of rank $1$.
Let $K\models T$ and let $X$ be a cell. Let $a\in X$ and the form $(a_1,f(a_1))$, where $a_1\in U$, an open subset of some $K^m$ and $f$ be a $C^1$-correspondence from $U\to  K^{n-m}$. Then $\partial(f(a_1))=\partial(a_2)=\sum_{j=1}^m \partial_j(f(a_1)) \partial(a_{1j})+f^{\partial}(a_1).$ 
\par In case $a_1$ is $\acl$-independent in $U$, we can extend $\partial$ on $a_1$ freely and then it is uniquely determined on $a_2.$ So in this case, 
\begin{align*}
\tau(X):=\{(a,b)\in X\times K^n\colon& a=(a_1,a_2)\wedge b=(b_1,b_2)\wedge a_1\in U\wedge f(a_1)=a_2\wedge b_1\in K^m\wedge\\ &b_2=\sum_{j=1}^m \partial_j(f(a_1)) b_{1j}+f^{\partial}(a_1)\}.
\end{align*}
\eex

\lem \label{section}
Let $X$ be an $\cL(k)$-definable set in $(\U,\partial)$ with $X\subset K^n$ for some $n\in \N$,  let $p(x)$ be an $\cL(k)$-generic type of $X$. 
Let $s$ be a section associated with $X$.
Then there is some  realization $a$ of $p$ such that $s(a) = \nabla(a)$. 
\elem
\begin{proof} 
Let $p(x)$ be a complete $\cL(k)$-type with $\dim(p) = \dim(X)=m\leq n$. Let $a= (a_{1},..,a_{n})$ realize $p(x)$ in $\U$. Recall that $s(a)$ was defined as follows.
Letting $a[m]$ be a maximal subtuple of $\acl$-independent elements of $a$, we choose $b\in K^m$ such that $\chi(a,b)$ holds for some $(\chi,\varphi)\in \cF(\nabla(X))$ with $\varphi$ be an $\cL$-formula such that (\ref{eqdcl}) holds.

Then $s(a):=(a, b, f_{\chi}(a,b))$ with $\varphi(a, b, f_{\chi}(a,b))$. Further the following property holds: if $b=\partial(a[m])$, then $b, f_{\chi}(a,b)=\partial(a)$.
 Let $\M$ be a small elementary submodel of $\U$ containing $a, b$ and $k$, but $\vert k\vert^+$-saturated. 
\par Extract from $b$ an $\acl$-independent subtuple $b[m']$ over $(k, a[m])$ with $m'\leq m$. Let $q:=tp_{\cL}(a,b)$.
 \par By the compatibility condition (C2), there is a derivation $\partial^*$ on the field sort of $\M$, extending $\partial$ on $k$, sending $a[m]$ to $b$, for $i\in [n-m]$, $a_{i}$ to $b_{i}=f_{q,i}(a[m],b, a_{i})$, and $b[m']$ to $c[m']$, where $c[m']$ has been chosen arbitrarily in $\acl(a[m],b[m'],k)=\acl(a, b,k)\subset M$.

 \par So $(\M,\partial^{*})$ is a model of $T(\partial)$ with $\partial^*(a[m])=b$, hence embeds in a model of $T_\partial$ which we may assume to be $(\U,\partial)$. Note that this embedding is the identity on $k$. Let $(\tilde a,\tilde b)$ be the image of $(a,b)$ by this (elementary) embedding. Then $\tilde a$ realizes $p$ and $(\tilde a[m], \tilde b)$ realizes $q$ (and so $\chi\wedge\varphi (\tilde a,\tilde b))$ holds. Therefore  $s(\tilde a)=\nabla(\tilde a)$. 
\end{proof}
\bigskip
\par We extend the dimension function $\dim$ induced by $\acl$ on the field sort  to certain $\cL_{\partial}$-definable subsets of the field sort in models of $T(\partial)$. When we wish to stress that we restrict ourselves to the language $\cL$, we use $\cL$ as a subscript.
\dfn \label{def-dimpartial} Let $a$ be a tuple in $K$ and $k$ a subset of $K$. Let $k$ be a small differential subfield of $K$. We define $\dim_{\partial}(a/k)$ to be $\max\{\dim(\nabla_{N}(a)/k)\colon N\in \N\}$, when this maximum is infinite, we set $\dim_\partial(a/k):=\infty$.
\par\noindent Let $X\subset K^n$ be an $\cL_{\partial}(k)$-definable set in $(\U,\partial)$, then $\dim_{\partial}(X):=\max\{\dim_{\partial}(a/k)\colon a\in X\}$, when this maximum is infinite, we set $\dim_\partial(X)=\infty$. We say that $X$ is finite-dimensional (over $k$) if $\dim_{\partial}(X)\in \N$.
\par Suppose $X$ is finite-dimensional, 
then $a\in X$ is {\it generic}  (w.r.to $\dim_\partial$) if $\dim_{\partial}(a/k)=\dim_{\partial}(X).$ Note that when $X$ is finite-dimensional and $a\in X$, then there is some $N\in \N$ such that $\dim_{\partial}(a/k)=\dim(\nabla_{N}(a)/k)$.
\par Suppose that $X$ is finite-dimensional and let $a, b\in X$, 
 then the analog of equation (\ref{dim}) is as follows:
\begin{align}\label{dim2}
\dim_{\partial} (a,b/k)=\dim_{\partial}(a/k \nabla_{\infty}(b))+\dim_{\partial}(b/k),
\end{align}
since for some $N\in \N$ and any $c\in X$, $\nabla_{\infty}(c)\subset \acl(\nabla_{N}(c))$,
\[
 \dim_{\partial} (a,b)/k=\dim(\nabla_{N}(a),\nabla_{N}(b)/k)=\dim(\nabla_{N}(a)/k\nabla_{N}(b))+\dim(\nabla_{N}(b)/k).
 \]
\par Let $b, d$ be two tuples of $K$. Then $b$ and $d$ are $\partial$-independent over $k$ if $\dim_{\partial}(b,d/k)=
\dim_{\partial}(b/k)+\dim_{\partial}(d/k)$.
\edfn
\medskip

\par In the following lemma we give a description of finite-dimensional $\cL_{\partial}$-definable sets (finite-dimensional as defined in Definition \ref{def-dimpartial}) which will be used to give a reasonably canonical description of finite-dimensional groups $\cL_{\partial}$ definable in a saturated model of $T_{\partial}$, adapted from well-known constructions in $\DCF_{0}$. 
\medskip
\lem\label{set} 
Let $X$ be a finite-dimensional $\cL_{\partial}(k)$-definable set in $\U$ with $X\subseteq K^{n}$, for some $n$. Then for suitable $N$, there is an $\cL(k)$-definable set $Y\subset K^{n(N+1)}$,
there are finitely many $(\chi,\varphi)\in \cF(\nabla(Y))$, finitely many $\cL(k)$-definable sets $Y_{\chi,\varphi}$ with $Y=\bigcup Y_{\chi,\varphi}$, such that 
there is an $\cL(k)$-definable section $s:Y\to \tau(Y)$ such that
\newline
(i) if $\nabla^N(b)\in Y$ with $b\in K^n$, then $b\in X$,
\newline
(ii) $\nabla_{N}(X) = \bigcup_{(\chi,\varphi)\in \cF(\nabla(Y))} \{a\in Y_{\chi,\xi}: \nabla(a) = s(a)\}$,
\newline
(iii) if $p$ is an $\cL(k)$-generic type of some $Y_{\chi,\xi}$, then it is realized by some $\nabla_{N}(a)$ with $a\in X$.

\elem
\pr  
(i) By (A1), there is an $\cL(k)$-formula $\psi$ such that $\psi(\bar \partial^{\bar m}(x))$ defines $X$.  For ease of notation we assume that $\bar m$ is a constant tuple that we denote by $m$. 
Note that $\nabla_{m}(X)$ can be defined as follows: $\nabla_{m}(X):=\{(\bar x_{0},\bar x_{1},\ldots,\bar x_{m})\colon \psi(\bar x_{0},\bar x_{1},\ldots,\bar x_{m})\wedge \partial (\bar x_{0},\bar x_{1},\ldots,\bar x_{m-1})=(\bar x_{1},\bar x_{2},\ldots,\bar x_{m})\}$, with $\bar x_{0}=x$ and $\bar x_j=(x_{j1},\ldots,x_{jn})$, $j\geq 0$.

\par\noindent Since $\dim_\partial(X)\in \N$, there is $N\in \N$ such that for each $a\in X$, $\partial^{N}({a})\subset \acl_{\cL}(k, a, \cdots, \partial^{(N-1)}(a))$. So by compactness there are finitely many $\cL$-formulas of the form $\bigwedge_{i=1}^n\;\xi_{i\ell}(\bar c_0,\cdots, \bar c_{N-1},x)$ with $\ell\in L$, $L$ finite, $\vert \bar c_i\vert=n$, $\vert x\vert=1$, $1\leq i\leq n$, such that for each $a\in X$ we have 
\[
\bigvee_{\ell\in L} \bigwedge_{i=1}^n\;\xi_{i\ell}(a, \cdots, \partial^{(N-1)}(a),\partial^N(a_i))\wedge  \bigwedge_{i=1}^n\;\exists^{=d_\ell(i)} x\;\xi_{i\ell}(a, \cdots, \partial^{(N-1)}(a),x),
\]
 $d_{i\ell}\in \N$. 
\par Set $\bar x:=({\bar x}_{0},....,{\bar x}_{N})$ and define $\tilde \psi(\bar x):=$
\[
\psi(\bar x))\wedge \bigvee_{\ell\in L} \bigwedge_{i=1}^n\;\xi_{i\ell}(\bar x_0, \cdots, \bar x_{N-1},x_{Ni})\wedge  \bigwedge_{i=1}^n\;\exists^{=d_\ell(i)} x\;\xi_{i\ell}(\bar x_0, \cdots, \bar x_{N-1},x)\}.
\]
Let $Y:=\tilde \psi(\U)$. By construction of the formula $\tilde \psi$, if $b\in K^n$ and $\nabla^N(b)\in Y$, then $b\in X$.
\par (ii) Let $\bar a\in Y$. We extract from $\bar a$ a maximal $\acl$-independent (over $k$) subtuple (note that it is a subtuple of $(\bar a_0,\ldots, \bar a_{N-1})$. We start with $\bar a_0$, denote this subtuple by $\bar a_0[m_0]$, for some $m_0\leq n$, then we repeat the process for each $\bar a_i$  $1\leq i\leq N-1$, namely let $\bar a_i[m_i]$ be a maximal $\acl$-independent subtuple of $\bar a_{i}$ (over $k, \bar a_0[m_0],\ldots, \bar a_{i-1}[m_{i-1}]$), $1\leq i\leq N-1$. Note that $m_{0}\geq \ldots\geq m_{N-1}$. (Indeed reasoning with $m_1$: we have $m_1\leq m_0$ since $a_{j1}=\partial a_{j0}\in \acl(k,\bar a_0[m_0])$, $j\in [n-m_0]$, writing $\bar a_i=(a_{1i},\ldots,a_{ni})$). Let $\bar a[\tilde N]:=(\bar a_{0}[m_{0}],\ldots,\bar a_{N-1}[m_{N-1}])$ be a maximal $\acl$-independent (over $k$) subtuple of $(\bar a_0,\ldots, \bar a_{N-1})$ (or equivalently of $(\bar a_{0},\bar a_{1},\ldots, \bar a_{N}))$, with $\tilde N=(m_{0},\ldots, m_{N-1})$. Define $g(\bar a[\tilde N])=(\bar a_{1}[m_{0}],\ldots,\bar a_{N}[m_{N-1}])$; it is a subtuple of $\bar a^1:=(\bar a_{1},\ldots,\bar a_{N})$.

\par Let $p_{i}$ be the $\cL(k)$-type of $(\bar a_0,\ldots, \bar a_i, \bar a_{i+1}[m_i])$ and let 
$f_{p_i}$ be the $\cL(k)$-definable function sending $(\bar a_0,\bar a_1,\ldots,\bar a_i, \bar a_{i+1}[m_i])$ to $\partial(\bar a_{i}[n-m_{i}])$, $0\leq i\leq N-1$ (see Notation \ref{notsection}). 

Then, let $p_N=p$ be the type $tp_{\cL}(\bar a)$ and $f_{p_N}$ be the $\cL(k)$-definable function sending $\bar a$ to $\partial(\bar a_N)$. 
\par Set $f_p(\bar a):=( f_{p_0}(\bar a_0,\bar a_1[m_0]),\ldots, f_{p_{N-1}}(\bar a_0,\ldots,\bar a_{n-1},\bar a_N[m_{N-1}]), f_{p_N}(\bar a))$. By abuse of notation denote $(g(\bar a[\tilde N]), f_p(\bar a))$ the tuple $$(a_1[m_0], f_{p_0}(\bar a_0,\bar a_1[m_0]),\ldots, a_N(m_{N-1}), f_{p_{N-1}}(\bar a_0,\ldots,\bar a_{n-1},\bar a_N[m_{N-1}]), f_{p_N}(\bar a)).$$
\par Moreover, if $g(\bar a[\tilde N])=\partial(\bar a[\tilde N])$, then $\partial(\bar a)=(g(\bar a[\tilde N]), f_p(\bar a))$ and since $\psi(\bar a)$ holds, we have that $\partial(\bar a)=(g(\bar a[\tilde N]), f_p(\bar a))$ implies that $\bar a_0\in X$.

\par Now let us define the following section $s$ on $Y$. Let $\bar a\in Y$
then we define $s(\bar a)=(\bar a, g(\bar a[\tilde N]), f_p(\bar a))$; it is $\cL$-definable and only depends on $p$.
\par By  Notation-Lemma \ref{notsection}, there is an $\cL$-formula $\chi\in p$, an $\cL$-formula $\varphi$ 
such that, with $\vert \bar x\vert=\vert \bar y\vert=n(N+1)$,
\begin{equation*}
T(\partial)\cup Diag(k,\partial)\models (\chi(\bar x,\bar y)\implies \varphi(\bar x, g(\bar x[\tilde N]), f_{p}(\bar x))\wedge \varphi(\bar x,\partial(\bar x))),
\end{equation*}

\par For $(\chi,\varphi)\in \cF(\nabla(Y))$, define $Y_{\chi,\varphi}:=\{\bar x\in Y: \chi(\bar x, g(\bar x[\tilde N]))$ and on each of these sets, we denoted (see Notation-Lemma \ref{notsection}) $f_p$ by $f_{\chi}$. Moreover, the projection onto the first $nN$ coordinates of $Y_{\chi,\varphi}$ has dimension $\tilde N$ (and $\bar x_N\in \acl(\bar x_0,\ldots, \bar x_{N-1})$. So if we take a generic point $\bar a$ of $Y_{\chi,\varphi}$, then $\bar a[\tilde N]$ is a tuple of $\acl$-independent elements and so by Lemma \ref{section}, there is a realization $\bar d$ of $tp_{\cL}(\bar a)$ with  $\nabla(\bar d)=s(\bar d)$. So $\bar d$ is a differential tuple belonging to $\nabla_N(X)$.
This concludes the lemma. \qed

\lem\label{function+} Let $X$ be a finite-dimensional $\cL_{\partial}(k)$-set in $(\U,\partial)$ included in some cartesian product of the field sort.
Let $f$ be an $\cL_{\partial}(k)$-definable function from $X\to K$.
Then for suitable $N$, there is an $\cL(k)$-definable set $Y$ in $\U$ with $Y\subset K^{n (N+1)}$ and $\nabla_{N}(X)\subset Y$ and  
there is an $\cL(k)$-definable function $F$ from $Y$ to $K$, 
such that for $x\in X$, 
\[f(x)=F(\nabla_{N}(x)).\]
\par Furthermore for $m\geq 1$, there are $\cL$-definable prolongations $F(m)$ of $F$ with domain a large subset of $\dom(F)$ and such that $\partial^m(f(x))=F(m)(\nabla_{N}(x))$.
\elem
\pr Apply Lemma \ref{set} to $\dom(f)$. Let $N$ be such that for each $\bar a\in X$, $\partial^N(\bar a)\subset \acl_{\cL}(\nabla_{N-1}(\bar a), k)$. By (C3), $\partial^{N+m}(\bar a)\subset \dcl_{\cL}(\nabla_{N}(\bar a),k)$, $m\geq 1$. Then we use assumption (A2), namely since $f(\bar a)\in \dcl_{\partial}(\bar a, k)$, then $f(\bar a)\in \dcl(\nabla_{\infty} (\bar a), k)$ and so 
$\partial^m(f(\bar a))\in \dcl(\nabla_{N}(\bar a), k)$, $m\geq 0$, with $\partial^0(f(\bar a))=f(\bar a)$. We can express by an $\cL(k)$-formula that $y\in \dcl(\bar z)$ with $\vert y\vert=1$ and $\vert\bar z\vert=(N+1)\vert \bar x\vert$.
\par Now let $x$ vary in $X$ and let $0\leq m\leq N$ and so we have a collection of $\cL(k)$-formulas $\xi_m(x,\bar z)$ expressing that  $\xi_m(\nabla_N(\bar x),\partial^m(f(\bar x)))$ holds in $T(\partial)\cup\{\bar x\in X\}$.
Applying compactness, for each $0\leq m\leq N$,  there are finitely many such $\cL(k)$-formulas $\xi_{m,j}(\bar z, y)$, $j\in J(m)$.
\par By Lemma \ref{set}, there are finitely many $(\chi,\varphi)\in \cF(\nabla(Y))$, with $Y=\bigcup Y_{\chi,\varphi}$, $Y_{\chi,\varphi}$ an $\cL(k)$-definable set, such that $\nabla_{N}(X) = \bigcup_{(\chi,\varphi)\in \cF(\nabla(Y))} \{a\in Y_{\chi,\xi}: \nabla(a) = s(a)\}$, where $s:Y_{\chi,\varphi}\to \tau(Y)$ is an $\cL(k)$-definable section. 
\par For each $0\leq m\leq N$, $j\in J(m)$, let $\psi_{m,j}(\bar z, y)$ be the $\cL(k)$-formula $\xi_{m,j}(\bar z,y)\wedge \bar z\in Y\wedge (\forall \bar z\exists^{=1} y\; \xi_{m,j}(\bar z,y))$. 
So $\psi_{m,j}(K)$ is the graph of a partial function $\tilde F_j(m)$ whose domain is included in $Y$ such that
\begin{enumerate}
 \item $\forall x\in X\;(f(\bar x)=y \leftrightarrow \bigvee_{j\in J(0)} \tilde F_{j}(\nabla_{N}(\bar x))=y)$,
 \item for $m\geq 1$, for each $j\in J(m)$, $\psi_{m,j}(\nabla_{N}(\bar x), y)\rightarrow \;(y=\partial^m(\tilde F_{j}(\nabla_{N}(\bar x)))$, 
 \item for $m\geq 1$, $\bigwedge_{j\neq j'\in J(m)}\dom(\tilde F_{j}(m))\cap \dom(\tilde F_{j'}(m))=\emptyset$,
 
 \end{enumerate}
\par For each $j\in J(m)$, we define the $m$-prolongation $\tilde F_{j}[m](\bar z)$ to be the $m+1$-tuple of (partial) functions $(\tilde F_{j}(\bar z),\ldots,\tilde F_{j}(m)(\bar z))$; let  $F(m)(\bar z):=\bigcup_{j\in J(m)} \tilde F_{j}(m)$, with $F=F(0)$ and accordingly $F[m](\bar z):=\bigcup_{j\in J(m)} \tilde F_{j}[m]$, $N\geq m\geq 0$. 
\par Note that for each $0\leq m\leq N$, $\dom(F(m))$ is a large subset of $Y$ since any $\cL(k)$ generic type of $Y$ is realized by some differential tuple $\nabla_N(\bar a)$ with $\bar a\in X$ (see Lemma \ref{set}).\qed
\cor\label{function++} Let $X, \tilde X$ be finite-dimensional $\cL_{\partial}(k)$-sets in $(\U,\partial)$ included in some cartesian product of the field sort and suppose that $\dim_{\partial}(X)$ finite.
Let $f$ be an $\cL_{\partial}(k)$-definable function in $(\U,\partial)$ from $X\to \tilde X$.
Then for suitable $N$, there is an $\cL(k)$-function $F$ in $\U$ with $\nabla_{N}(X)\subset \dom(F)$ and 
such that for $x\in X$, 
\[f(x)=F(\nabla_{N}(x)).\]
\par Furthermore for $m\geq 1$, there are $\cL$-definable prolongations $F(m)$ of $F$ with domain a large subset of $\dom(F)$ and such that $\partial^m(f(x))=F(m)(\nabla_{N}(x))$.
\ecor
\pr Assume that $\tilde X\subset K^n$ and
  compose $f$ with the $n$ coordinate projections $\pi_{i}: K^n\to K$ and set $f_{i}:=\pi_{i}\circ f$, $1\leq i\leq n$ and we may assume that $\dom(f_{i})=\dom(f)$.
\par We proceed as in Lemma \ref{function+} to get a natural number $N$ and an $\cL(k)$-definable set $Y$ such that $\nabla_N(X)$ is a dense subset of $Y$ and we obtain $n$ $\cL(k)$-definable functions $F_i:Y\to K$ such that for $F:=(F_1,\ldots,F_n)$, we get $f(x)=F(\nabla_N(x))$. Further for each $m\geq 1$, wet get $m$ prolongations $F(m):=(F_{1}(m),\ldots, F_{n}(m))$ with $\partial^m(f(x))=F(m)(\nabla_N(x))$. Since later on we will have to compose $f$ with itself, we introduce the notation: $F[m]:=(F_{1}[m],\ldots, F_{n}[m])$. \qed

\section{Weil pre-group construction}\label{chunk}
In this section we recall the Weil pre-group construction. We fix a complete $\cL$-theory $T$ (expanding the theory of fields of characteristic $0$ with $\cL$ containing the language of rings $\cL_{rings}=\{+, -, \cdot, 0, 1\}$), possibly many-sorted and throughout this section $T$ will be a geometric theory on the field sort.

\rem For convenience of the reader (and future use), let us recall two well-known properties of generic elements, which are direct consequences of equation (\ref{dim}) in section 2.1 and the fact that if $a$ and $c$ are interalgebraic over $b$, then $\dim(a/b)=\dim(c/b)$. 
\medskip
\par Let $I, F$ be two $\cL(C)$-definable functions with $\dom(I)\subset X$ and $\dom(F)\subset X\times X$ and images in $X$ . 
\cl \label{cl1} Let $a\in X$ and suppose that $a$ is  generic in $X$ over $b$ with $b$ some tuple of elements in the field sort. Suppose that $a, I(a)\in \dom(I)$ and $I(I(a))=a$. Then $I(a)$ is again a generic element of $X$ over $b$. \qed
\ecl

\cl \label{cl2} Let $a$ be generic in $X$ over $b$, with $b\in \dom(I)$. Suppose that $(a,b)\in \dom(F)$ and $a=F(F(a,b), I(b))$. Then $F(a,b)$ is generic in $X$ over $b$. Similarly if $(b,a)\in \dom(F)$ and $a=F(I(b), F(b,a))$, then $F(b,a)$ is generic in $X$ over $b$. \qed
\ecl
\erem

\fct \label{def-large} {\rm \cite[Proposition 1.13, Remark 1.14]{P88}, \cite[Lemma 2.3]{Hrushovski-Pillay-groupslocalfields}} 

Let $X$ be an $\cL(B)$-definable set in $\U$, with $X\subset K^n$ and $B\subset K$ and let $\varphi( x, y)$ be an $\cL$-formula with $\vert y\vert=m$. Then $\{c\in K^m:\;\varphi(K, c)\cap X$ is large in $X\}$ is $B$-definable.

Moreover $\{c\in K^m:\;\varphi(K, c)\cap X$ is large in $X\}=\{c\in K^m:$ for every generic element $u$ of $X$ over $c$, $\U\models \varphi(u, c)\}$.
\efct

\subsection{The construction}

\par For the convenience of the reader, we recall below the definition of  a pre-group that we will be working with and give the proof that from a pre-group one obtains a definable group, often presented as an interpretable group. This construction was introduced by A. Weil \cite[Proposition 1]{W} and adapted in many different contexts by model-theorists from the late eighties on (see for instance \cite{D90}, \cite{E}). 
\par The second aspect in A. Weil construction, is not only to obtain a definable group but recover, from generic data, in his case, an algebraic group. An analogue of this aspect, when $T$ is either the theory of real-closed fields or the theory of $p$-adically closed fields is about constructing a Nash group which was done in \cite{P88}, \cite{P89} (for $\Q_{p}$ and extended for $p$-adically closed fields in the paper on Nash D-groups \cite{PPP2}).

\medskip

\par For the reader familiar with the Weil's construction, we give first a quick sketch, so one can skip the remainder of this section (and possibly later returns to it (for checking the notations) when reading section 4). We also compare what is done here and a recent account in \cite[section 2]{E}, since there the author works in the abstract context of a structure $\mathcal R$ endowed with a (definable, fibered) dimension function \cite[]{E}, but assuming that  $\mathcal R$ has elimination of imaginaries. He considers a definable subset $X\subset R^n$ endowed with an inverse function $i$ and a binary group operation $F$ defined on large subsets of $X$, respectively $X\times X$ and satisfying what he calls {\it Axioms} and are similar to the pre-group properties listed in Lemma \ref{lem:pregroup}. He calls this data $(X, i, F)$ {\it a definable group chunk} \cite[Definition 2.1]{E} and follows the now classical strategy. On the set of $2$-tuples of elements of $X$, he defines an equivalence relation on the corresponding set of germs on $X$ and obtains a group structure on this set of equivalence classes. Since $\mathcal R$ has elimination of imaginaries, he finally obtains a definable group from the above data \cite[Theorem 2.1]{E}.
\par Here following closely \cite{PPP}, we define a pre-group in terms of properties of its generic elements, then from that we obtain the so-called pre-group properties. Then indexed by a large subset $Z$ of $X$, we consider the corresponding family of germs $f_{a}$ on $X$ (also called generic permutations), $a\in Z$, and similarly we show that the (definable) set $\{f_a\circ f_b: a, b\in Z\}$ is a subgroup of the group of generic permutations on $X$ (this last group is a priori not definable). (Note that a germ is an equivalence class but this construction does not depend on the representative we choose in the equivalence class; further the map sending $a\in Z$ to $f_{a}$ is an injection.)
\par Then in that last step, as in \cite{PPP}, we differ slightly from what was done in \cite{E}, since $X$ is a definable subset, it has finite dimension, say $d$ and on the set of $2d+1$ translates of $Z$, we define a group
operation which coincides with initial operation $F$ on generic independent elements of $Z$. Note that the conclusion is weaker than the one obtained in \cite{E}, but we don't assume elimination of imaginaries and we have a weaker version of  Axiom (4) in \cite[Definition 2.1]{E}). 
\dfn\label{def:pregroup} 
Let $\tilde G$ be an $\cL(A)$-definable set in $\U$, included in the field sort with partially defined $\cL(A)$-definable functions $F$ and $F_{-1}$.  Then $(\tilde G, F, F_{-1})$ is called an $\cL(A)$-definable pre-group, if the following properties hold:
\begin{enumerate}[label={(P\arabic*)}]
\item for generic elements $a\in \tilde G$, we have $a\in \dom(F_{-1})$, $F_{-1}(a)$ generic in $\tilde G$ and $F_{-1}(F_{-1}(a))=a$,
\item for generic elements $(a, b)\in \tilde G^2$, we have $(a, b)\in \dom(F)$ and $F(a,b)$ generic in $\tilde G$,
\item for generic and pairwise independent elements $a, b , c\in \tilde G$, $F(a, b)$ is independent of $c$ over $A$, $a$ is independent of $F(b,c)$ over $A$ and  $F(a,F(b, c))=F(F(a,b),c)$,
\item for generic elements $(a,b)\in \tilde G^2$, we have  $F(F_{-1}(a),F(a,b))=b$ and $F(F(a,b),F_{-1}(b))=a$.
\end{enumerate}
\edfn
\par First we use that in a (sufficiently saturated) geometric structure $\U$, when an  $\cL(A)$-formula $\varphi$ is satisfied by all $\cL(A)$-generic elements in an $\cL(A)$-definable subset $Y\subset K^n$, then $\varphi(K)$ is a large subset of $Y$ (see Definition \ref{dlarge}). We gather in lemma below what we call the {\it pre-group properties} of large subsets.
\lem \label{lem:pregroup} 
Let $(\tilde G, F, F_{-1})$ be an $\cL(A)$-definable pre-group. Then $\dom(F_{-1})$ is a large $\cL(A)$-definable subset of $\tilde G$, $Z_{2}:=\dom(F)$ is a large $\cL(A)$-definable subset of $\tilde G\times \tilde G$, and there is a large $\cL(A)$-definable subset $Z_{3}$ of $\dom(F_{-1})$ with the following properties:  
\begin{enumerate}[label={(g\arabic*)}]
\item for $a\in Z_3$ we have  $F_{-1}(a)\in Z_3$ and $F_{-1}(F_{-1}(a))=a$,
\item for $a\in Z_3$, let 
$$X_a:=\{y\in \tilde G: (a,y)\in Z_2\;\&\;(F_{-1}(a), F(a,y))\in Z_2\,\&\,F(F_{-1}(a),F(a,y))=y\},$$
then $X_a$ is an $\cL(A)$-definable large subset of  $\tilde G$,
\item for $a\in Z_3$, let $$X_a^r:=\{y\in \tilde G: (y,a)\in Z_2\;\&\;(F(y,a),F_{-1}(a))\in Z_2\;\&\;F(F(y,a),F_{-1}(a))=y\},$$ then $X_a^r$ is an $\cL(A)$-definable large subset of $\tilde G$,
\item for $x\in Z_{3}$, for any generic element $a\in \tilde G$ independent over $x$, we have $x\in X_a$,
\item for $x\in Z_{3}$, for any generic element $a\in \tilde G$ independent over $x$, we have  $x\in X_a^r$,
\item for any $x, y\in Z_{3}$, the set $\{a\in \tilde G:\; (a,x), (a,y)\in \dom(F) \wedge ((F(a,x)=F(a,y)\rightarrow x=y)\;\&\; (F(x,a)=F(y,a)\rightarrow x=y))\}$ is an $\cL(A)$-definable large subset of $\tilde G$,
\item for any $z\in Z_{3}$, there is an $\cL(A)$-definable large subset $Y(z)$ of $Z_2$ such that if $(x,y)\in Y(z)$,  then $(y,z)\in \dom(F)$, $(x,y)\in \dom(F)$,  $(x,F(y,z))\in \dom(F)$, $(F(x,y),z)\in \dom(F)$ and $F(x,F(y,z))=F(F(x,y),z)$.
\end{enumerate}
\elem
\pr Even if this is well-known in order to give a flavour of how these properties follow from the properties of the generic elements, we will indicate why some of these properties hold.
\par First the properties that $\dom(F_{-1})$ is a large $\cL(A)$-definable subset of $\tilde G$, and $\dom(F)$ is a large $\cL(A)$-definable subset of $\tilde G\times \tilde G$, follow immediately from (P1) and (P2) in the pre-group definition.
\par Then we will put requirements on an $\cL(A)$-formula $\varphi(x)$ by stages and $Z_{3}$ will be defined as $\varphi(K)$.
\par (i) First let $\varphi(x)$ be an  $\cL(A)$-formula expressing that $x\in \dom(F_{-1})$, that $F_{-1}(x)\in \dom(F_{-1})$ and $F_{-1}((F_{-1}(x))=x$.
Then $\varphi(K)$ is an $\cL(A)$-definable large subset of $\tilde G$ since by (P1), if $x$ is generic, then $\varphi(x)$ holds.
\par (ii) By Fact \ref{def-large}, $\{x\in \tilde G\colon \forall y$ $y$ generic over $x, y\in X_{x}\}$ is an $\cL(A)$-definable subset. So we may assume that $\varphi(x)$ also expresses that $X_{x}$ is large in $\tilde G$. By (P4) for generic $x$, if $y$ is generic and independent from $x$, then $y\in X_{x}$. So $\varphi(x)$ holds for any generic element $x$. 
 \par By a similar reasoning we may assume that $\varphi(x)$ also expresses that $X_{x}^r$ is large in $\tilde G$ and still get that $\varphi(K)$ is an $\cL(A)$-definable large subset of $\tilde G$.
\par (iii) Again by Fact \ref{def-large}, $\{z\in \tilde G\colon Y(z)$ is a large subset of $\dom(F)\}$ is an $\cL(A)$-definable subset. So we may suppose that $\varphi(z)$ expresses in addition that $Y(z)$ is a large subset of $\dom(F)$. 
Now let us show that if $c\in \tilde G$ is generic then $\varphi(c)$ holds.
Let $a, b\in \tilde G$ with $a, b$ generic independent over $c$. Then by (P2), $F(a, b)$ is generic and independent of $c$ over $A$ and $F(b, c)$ is generic and independent over $a$ (and so $a$ is independent of $F(b,c)$) over $A$. So by (P3),  $F(a,F(b, c))=F(F(a,b),c)$. So $Y(c)$ is a large $\cL(A)$-definable subset of $\tilde G\times \tilde G$ and $\varphi(K)$ is a large subset of $\tilde G$.
\medskip
\par (iv) Again by Fact \ref{def-large}, $\{x\in \tilde G\colon \forall a\in \tilde G$, if  $a$ generic over $x$, then $x\in X_a\}$ is an $\cL(A)$-definable subset. So we may assume that $\varphi(x)$ also expresses that $\forall a\in \tilde G$, if  $a$ generic over $x$, then $x\in X_a$. Let us show that if $x$ is generic in $\tilde G$, then $\varphi(x)$ holds.
Since $a$ is generic, $X_{a}$ is large and so $x\in X_{a}$. 
\par We have a similar statement with $X_{a}^r$. Let $Z_{3}:=\varphi(K)$ and note that $Z_{3}\subset \tilde G$.
\par Now let us show that for any $x, y\in Z_{3}$, the set $\{a\in \tilde G: ((F(a,x)=F(a,y)\rightarrow x=y)\;\&\; (F(x,a)=F(y,a)\rightarrow x=y))\}$ is a large $\cL(A)$-definable subset of $\tilde G$. Let $a\in \tilde G$ be any generic element independent over $x, y$. By construction of $\varphi$,  $x, y\in X_{a}$ ($X_{a}$ is large). So by the definition of $X_{a}$, $F(F_{-1}(a),F(a,x))=x$ and $F(F_{-1}(a), F(a,y))=y$.
Since $x, y\in X_{a}^r$, we also have the second part of the statement.
\qed
\bigskip
\dfn\label{gen-inj} Let $X$ be an $\cL(A)$-definable set in $\U$ with $X$ included in the field sort and let $f: X\to X$ be a partially defined definable map. Then we say that $f$ is a generic $\cL(A)$-definable permutation of $X$ if the domain and the image of $f$ are large $\cL(A)$-definable subsets of $X$ and on its domain, $f$ is injective. 
Let $f_1, f_2$ be generic $\cL(A)$-definable permutations of $X$, then $f_1, f_2$ have the same germ on $X$, denoted by $f_1\sim  f_2,$ if $f_1, f_2$ agree on a large $\cL(A)$-definable subset of $X$. We will denote the equivalence class containing a generic $\cL(A)$-definable permutation $f$ of $X$ by $\f$ and call this equivalence class a germ on $X$. Note that $\sim$ is a congruence for the composition operation, namely if $f_1\sim g_1$ and $f_2\sim g_2$, then $(f_1\circ f_2)\sim(g_1\circ g_2)$, with $f_{1}, f_{2}, g_1, g_2$ generic $\cL(A)$-definable permutations of $X$. 
\par Denote by $Perm_g(X)$ the group of all germs on $X$ endowed with composition.
\edfn
\nota Let $(\tilde G, F, F_{-1})$ be a pre-group and using Lemma \ref{lem:pregroup} together with its notations, for $c\in Z_{3}$, we associate the map (acting on a large subset of $\tilde G$): $$f_{c}:X_{c}\to \tilde G: x\mapsto F(c,x).$$
\enota
\par Note that if for $a, b\in Z_{3}$, $f_a \sim f_b$, namely on a large subset of $\tilde G$, $f_{a}$ and $f_{b}$ coincide, then $a=b$ (see Claim \ref{cl-inj} in the proof of Proposition \ref{prop:chunk}). Furthermore, suppose that $Z_{3}$ is defined by an $\cL$-formula of the form $\varphi(x,c)$ with $c$ a tuple from $A$ and $\vert x\vert=n$. Let $n\geq m=\dim(\tilde Z_{3})$. By the property of definability of dimension in geometric theories \cite[Lemma 2.3 (ii)]{Hrushovski-Pillay-groupslocalfields}, there is an $\cL$-formula $\theta(y)$ such that $\dim(Z_{3})=m$ if and only if $\theta(c)$ holds. Now let $a, b\in Z_{3}$ and express that $Z(a,b):=\{x\in Z_{3}: F(a,x)=F(b,x)\}$ is a large subset of $Z_{3}$, namely $\dim(Z(a,b))=m$ and $\dim(Z_{3}\setminus Z(a,b))<m$. Again, by definability of dimension in geometric theories, there exits an $\cL(A)$-formula $\xi(u,v)$ such that $\xi(a,b)$ holds if and only if $f_{a}\sim f_{b}$.
 \medskip
\par For future use in this paper, we detail the main steps of a proof of the following proposition which explain how to construct a definable group from generic data. 
\prop\label{prop:chunk} Let $(\tilde G, F, F_{-1})$ be an $\cL(A)$-definable pre-group with $\dim(\tilde G)=d$.
 Then there is an $\cL(A)$-definable large subset $Z$ of $\tilde G$ such that $\{\mathfrak f_{c}\colon c\in Z\}\subset Perm_{g}(\tilde G)$
 and the set $H:=\{\f_{b}\circ \f_{a}\colon a, b\in Z\}$ is the domain of a subgroup of $Perm_g(\tilde G)$.
 \par There are finitely many parameters $a_1,\ldots,a_{2d+1}$ and an $\cL(A\cup\{a_1,\ldots,a_{2d+1}\})$-definable group $(G,\times)$ definably isomorphic to $(H,\circ)$. Moreover there is an $\cL(A\cap\{a_1,\ldots,a_{2d+1\}})$-definable injection $\varsigma$ of $Z$ into $G$ such that
for $a, b\in Z$ with $a$ generic over $b$, we have that $\varsigma(F(a,b))=\varsigma(a)\times\varsigma(b)$.  
\eprop
\pr We will refer the pre-group properties stated in Lemma \ref{lem:pregroup} by their numbering (g1) up to (g7) and we use the notations of the preceding lemma. For ease of notation, since all data is over a set $A$ of parameters, we use {\it independent} to mean $\acl$-independent over $A$ and we replace $\dim(\;/A)$ by $\dim(\;)$.
\par First note that $Z_3$ acts "generically" on $\tilde G$ by left multiplication: $x\mapsto F(a,x)$ with $a\in Z_{3}$. 
Indeed, by (g2), $X_{a}$ is a large subset of $\tilde G$, by (g1), (g2), $F_{-1}(a)\in Z_{3}$, the image of $f_{a}$ is a large subset and again by (g2) $f_{a}$ is injective on $X_{a}$. So $f_a$ is a generic bijection and its inverse is given by $f_{F_{-1}(a)}$. Moreover, as noted before, $\sim$ is a congruence relation for composition. So we have shown:
\cl\label{cl-1} The set $\{\f_a \colon a\in Z_3\}\subset Perm_g(\tilde G)$.\qed
\ecl
\cl\label{cl-inj} The map from $Z_{3}$ to $Perm_{g}(\tilde G)$ sending $a\mapsto \f_{a}$  is injective. 
\ecl
\pr For any  generic element $x\in \tilde G$ independent from $a,\; b$, we have that $a, b\in X_{x}^r$ by (g5) and by assumption $F(a,x)=F(b,x)$ (by Fact \ref{def-large}). So again by (g5), $a=b$ since
 $a=F(F(a,x),F_{-1}(x))=F(F(b,x),F_{-1}(x))=b$. In other words, for $a\neq b\in Z_3$, $f_{a} \not\sim f_{b}$.\qed
\medskip
\par Let us observe that generically, composition of maps $f_{a}\circ f_{b}$ corresponds to the map $f_{F(a,b)}$.
\cl \label{cl0} There is a large $\cL(A)$-definable subset $Z$ of $Z_{3}$ such that if  $a\in Z$,  we have: for every independent generic element $b$ over $a$, for every $x\in Z_3$ for which the following operations are defined, $F(F(b,a),x)=F(b,F(a,x))$. 
\par So for $a\in Z$, for every independent generic element $b$ over $a$, we have $f_b\circ f_a\sim f_{F(b,a)}$.\qed
\ecl

\medskip
\par From now on we will replace $Z_{3}$ by $Z$ as in Claim \ref{cl0} and then $Z$ by $Z\cap F_{-1}(Z)$. Note that $Z$ is still an $\cL(A)$-definable large subset of $\tilde G$.
\bigskip

\par For convenience of the reader, even if the argument is well-known, we give the proof of the following claim since it will be used later in the proof of our main (see Claim \ref{gd}).
\cl \label{cl3} Let $H:=\{\f_a\circ \f_b\colon a, b\in Z\}$. Then $H$ endowed with composition is a subgroup of $Perm_g(\tilde G)$. 
\ecl
\pr 
We need to show that for $a, b, c, d\in Z$, there are $r, s\in Z$ such that 
\[f_a \circ f_b\circ  (f_{d}\circ f_{c})^{-1}\sim f_r \circ f_s.
\]
 Since $(f_d)^{-1}\sim f_{F_{-1}(d)}$, we have 
\[f_a\circ f_b\circ (f_{d}\circ f_{c})^{-1}\sim f_a\circ f_b\circ f_{F_{-1}(c)}\circ f_{F_{-1}(d)}.
\]
Choose $e\in Z$ generic over $a, b, c, d$. By (P1), $F_{-1}(e)$ is generic over $a, b, c, d$ and by Claim \ref{cl2}, $F(e,F_{-1}(c)), F(b,F_{-1}(e))$ are generic over $a, b, c, d$.
By Claim \ref{cl0}, 
\begin{align*}
f_b \circ f_{F_{-1}(e)}\sim& f_{F(b,F_{-1}(e))}, \\
f_e \circ f_c\sim& f_{F(e,c)}, \\
f_a\circ f_{F(b,F_{-1}(e))}\sim& f_{F(a,F(b,F_{-1}(e)))},\\ 
f_{F(e,F_{-1}(c))}\circ f_{F_{-1}(d)}\sim& f_{F(F(e,F_{-1}(c)),F_{-1}(d))}.
\end{align*}
So we get:
\begin{align*}
f_a \circ f_b \circ f_{F_{-1}(c)}\circ f_{F_{-1}(d)}\sim& f_a \circ f_b\circ f_{F_{-1}(e)}\circ f_e\circ f_{F_{-1}(c)}\circ f_{F_{-1}(d)}\\
\sim& f_a \circ f_{F(b,F_{-1}(e))} \circ f_{F(e,F_{-1}(c))}\circ f_{F_{-1}(d)}\\
\sim& f_{F(a,F(b,F_{-1}(e)))} \circ f_{F(F(e,F_{-1}(c)),F_{-1}(d))}.
\end{align*}
\qed
\medskip
\par The following calculation may be found in \cite{W}. 
\cl Let $d=\dim(\tilde G)$ and let $a_1,\ldots,a_{d+1}$ be generic independent elements in $\tilde G$. Let $b\in \tilde G$. Then there $1\leq i\leq d+1$ such that $a_i$ is independent from $b$.\qed
\ecl

\cl \label{cl:dim1} Let $d=\dim(\tilde G)$ and let $a_1,\ldots,a_{2d+1}$ be generic independent elements in $\tilde G$. Let $a, b\in \tilde G$. Then there $1\leq i\leq d+1$ such that $a_i$ is independent from $a, b$. \qed
\ecl 

\cl \label{cl:dim2} Let $a_1,\ldots,a_{2d+1}$ be generic independent elements in $\tilde G$. Then the set $\{\f_b\colon b\in Z\}\subseteq\bigcup_{i=1}^{2d+1} \{\f_{F_{-1}(a_i)}\circ \f_{d}\colon d\in Z\}.$ \qed
\ecl
\medskip
\cl \label{cl:dim3} Let $a_1,\ldots,a_{2d+1}$ be any $2d+1$ generic independent elements in $\tilde G$. Let $G:=\bigcup_{i=1}^{2d+1} \{\f_{F_{-1}(a_i)}\circ \f_{d}\colon d\in Z\}$.
Then there is a bijection between  $H$ and $G$, definably given over the parameters $A\cup\{a_1,\ldots,a_{2d+1}\}$. Further $G$ can be endowed with an $\cL(A\cup\{a_1,\ldots,a_{2d+1}\})$-definable group operation $\times$ such that $(G,\times)\cong (H,\circ)$.
\par In addition, there is an injection $\varsigma$ of $Z$ into $G$ such that for $a, b\in Z$ with $a$ generic over $b$, we have that $\varsigma(F(a,b))=\varsigma(a)\times\varsigma(b)$.  
\ecl
\pr By Claim \ref{cl:dim2}, $\{\f_b\colon b\in Z\}\subseteq\bigcup_{i=1}^{2d+1} \{\f_{F_{-1}(a_i)}\circ \f_{d}\colon d\in Z\}.$ Note that by Claim \ref{cl-inj}, the map sending $Z\to \{\f_b: b\in Z\}: b\mapsto \f_b$ is injective. So we have an injective map $\varsigma$ from $Z$ to $G$.
\par Now let us show that $H=G$. Let $a, b\in Z$ and choose $a_i$ generic over $a,\,b$, for some $1\leq i\leq 2d+1$. Let us check that $f_a\circ f_b\sim f_{F_{-1}(a_{i})}\circ f_d$ for some $d\in Z$. Write  $a$ as $F(F_{-1}(a_i),F(a_i, a))$ (Claim \ref{cl:dim1}). Then $f_{F(a_i, a)}\circ f_b\sim f_{F(F(a_i, a),b)}$, since $F(a_i,a)$ is generic over $b$ (by (g5), $a\in X_{a_{i}}^r$). Since $f_a\sim f_{F_{-1}(a_i)}\circ f_{F(a_i,a)}$, we get $f_a\circ f_b\sim (f_{F_{-1}(a_i)} \circ f_{F(a_i,a)})\circ f_b\sim f_{F_{-1}(a_i)} \circ (f_{F(a_i,a)}\circ f_b)\sim f_{F_{-1}(a_i)}\circ f_{F(F(a_i,a),b)}$.
\par So $\bigcup_{i=1}^{2d+1} \{\f_{F_{-1}(a_i)}\circ \f_{d}\colon d\in Z\}=\{\f_a \circ \f_b: a, b\in Z\}.$ 
\par Set $Z^{(0)}:=Z$. Let us identify $G$ with the (disjoint) union: $\bigsqcup_{i=1}^{2d+1} \{\f_{F_{-1}(a_i)} \f_{d}\colon d\in Z^{(i)}\}$, where 
 $Z^{(i)}$, $1\leq i\leq 2d$, is defined as $\{d\in Z\colon f_{F_{-1}(a_{i})}\circ f_{d}\not\sim f_{F_{-1}(a_{j})}\circ f_{e}$ for some $e\in Z^{(j)}$ and some $0\leq j<i \}$.
 
\par  So we define a map $\varrho$ from $H$ into $G$ by choosing $i$ minimum such that $f_a\circ f_b\sim f_{F_{-1}(a_i)}\circ f_d$, for some $d\in Z$. 
\par By Claim \ref{cl3}, $H$ is a subgroup of $Perm_g(\tilde G)$. Then we define $\times$ on $G$ as follows: $(\f_{F_{-1}(a_i)}\circ \f_a)\times \f_{F_{-1}(a_j)}\circ \f_b):= \f_{F_{-1}(a_k)}\circ \f_c$, with $k$ minimum such that $f_{F_{-1}(a_k)}\circ f_c\sim (f_{F_{-1}(a_i)}\circ f_a)\circ (f_{F_{-1}(a_j)}\circ f_b)$.
\par This is well-defined since by Claim \ref{cl3}, there is $u, v\in Z$ such that $(f_{F_{-1}(a_i)}\circ f_{a})\circ (f_{F_{-1}(a_j)}\circ f_{b})\sim f_{u}\circ f_{v}$. So $\f_u\circ \f_v=\f_{F_{-1}(a_k)}\circ \f_c$ with $k$ chosen minimum such.
\par The elements $u, v$ are defined as $u=F(F_{-1}(a_i),F(a,F_{-1}(e))$, $v=F(F(e,F_{-1}(a_j)),b)$ for some $e$ generic over $a_i, a, a_j, b$. 
\par Now take $a, b\in Z$ with $a$ generic over $b$,  so $F(a,b)$ is defined and belongs to $Z$. By Claim \ref{cl0}, $f_{a}\circ f_{b}\sim f_{F(a,b)}$. So, $\varsigma(a)\times\varsigma(b)=\varsigma(F(a,b))$.
\qed

\rem \label{ident} Keeping the notations of the above proposition, let us note the following. 
There is a bijection from  $\bigsqcup_{i=1}^{2d+1} \{\f_{F_{-1}(a_{i})}\circ \f_{a}\colon a\in Z^{(i)}\}$ to $\bigsqcup_{i=1}^{2d+1} \{F_{-1}(a_i)\}\times Z^{(i)}$. (By definition, $\f_{F_{-1}(a_{i})}\circ \f_{a}= \f_{F_{-1}(a_{j})}\circ \f_{b}$ with $a\in Z^{(i)}$, $b\in Z^{(j)}$, implies that $i=j$ and by composing with $\f_{a_i}$, we get that $\f_{a}=\f_{b}$ which implies by Claim \ref{cl-inj}, that $a=b$.)
\erem
 \section{$\cL_\partial$-definable groups in models of $T_\partial$}\label{sec:groups}
\par We work under the following assumptions on $T(\partial)$ and $T_{\partial}$ described in sections \ref{sec:diff}, \ref{sec:defsets}, namely:
\begin{enumerate}
\item $T$ is geometric, model-complete,
\item $T(\partial)$ has a model-companion $T_\partial$ and $T(\partial)$ satisfies the compatibility conditions (C1) up to (C3),
\item $T_\partial$ satisfies the compatibility conditions (A1) till (A2). 
\end{enumerate}
\par Let $(\U,\partial)$ be a sufficiently saturated model of $T_{\partial}$ and let $k$ be a small differential subfield of $K$.

\par In this section we describe finite-dimensional $\cL_{\partial}(k)$-definable groups in $(\U,\partial)$ included in the field sort $K$.
In section \ref{sec:defsets}, we showed how to associate with a finite-dimensional  $\cL_\partial(k)$ definable set $X$ in $K$ an $\cL(k)$-definable set $Y$ (possibly living in a different cartesian product of $K$ and containing $\nabla_N(X)$) with the property that any $\cL(k)$-generic type in $Y$ is realized by a differential tuple in $Y$ (see Lemma \ref{set}). Further  we showed how to associate with an $\cL_\partial(k)$-definable function $f$ with domain $X$, an $\cL(k)$-definable function $F$ with domain a large subset of $Y$ (see Lemma \ref{function+} and Corollary \ref{function++}).
\subsection{Finite-dimensional groups}\label{sec:fd}

\prop \label{largeV}  Let $\cG:=(\G,f_{\times} ,f_{-1},e )$ be an $\cL_{\partial}(k)$-definable finite-dimensional group in $(\U,\partial)$ and $\G\subset K^n$. Then, there are $N\in \N$ and
 \begin{itemize}
 \item an $\cL(k)$-definable set $Y$ with $\nabla_N(\Gamma)\subset Y$ with the property that any $\cL(k)$-generic type in $Y$ is realized by a tuple $\nabla_N(a)$ with $a\in \Gamma$,
 \item an $\cL(k)$-definable function $F_{-1}$ on a $\cL(k)$-definable large subset $Z\subset Y$ containing $\nabla_N(\Gamma)$ such that
for every $x\in \G$ 
\[	
f_{-1}(x)=F_{-1}(\nabla_N(x)),
\]

\item an $\cL(k)$-definable function $F_{\times}$
on an $\cL(k)$-definable large subset $D\subset Y\times Y$ containing $\nabla_N(\Gamma)\times \nabla_N(\Gamma)$ such that for every $x_{1}, x_{2}\in \G$ 
\[
f_{\times}(x_{1},x_{2})=F_{\times}(\nabla_N(x_{1}),\nabla_N(x_{2})),
\]
\item 
$\nabla_{N}(f_{-1}(x))=F_{-1}[N](\nabla_N(x))\;{\rm and\;} \nabla_{N}(f_{\times}(x_{1},x_{2}))=F_{\times}[N](\nabla_N(x_{1}),\nabla_N(x_{2})).$

\end{itemize}
Moreover $(Y, F_{-1}[N], F_\times[N])$ is an $\cL(k)$-definable pre-group.
\eprop
\pr The first part of the statement follows from Lemmas \ref{set},  \ref{function+} and Corollary \ref{function++}, with $X = \Gamma$.
\par Now to check that $(Y, F_{-1}[N], F_\times[N])$ satisfies the pre-group properties, is rather immediate since properties (P1) till (P4) have to be verified for generic elements of $Y$. But any $\cL(k)$-generic type in $Y$ is realized by a differential tuple $\nabla_N(a)$ with $a\in \Gamma$, and by the first part the operations $F_{\times}$ and $F_{-1}$ coincide with the group operations in $\Gamma$, when applied to differential tuples. Then we have to use their prolongations to compose them in $Y$. \qed
\thm \label{thm:def}  Let $\cG:=(\G,\times ,{\;}^{-1})$ be an $\cL_{\partial}(k)$-definable finite-dimensional group  with $\G\subset K^n$. Then there is an integer $N$ and an $\cL(k)$-definable group $G$ over some finitely many additional parameters $\nabla_{N}(\bar u)$ with $\bar u\in \G$ and an $\cL_{\partial}(k)$-definable group embedding of $\nabla_{N}(\G)$ into $G$ such that:
\begin{enumerate}
\item any $\cL(k)$-generic type in $G$ is realized by a tuple $\nabla_{N}(a)$ with $a\in \Gamma$,
\item there are finitely many $\cL(k)$-definable sets $Y_{\ell}$ in $\U$, $\ell\in L$, with $Y_\ell \subset K^{n(N+1)}$, such that for $Y:=\bigcup_{\ell\in L} Y_\ell$, there is an $\cL(k)$-definable section $s:Y\to \tau(Y)$, such that
$\nabla_{N}(\G) = \bigcup_{\ell\in L} \{a\in Y_{\ell}: \nabla(a) = s(a)\}$.
\end{enumerate}

\ethm
\pr (1) We apply Proposition \ref{prop:chunk} to $(Y, F_\times, F_{-1})$ as defined in Proposition \ref{largeV}, with in particular $\nabla_{N}(\Gamma)\subset Y$. For ease of notation note that in the statement we denoted the group operations in $\G$ by $\times$ and ${\;}^{-1}$ (instead of $f_\times$ and $f_{-1}$ used in the proposition above).
\par In Proposition  \ref{prop:chunk}, we constructed an $\cL(k)$-definable group $G$ isomorphic to the subgroup $H$ of $Perm_g(Y)$ with domain $\{\f_a\circ \f_b: a, b\in Z\}$ with $Z\subset \dom(F_{-1})$ a large $\cL(k)$-definable subset of $Y$.  Moreover inspecting the construction of $Z$, we note that $\nabla_{N}(\Gamma)\subset Z$ (indeed we possibly restricted $\dom(F_{-1})$ (to a large $\cL(k)$-definable subset $Z$ of $Y$) in order to ensure that generically the group properties hold on $Z$, using that these properties hold on $\nabla_{N}(\Gamma)$). Let $d:=\dim(Y)$. 
 Using $2d+1$ $\cL(k)$-generic independent elements $(c_{i})_{i=1}^{2d+1}$ of $Y$, we showed that the domain of $G$ is equal to $\bigcup_{i=1}^{2d+1} \{\f_{c_{i}}\circ \f_{d}: d\in Z\}$. 
 Since any $\cL(k)$-generic type in $Y$ is realized by a tuple $\nabla_N(a)$ with $a\in \Gamma$, these generic elements can be taken as differential tuples coming from $\G$. 
 Let $a_{1},\ldots, a_{2d+1}$ be $2d+1$ $\cL_{\partial}(k)$-generic independent elements of $\Gamma$. Then $\nabla_{N}(a_{1}),\ldots, \nabla_{N}(a_{2d+1})$ are $\cL$-independent generic elements in $Z$ (over $k$).
\par  Let $a, b\in \Gamma$ be two $\cL_\partial(k)$-generic independent elements, then $\nabla_N(a),\nabla_N(b)$ are $\cL(k)$-generic independent elements of $Z$. 
Let $f_{\nabla_N(a)}, f_{\nabla_N(b)}$ be the associated generic permutations acting on $Y$. By Proposition \ref{prop:chunk} (Claim \ref{cl0}), $f_{\nabla_N(a)}\circ f_{\nabla_N(b)}\sim f_{F_{\times}[N](\nabla_N(a),\nabla_{N}(b))}$. By Proposition \ref{largeV}, $F_{\times}[N](\nabla_N(a),\nabla_{N}(b))=\nabla_{N}(a\times b)$. Further $\nabla_N(a^{-1})\in Z$, it is $\cL(k)$-generic and by Proposition \ref{largeV},  $\nabla_N(a^{-1})=F_{-1}[N](\nabla_{N}(a)).$
\par Now $\cG$ is an $\cL_{\partial}(k)$-definable group, so any element $b\in\Gamma$ is the product of two $\cL_{\partial}(k)$-generic elements, $b_{1}, b_{2}:=b_{1}^{-1}\times b$ (with $b_{1}$ $\cL_{\partial}(k)$-generic over $b$).
So $\nabla_N(b_1), \nabla_N(b_2)\in Z$ and we send $\nabla_N(b)$ to $\f_{\nabla_N(b_1)}\circ \f_{\nabla_N(b_2)}$. 
\cl \label{wd} The map $\rho$ from $\Gamma\to Perm_{g}(Y): b\mapsto f_{\nabla_N(b_1)}\circ f_{\nabla_N(b_2)}$ is well-defined, where $b=b_{1}\times b_{2}$ and $b_{1}, b_{2}$ are two $\cL_{\partial}(k)$-generic elements in $\G$.
\ecl
\prcl Let $b_{1}, b_{2}, c_{1}, c_{2}$ be $\cL_{\partial}(k)$-generic elements in $\G$ such that $b=b_{1}\times b_{2}= c_{1}\times c_{2}$, we have to show that $\f_{\nabla_{N}(c_{1})}\circ \f_{\nabla_{N}(c_{2})}= \f_{\nabla_N(b_1)}\circ \f_{\nabla_N(b_2)}$ in $Perm_{g}(Y)$.
\newline So we choose $\nabla_{N}(x)\in \nabla_{N}(\Gamma)$ generic over $\nabla_{N}(b_{1}), \nabla_{N}(b_{2}), \nabla_{N}(c_{1}), \nabla_{N}(c_{2})$.
\newline We have  $f_{\nabla_{N}(b_{2})}(\nabla_{N}(x))=F_{\times}[N](\nabla_{N}(b_{2}), \nabla_{N}(x))=\nabla_{N}(b_{1}\times b_{2})$ since $\nabla_{N}(x)\in X_{\nabla_{N}(b_{2})}$ and by Proposition \ref{largeV}. Then, 
\begin{align*}
f_{\nabla_N(b_1)}\circ f_{\nabla_N(b_2)}(\nabla_{N}(x))=&f_{\nabla_{N}(b_{1})}(\nabla_{N}(b_{2}\times x))\\
=&\nabla_{N}(b_{1}\times (b_{2} \times x))\\
=&\nabla_{N}((b_{1} \times b_{2})\times x)\\
=&\nabla_{N}((c_{1}\times c_{2})\times x)\\
=&f_{\nabla_{N}(c_{1})}\circ f_{\nabla_{N}(c_{2})}(\nabla_{N}(x)).
\end{align*}
\qed
\medskip
\newline So we send $\Gamma$ into $G$ by first sending $\Gamma$ to $H$ and then use the isomorphism between $H$ and $G$. Let $b\in \Gamma$, send it to $\f_{\nabla_N(b_1)}\circ \f_{\nabla_N(b_2)}$.
\cl \label{gd} There is an $\cL_{\partial}(k)$-definable group embedding from $\cG$ to $H$. 
\ecl
\prcl We use the map $\rho$ defined in Claim \ref{wd}.  
So given $b, c\in \Gamma$, with $b=b_{1}\times b_{2}$, $c=c_{1}\times c_{2}$ and $b\times c=d_{1}\times d_{2}$, where $b_{1}, b_{2}, c_{1}, c_{2}, d_{1}, d_{2}$ are $\cL_{\partial}(k)$-generic elements, we have to show that $$(f_{\nabla_{N}(b_{1})}\circ f_{\nabla_{N}(b_{2})})\circ (f_{\nabla_{N}(c_{1})}\circ f_{\nabla_{N}(c_{2})})\sim f_{\nabla_{N}(d_{1})}\circ f_{\nabla_{N}(d_{2})}.$$
As in Claim \ref{cl3}, we choose $e$ generic over $c_{1}, c_{2}, b_{1}, b_{2}$ and we write $(f_{\nabla_{N}(b_{1})}\circ f_{\nabla_{N}(b_{2})})\circ (f_{\nabla_{N}(c_{1})}\circ f_{\nabla_{N}(c_{2})})$ as $f_{F_{\times}[N](\nabla_{N}(b_{1}),F_{\times}[N](\nabla_{N}(b_{2}), \nabla_{N}(e^{-1})))}\circ f_{F_{\times }[N](F_{\times}[N](\nabla_{N}(e) ,\nabla_{N}(c_{1})), \nabla_{N}(c_{2}))}$.
Then by the above claim,  its action on $\nabla_{N}(x)$ with $x\in Z$ $\cL_{\partial}(k)$-generic independent over $e, b_{1},b_{2},c_{1},c_{2}, d_{1},d_{2}$ is equal to $\nabla_{N}((b_{1}\times b_{2}\times e) \times ((e^{-1}\times c_{1})\times c_{2})\times x)$ and it is
the same as the action of $f_{\nabla_{N}(d_{1})}\circ f_{\nabla_{N}(d_{2})}$ on $\nabla_{N}(x)$.
\par Finally let us show that it preserves the inverse. We have $b^{-1}=b_{2}^{-1}\times b_{1}^{-1}$ and the image of $b^{-1}$ in $H$ is equal to
$f_{\nabla_N(b_2^{-1})}\circ f_{\nabla_N(b_1^{-1})}$. Since $\nabla_N(b_1^{-1})=F_{-1}[N](\nabla_{N}(b_{1}))$ and that in $Perm_{g}(Y)$, we have that $f_{F_{-1}[N](\nabla_{N}(b_{1}))}=f_{\nabla_N(b_1)}^{-1}$ (property (g2)), we get the result.

\qed
\cl \label{cl:union} Any $\cL(k)$-generic type in $G$ is realized by a tuple $\nabla_{N}(a)$ with $a\in \Gamma$.
\ecl
\pr Set $Z^{(0)}:=Z$;  the domain of $G$ is equal to the (disjoint) union $\bigsqcup_{i=1}^{2d+1} \{\f_{\nabla_N(a_{i})}\circ \f_{d}\colon d\in Z^{(i)}\}$, where $Z^{(i)}$, $1\leq i\leq 2d$, has been defined as $\{d\in Z\colon f_{\nabla_N(a_{i})}\circ f_{d}\not\sim f_{\nabla_N(a_{j})}\circ f_{e}$ for some $e\in Z^{(j)}$ and some $0\leq j<i \}$. Further we saw that we may identify the domain of $G$ with $\bigsqcup_{i=1}^{2d+1} \{\nabla_N(a_{i})\}\times Z^{(i)}$ (see Remark \ref{ident}).
 \par As already noted, $\nabla_{N}(\Gamma)\subset Z$. Let $p$ be an $\cL(k)$-generic type in $G$; it is realized by an element of the form $(\nabla_{N}(a_{i}),u)$ with $u\in Z^{(i)}$. Moreover by Claim \ref{cl:dim1}, for each $u\in Z$, there $1\leq j\leq d+1$ such that $\nabla_{N}(a_{j})$ is independent from $u$, so we may assume that $(\nabla_{N}(a_{i}),u)$ is generic in $Z\times Z$ and that $u$ is of the form $\nabla_{N}(g)$ with $g\in \Gamma$. By Claim \ref{cl0}, $\f_{\nabla_{N}(a_{i})}\circ \f_{\nabla_{N}(g)}= \f_{\nabla_{N}(F(a_{i},g))}$. \qed

 \medskip
 \par \noindent (2) The second statement follows from Lemma \ref{set}. \qed
\medskip
\rem \label{rem:largeZ} The set $(\{\nabla_N(a_{1})\}\times Z)$ is a large $\cL(k)$-definable subset of $\dom(G)$.
\par Indeed let $u\in Z$ and let $2\leq i\leq 2d+1$. First assume that $F_{\times}(F_{-1}(\nabla_{N}(a_{1})), \nabla_{N}(a_{i}))$ is independent from $u$. 
\par\noindent Then by (g2), $u\in X_{F_{\times}(F_{-1}(\nabla_{N}(a_{1})), \nabla_{N}(a_{i}))}$, by Claim \ref{cl0}, $\f_{F_{\times}(F_{-1}(\nabla_{N}(a_{1})), \nabla_{N}(a_{i}))}\circ \f_{u}=\f_{F_{\times}(F_{\times}(F_{-1}(\nabla_{N}(a_{1})), \nabla_{N}(a_{i})), u)}$ and $F_{\times}(F_{\times}(F_{-1}(\nabla_{N}(a_{1})), \nabla_{N}(a_{i})), u)\in Z$ (by Claim \ref{cl2}). 
So, 
\begin{align*}
\f_{\nabla_{N}(a_{i})}\circ \f_{u}= &(\f_{\nabla_{N}(a_{1})}\circ \f_{F_{\times}(F_{-1}(\nabla_{N}(a_{1})), \nabla_{N}(a_{i}))})\circ \f_{u}\\
= &\f_{\nabla_{N}(a_{1})}\circ (\f_{F_{\times}(F_{-1}(\nabla_{N}(a_{1})), \nabla_{N}(a_{i}))}\circ \f_{u})\\
=&\f_{\nabla_{N}(a_{1})}\circ \f_{F_{\times}(F_{\times}(F_{-1}(\nabla_{N}(a_{1})), \nabla_{N}(a_{i})), u)}.
\end{align*}
Now the set of elements $u\in Z$ such that $F_{\times}(F_{-1}(\nabla_{N}(a_{1})), \nabla_{N}(a_{i}))$ is independent from $u$
is large in $Z$.
\erem
\medskip
\par When $T$ is a theory of geometric fields, we may replace in the statement of Theorem \ref{thm:def} $Y_{\ell}$, by a $k$-variety intersected with an open definable set, with associated $D$-variety.
\par When $T$ is a geometric theory of dp-minimal topological fields,
 then we may replace $Y_{\ell}$ by a $\partial$-compatible cell with associated $D$-cell.

\subsection{Prolongations}\label{torsor}
In this subsection, we will place ourselves into two settings: either $T$ is dp-minimal not strongly minimal, or $T$ is an open theory of topological fields with the property that $T$ has a model which is a complete nondiscrete valued field of rank $1$.
\par Let $\U$ be a model of $T(\partial)$ and let $A$ be a differential subfield of $K$.
\par Recall that any $\cL(A)$-definable subset (in the field sort) can be partitioned into a finite union of $\cL(A)$-definable cells $C_{i}$ which are graphs of $\partial$-compatible $C^1$-correspondences ($\partial$-compatible cells) (see Propositions \ref{C1-decomp-dp}, \ref{C1-decomp}). 
\par Let us first define the prolongation of a $\partial$-compatible cell $C$. Then we will define the prolongation of a definable set given with a partition into $\partial$-compatible cells.
\dfn Let $C:=(\rho_{C}(C),h)$ with $C\subset K^n$, $\rho_{C}:K^n\to K^d$, a coordinate projection with $\rho_{C}(C)$ an open set, $d:=\dim(C)$,  $h$ is a $\partial$-compatible $C^1$-correspondence from $\rho_{C}(C)\to K^{n-d}$. Let $h=(h_{1},\ldots, h_{n-d})$ and $h^{\partial}:=(h_{1}^{\partial},\ldots, h_{n-d}^{\partial})$ with $h_{i}:\rho_{C}(C)\to K^{n-d}$.  First we define $\tau(h)(b,u)$ with $b\in \rho_{C}(C)$ and $u\in K^d$. 
\par Let $d(h_{j})_b(u):=\sum_{t=1}^d \partial_{x_t} h_{j}(b) u_t$ and so $d(h_{b})(u)=(d(h_{1})_{b}(u),\ldots, d(h_{n-d})_{b}(u))$. 
\par Then define $T(h)$ by
\[
T(h)(b,u)=(h(b), d {h}_{b}(u)),
\]
and
 $\tau(h)$ by
\[
\tau(h)(b,u)=(h(b),d {h}_{b}(u)+h^{\partial}(b)).
\]
\par Note $T(h)$ and $\tau(h)$ are $\cL(A)$-definable $C^0$-maps since the maps $\partial_{x_{t}} h$ are $\cL$-definable and continuous, $1\leq t\leq n$, as well as the map $h^{\partial}$.
\par We have $\tau(h)(b,\partial(b))=(h(b),\partial(h(b)))$ (since $
\partial$ is compatible with $h$).

The prolongation  $\tau(C)$ is defined as $\{((b,h(b)),(u, d {h}_{b}(u)+h^{\partial}(b)))\colon b\in \rho_{C}(C), u\in K^d\}$. (Note that in the special case where $C$ is an open subset of $K^n$, then $\tau(C)=C\times K^n$ and when $C$ is finite we use the convention that a finite set is an open subset of $K^0$.)
\edfn
\par Now let us consider the case of a definable set.
\dfn Let $W$ be an $\cL(A)$-definable subset of $\U$ together with a cell decomposition, namely $W=\bigsqcup_{i\in I} C_{i}$, where $C_{i}$ are $\partial$-compatible cells, $I$ is finite. 
Then the prolongation $\tau(W)$ of $W$ is defined as $\bigsqcup_{i\in I} \tau(C_{i})$. 
\edfn

\par Now following \cite[1.6]{M}, but replacing regular maps by definable maps, we extend definable maps on $\cL(A)$-definable sets to maps on their prolongations. 

\

\par Let $f: K^n\to K^m$ be an $\cL(A)$-definable map in $\U$, $f=(f_{1},\ldots,f_{m})$. By Propositions \ref{C1-decomp-dp}, \ref{prop:comp1}, there is a finite partition of the domain of each $f_{j}$, $1\leq j\leq m$, into $\partial$-compatible cells $C_{ji}=(\rho_{C_{ji}}(C_{ji}),h_{ji})$ where $h_{ji}$ is a $C^1$-$\partial$-compatible $\cL(A)$-definable correspondence and $\rho_{C_{ji}}(C_{ji})$ an open subset of $K^{\dim(C_{ji})}$, $i\in I(j)$ a finite set.

Further there are finitely many open disjoint subsets $W_{\ell}$ of $K^n$, finitely many $C^1$-$\partial$-compatible correspondences $g_{\ell}$ on $W_{\ell}$ such that for $x=(x_{1},\ldots, x_{n})\in W_{\ell}\cap C_{ji}$, $f_{j}(x)=g_{\ell}(\rho_{C_{i}}(x), h_{ji}(\rho_{C_{ji}}(x)))$. So, letting $d=\dim(C_{ji})$, $z=(z_{1},\ldots,z_{n})$, $d(f_{j})_x(z):=\sum_{i=1}^d \partial_{x_i} f_{j}(x) z_i$ since $f$ only depends on $\bar x:=(x_{1},\ldots, x_{d})$.
\par First $\tau(g_{\ell})(x,z)=(g_{\ell}(x), d(g_{\ell})_{x}(z)+g_{\ell}^{\partial}(x)$ with $x\in W_{\ell}$, $z\in K^n$ and we have already defined $\tau(h_{ji})(b,u)=(h_{ji}(b),d {h_{ji}}_{b}(u)+h_{ji}^{\partial}(b))$ with $b\in \rho_{C_{ji}}(C_{ji})$ and $u\in K^{n-d}$.
Then,
\[
\tau(f_{j})(x,z)=(f_{j}(x),d(f_{j})_{x}(z)+f_{j}^{\partial}(x)).
\]
\par Since $\partial$ is compatible with $f_{j}$, we have $\tau(f_{j})(x,\partial(x))=(f_{j}(x), \partial(f_j(x)))$. 
\par Moreover on $(W_{\ell}\cap C_{ji})\times K^d$, we have that $\tau(f_{j})=\tau(g_{\ell})\circ \tau(h_{ji})$.
If we express the above expression in terms of $g_{\ell}$ and $h_{ji}$ with , we get that $d(f_{j})_{(x)}(z)+f_{j}^{\partial}(x)=$
\[ \sum_{i =1}^d \big{(}(\partial_{x_i} g_{\ell})(\bar x,h_{ji}(\bar x)) +\sum_{k=1}^{n-d} (\partial_{x_{k+d}} g_{\ell})(\bar x,h_{ji}(\bar x))\partial_{x_i} h_{jik}(\bar x) \big{)}  z_i+\sum_{k=1}^{n-d} h_{jik}^{\partial}(\bar x)+g_{\ell}^{\partial}(\bar x,h_{ji}(\bar x)),
\]
where $h_{ji}(\bar x)=(h_{ji1}(\bar x,\ldots,h_{ji(n-d)}(\bar x))$.

\bigskip
\par Let us recall the following result which can be found in \cite{M}.
\fct \cite[section 2]{M}\label{marker} Let $G$ be an algebraic group with multiplication $m$, the prolongations $T(G)$ and $\tau(G)$ are algebraic groups with multiplication $T(m)$, $\tau(m)$ respectively. Moreover the differential section $\nabla: G\to \tau(G): g\mapsto (g,\partial(g))$ is a group homomorphism: we have $\tau(m)\circ \nabla_{G\times G}=\nabla_{G}\circ m$.
\efct
\par On $T(G)$, Marker gave an explicit formula for the group law on $T(G)$. Suppose $G\subset K^n$. Let $g\in G$ and let $\lambda^g, \rho^g$ be respectively left and right multiplication by $g$ on $G$. For $g, h\in G$ and $v:=(v_{1},\ldots, v_{n})\in K^n$, let $d\lambda_{h}^g v=(\sum_{i=1}^n \partial_{x_i} \lambda_{1}^g(h) v_{i}, \ldots, \sum_{i=1}^n \partial_{x_i} \lambda_{n}^g(h) v_{i}).$
\par Then given $(g,u),\; (h,v)\in T(G)$, we have $(g,u) (h,v)=(gh, d\lambda_{h}^g v+d \rho_{g}^h u)$.\\ For $(g,u),\; (h,v)\in \tau(G)$, we have $(g,u) (h,v)=(gh, d\lambda_{h}^g v+d \rho_{g}^h u+(\lambda^g)^{\partial}(h)+(\rho^h)^{\partial}(g))$.
\medskip
\par Proceeding exactly in Fact \ref{marker} (see also \cite{PPP2}), and using the discussion above, we obtain the following proposition.
\prop \label{prop:torsor} Let $\U$ be a sufficiently saturated model of $T(\partial)$ and let $(G,F_{\times},F_{-1})$ be a $\cL(A)$-definable group in $K$, where $A$ is a small differential subfield of $K$. 
\par Then $T(G)$, respectively $\tau(G)$, can be endowed with an $\cL$-definable topological group structure corresponding to the $\cL$-definable continuous maps $T(F_{\times}), T(F_{-1}), \tau(F_{\times}), \tau (F_{-1})$. Moreover the differential section $\nabla: G\to \tau(G): g\mapsto (g,\partial(g))$ is a group homomorphism.
\eprop
\pr By Proposition \ref{prop:topo}, $G$ is a $C^1$-group, namely it can be endowed with the $\cL$-definable structure of a topological group where in addition the group operations are $\partial$-compatible $C^1$-maps. Then as in Fact \ref{marker}, we use that the property that $T(G)$, respectively $\tau(G)$ are groups, is translated through the commutativity of certain diagrams.
\qed
\dfn \label{D-groups} Let $(G,F_{\times},F_{-1})$ be a $C^1$-group and suppose there an $\cL$-definable section $s: G\to \tau(G)$, then the pair $(G,s)$ is called an $\cL$-definable $D$-group.
\edfn
\subsection{$D$-groups} Now let us state the analog of \cite[Theorem 4.10]{PPP2}, stated for $p$-adically closed fields.
\medskip
\thm \label{thm:D-group} Let $T$ be either a model-complete geometric dp-minimal not strongly minimal theory and assume that $T(\partial)$ has a model-companion $T_\partial$, or let $T$ be a model-complete open theory of topological fields with the property that $T$ has a model which is a complete nondiscrete valued field of rank $1$.
\par Let $\cG:=(\G,\times ,{\;}^{-1})$ be a finite-dimensional  $\cL_{\partial}(k)$-definable group in $\U\models T_\partial$, with $\G\subset K^n$ and $k$ a small differential subfield of $K$.
Then there is an $\cL$-definable $D$-group $(G,s)$, a natural number $N$ and an $\cL_\partial$-embedding from $\nabla_{N}(\Gamma)$ to $G$ whose image is $\{g\in G\colon s(g)=\nabla(g)\}$.
\ethm
\pr In Proposition \ref{largeV}, we showed that for some natural number $N$, there is an $\cL(k)$-definable pre-group $(Y,F_{-1}[N], F_{\times}[N])$, in which $\nabla_N(\G)$ embeds densely. 
\par By Lemma \ref{set}, $Y$ is a finite union $\cL(k)$-definable sets $Y_{\chi,\varphi}$ with $(\chi,\varphi)\in \cF(\nabla(Y))$ and is equipped with
an $\cL(k)$-definable section $s:Y\to \tau(Y)$ with the property that  $\nabla_{N}(\G) = \bigcup_{(\chi,\varphi)\in \cF(\nabla(Y))} \{a\in Y_{\chi,\xi}: \nabla(a) = s(a)\}$.
\par In Theorem \ref{thm:def}, we made explicit the construction of the associated $\cL(k)$-definable group $G$; in Claim \ref{cl:union}, we identified the domain of $G$ with the (disjoint) union: $\bigsqcup_{I=1}^{2D+1}\{\f_{\nabla_N(a_{i})}\circ \f_d: d\in Z^{(i)}\}$, where $Z^{(i)}$, $1\leq i\leq 2d$, has been defined as $\{d\in Z\colon 
\f_{\nabla_N(a_{i})}\circ \f_d\not\sim \f_{\nabla_N(a_{j})}\circ \f_e$ for some $e\in Z^{(j)}$ and some $0\leq j<i \}$, with $Z^{(0)}:=Z$, $Z$ a large subset of $Y$ (see Remark \ref{rem:largeZ}). Moreover, by Remark \ref{ident}, there is a bijection between 
$\{\f_{\nabla_N(a_{1})}\circ \f_d: d\in Z\}$ and $\{\nabla_N(a_1)\}\times Z$. 
\par By Propositions \ref{C1-decomp}, \ref{C1-decomp-dp}, an $\cL(k)$-definable set can be partionned into finitely many $\partial$-compatible $C^1$ $\cL(k)$-definable cells. In the proof of proposition \ref{prop:topo}, in order to put on $G$ the structure of a $C^1$-group, we could only consider the cells of maximal dimension (and so in this case the cells of maximal dimension that partition $Z$. By Proposition \ref{prop:torsor}, we have a differential section $\nabla$ from $G$ to $\tau(G)$ (which is a group homomorphism) and this differential section coincides with an $\cL$-definable section $s$ (see Examples \ref{ex:section} (B)). (We simply have to extend $s$ on $Z$ to $s$ on $\{\nabla_N(a_1)\}\times Z$). \qed 

\section{Appendix}
\subsection{Appendix 1}
\lem {\rm See \cite[Lemma 2.2]{simon-walsberg2016}} Let $\K$ be a topological field of characteristic $0$ with a geometric theory and suppose $\K$ is $\aleph_{0}$-saturated. Suppose that any infinite definable set $X$ in $\K$ with $X\subset K$, has non-empty interior.
 Then the acl-dimension and topological dimension coincide on $\K$.
 \elem
\pr By definition of these dimensions, it amounts to show that if $X\subset K^n$ is definable in $\K$, then $\dim(X)=n$ iff $X$ has non-empty interior.
\par The right to left direction is immediate. By induction on $n$, let us prove the left to right direction. If $n=1$, this is the hypothesis.
\par Assume $n>1$ and let $\pi$ be the projection onto the last $n-1$-coordinates. Let $X_{b}:=\{y\in K^{n-1}: (b, y)\in X\}$.
\par Let us show that there is $\bar c\in \pi(X)$ such that $\pi^{-1}(\bar c)$ is infinite. (Suppose not, then since
 $\dim(X)\leq \dim(\pi(X))+\max_{\bar c\in \pi(X)}\{\dim\{z\in K: (z,\bar c)\in X\}\leq n-1$, $\dim(X)\leq n-1$, a contradiction.)
 \par The same inequality shows, since $\dim(K)=1$, that $\dim(\pi(X))=n-1$ and so by induction hypothesis, $\pi(X)$ has non-empty interior. Further we may assume that $\bar c$ belongs to that interior.
 \par By the case $n=1$, the interior of $\pi^{-1}(\bar c)$ is not empty (and this is a definable condition on $\bar c$).
 \par So $Q:=\{b\in K: (b,\bar c)\in X$ and there is an open set $U$ containing $\bar c$ such that $\{b\}\times U\subset X\}$ is infinite.
 \par Let $b$ belong to the interior of $Q$ be such that $\{b\}\times U_{b}\subset X$, and consider $P_{U_{b}}:=\{b'\in Q: \{b'\}\times U_{b}\subset X\}$. The family of $P_{U_{b}}$ is a directed set and uniformly definable. Since $Q$ is infinite, for some $U$ (an open neighbourhood of $\bar c$), $P_{U}$ is infinite. Let $O$ be an open subset of $P_{U}$. Then $O\times U\subset X$, i.e. $X$ has non-empty interior.
\qed
\medskip
\par Let $f(x)$ be an $\cL(K)$-definable function, namely there is a $\cL$-definable function $F(x,y)$ with $f(x)=F(x,\bar c)$.  
\par Now we want to show that if a function is continuous on an open fiber indexed by a generic tuple, then it is continuous on an open set. 
\lem Let $g:W\to K$ with $W\subset K^m\times K^n$. Let $\pi$ be the projection on the last $n^{th}$ coordinates. Let $b\in \pi(W)$ be generic and suppose that $W_{b}:=\{x\in K^m: (x,b)\in W\}$ is open and that $g(x,b)$ is continuous on $W_{b}$. Let $(a,b)\in W_{b}$, then $(a,b)$ belongs to the interior of $W$ and $g$ is continuous at $(a,b)$. 
\elem
\pr Let $(a,b)\in W_{b}$. Let $V$ is an open set containing $g(a,b)$ and $U_{a}$ an open set containing $a$ such that $g(U_{a},b)\subset V$.
\cl There are some parameters $\bar d$ and an open subset $U_{\bar d}$ of $O_{a}$ containing $a$ and defined over $\bar d$ such that $b$ is generic over $\bar d$.
\ecl
Let $Y:=\{b'\in \pi_{n}(W): g(U_{\bar d},b')\subset V \}$. Since $b\in Y$ and $b$ is generic over $\bar d$, there is an open subset $O\subset K^n$ such that $b\in O$ and $O\subset Y$. So $g$ is continuous at $(a,b)$ and $(a,b)\in U_{\bar d}\times O\subset W$.
\medskip
\par Now let us show the claim. Since $K$ is dp-minimal (and not strongly minimal), the topology on $K$ is either given by a valuation or an ordering. Moreover taking $K$ sufficiently saturated we may assume that the topology on $K$ is given by a valuation $v$. Suppose that $O_{a}$ is given by $B(a,\gamma)$, where $\gamma$ is in the value group of $K$, namely $B(a,\gamma)=\{u\in K: v(a-u)>\gamma\}$. Then choose $c\in K$ such that $\gamma=v(c)$. So the former formula translates to: $v(a-u)>v(c)$, equivalently $div(c,a-u)$. First choose $c'$ such that $\gamma'=v(c')$ is such that $B(a, \gamma')\subset B(a,\gamma)$ and $b$ is still generic over $c'$. Then choose $a'\in B(a,\gamma')$ with $b$ generic over $c', a'$. Since $B(a',\gamma')=B(a,\gamma)$, we get the required open set $U_{\bar d}$. \qed
\subsection{Appendix 2}
Now let us discuss that the hypothesis for a valued field $(K,v)$ of characteristic $0$ to be elementary equivalent to a complete valued field of rank $1$. 
 The following is probably well-known, but we include it for sake of completeness.
 Let $k$ be the residue field of $(K,v)$ and let $\Gamma$ be its value group.
\lem \label{regular} Suppose that a valued field $(K,v)$ of characteristic $0$ is elementary equivalent to a complete valued field of rank $1$. Then in residue characteristic $0$,  this is equivalent to $(K,v)$ henselian and $\Gamma$ regularly ordered. In mixed characteristic, assuming further we are in the unramified case or in the finitely ramified one, then this is again equivalent to $(K,v)$ henselian and $\Gamma$ regularly ordered. 
\elem
\pr  Suppose $(K,v)$  is elementary equivalent to a complete valued field of rank $1$. Then $(K,v)$ is henselian (a complete valued field is henselian and being henselian is an elementary property). 
\par  In the following we will place ourselves in settings where the Ax-Kochen/Ershov theorem is applicable (there are many references and variants for that result; we have used a paper of Kochen \cite[Theorem 1 and section 8]{K} and a recent paper of Anscombe, Dittmann and Jahnke \cite{ADJ} where a lot of other references are given). 
\par First assume that $K$ is of residue characteristic $0$. Then the Ax-Kochen/Ershov theorem is applicable. Let $(K_0,v_0)$ be a complete valued field of rank $1$ with residue field $k_0$ and value group $\Gamma_0$. Then,
\[
K\equiv K_0 \leftrightarrow (k\equiv k_0\;\; \& \;\;\Gamma\equiv \Gamma_0), 
\]
where the elementary equivalence in the left hand side is in the language of rings enriched with a predicate for the valuation ring and in the right hand side, respectively in the language of rings and in the language of ordered abelian groups.
Since $\Gamma_0$ is of rank $1$, namely it is archimedean, to be elementary equivalent to $\Gamma_0$ is equivalent, by a result of A. Robinson and E. Zakon, to be regularly ordered \cite{RZ}.
\par Conversely if $\Gamma$ is regularly ordered, let $\Gamma_0$ be an archimedean group elementary equivalent to $\Gamma$. Consider the field of Hahn series $k((\Gamma_0))$. So by the Ax-Kochen/Ershov theorem, $K$ is elementary equivalent to $k((\Gamma_0))$. Note that this last field is maximally complete and so in particular it is complete.
  \par Second assume that $K$ is of residue characteristic $p$ and first assume that the residue field is perfect.
   \par In \cite[Corollary 3]{Poo}, given an ordered abelian group $\Gamma$ and a perfect field $k$, B. Poonen constructed a maximally complete valued field with value group $\Gamma$ and  residue field $k$. He called such a field a $p$-adic Mal'cev-Neumann. Since we wish to point out that his construction may be slightly generalized, let us quickly describe it. 
   \par Fix in $\Gamma$ a positive element, say $1$, and consider the cyclic subgroup of $\Gamma$ generated by $1$, say $\Z$, then construct the Witt ring $W(k)$, then the ring of Hahn series $W(k)((\Gamma))$ with coefficients in $W(k)$. In that ring, define the set $N$ of {\it null series} as the set of elements $\alpha\in W(k)((\Gamma))$ of the form $\sum_{\gamma\in \Gamma} a_\gamma t^{\gamma}$ with $a_\gamma\in W(k)$ such that for all $\gamma$ in the support of $\alpha$, $\sum_{n\in \Z} a_{\gamma+n} p^n=0$ in $W(k)$ (for instance $-p+t^1\in N$). One shows that the set $N$ is an ideal \cite[Proposition 3]{Poo}.) Then one observes that the quotient $W(k)((\Gamma))/N$ is a field $L$ with residue field $k$ and value group $\Gamma$ \cite[Corollary 3, Proposition 5]{Poo}; moreover $L$ is maximally complete (and so complete) \cite[Theorem 1]{Poo}. 
 \par In the case $k$ is not necessarily perfect, one replaces the Witt ring in the construction above by a subring: the Cohen ring $C(k)$ \cite{S}. (It is still complete.) Then one checks that the set $N'$ of series $\alpha'$ of the form $\sum_{\gamma\in \Gamma} a_\gamma t^{\gamma}$ with $a_\gamma\in C(k)$ such that for all $\gamma$ in the support of $\alpha'$, $\sum_{n\in \Z} a_{\gamma+n} p^n=0$ in $C(k)$ is an ideal of $C(k)((\Gamma))$. Then one verifies that the quotient is still a maximally complete field (and so complete) with residue field $k$ and value group $\Gamma$.
\par Now further assume that we are in the unramified or finitely ramified case, in order to apply again the Ax-Kochen-Ershow theorem. (In the unramified case, one can appeal to \cite[Theorem 1.2]{AJ} and in the finite ramification case, one has to require that the residue fields are elementary equivalent in languages expanding the language of rings \cite[Definition 3.7]{ADJ} for a detailed discussion. Note that in our setting it has no consequence since we don't change the residue field.)
\par As before, if $(K,v)$ is elementary equivalent to a complete valued field of rank $1$, then it implies that both value groups are elementary equivalent. So $\Gamma$ is regularly ordered. Conversely, suppose the value group $\Gamma$ of $K$ is regularly ordered. So, there is an archimedean group $\Gamma_{0}$ elementary equivalent to $\Gamma$. Then  the p-adic Mal'cev-Neumann fields  $W(k)((\Gamma_0))/N$ in case the residue field $k$ is perfect or $C(k)((\Gamma_0))/N'$ when $k$ is non-perfect, are, respectively, elementarily equivalent to $(K,v)$, by the Ax-Kochen/Ershov theorem. \qed
\rem 
Let $(K,v)$ be a valued field with residue field $k$. First suppose that $(K,v)$ is a model of $ACVF_{0}$. Note that in this case we have an Ax-Kochen/Ershov theorem (this follows from the quantifier elimination result of A. Robinson.). Either $K$ is of residue characteristic $0$ and so $K\equiv k((\Q))$ with $k$ the residue field of $K$, or $K$ has residue characteristic $p$ and so $W(k)((\Q)/N$ the p-adic Mal'cev-Neumann field is algebraically ckosed \cite[Corollary 4]{Poo} and so $K\equiv W(k)((\Q)/N$.
\par Second suppose that $(K,v)$ be a model of $RCVF$, then $k((\Q))$ is again a model of $RCVF$ and $K\equiv k((\Q))$.
\par Finally let $K$ be a real-closed field and assume it is $\aleph_{1}$-saturated. Then $(K,v_{a})$ is an henselian field with divisible group, where $v_{a}$ the valuation induced by the archimedean equivalence relation on $K$. So $(K,v_{a})$ is elementary equivalent to a complete valued field of rank $1$. Moreover the order topology and the valuation topology coincide.
\erem
\subsection{Appendix 3}
\par For convenience of the reader, let us show that if $\partial$ is compatible with two correspondences, then it is also compatible with their composition (see \cite[Lemma 2.6]{FK}).
\lem Let $h: K^n\rightrightarrows K$ be a $C^1$ $m$-correspondence  on $U$ an open subset of $K^n$. Let $g:K^{n+m}\to K$ be a $C^1$ $1$-correspondence on an open subset $V$ of $K^{n+m}$ containing $h(U)$.
Suppose that $\partial$ is compatible with $h$ and $g$. Then $\partial$ is compatible with the correspondence $f:K^n\rightrightarrows K: x\mapsto g(x,h(x))$ defined on $U$.
\elem
\par By Lemma \ref{function}, each $a\in U$ has an open $U_a\subset U$ such that there are $m$-continuous functions $h_i:U_a\to K$ with $\graph(h)=\bigsqcup_{i=1}^m \graph(h_i)$. By assumption on $h$, the $\partial$ derivatives of $h_i$, $1\leq i\leq m$, exist and are continuous on $U_a.$ W.l.o.g. we may assume that  $V_a:=\{(x,h_1(x),\ldots,h_m(x)): x\in U_a\}$ is included in $V$. 
\par By the chain rule for partial derivatives of $C^1$-functions, we have for $1\leq i\leq n$:
\[
\partial_{x_i} f(x)=(\partial_{x_i} g)(x,h(x))+\sum_{j=1}^m (\partial_{y_j} g)(x,h(x)) \partial_{x_i} h_j(x).
\]
Then, we consider the scalar product with $(\partial(x_1),\ldots,\partial(x_n))$ and we get:
\[
\sum_{i =1}^n \partial_{x_i} f(x)\partial(x_i)=\sum_{i=1}^n ((\partial_{x_i} g)(x,h(x))+\sum_{j=1}^m (\partial_{y_j} g)(x,h(x)) \partial_{x_i} h_j(x)) \partial(x_i).
\]
Since $\partial$ is compatible with $h$ and $h=(h_1,\ldots,h_m)$ on $U_a$, we have for each $1\leq j\leq m$,
\begin{equation}\label{eqh}
\sum_{i=1}^n \partial_{x_i} h_j(x) \partial(x_i)=\partial(h_j(x))-h_j^{\partial}(x).
\end{equation}
Since $\partial$ is compatible with $g(x,y)$ on $V$, we have
\[
\partial(g(x,y))=\sum_i (\partial_{x_i} g)(x,y) \partial(x_i)+\sum_j (\partial_{y_j} g)(x,y)\partial(y_j)+g^{\partial}(x,y).
\]
Evaluating the above expression at $(x,h_1(x),\ldots, h_m(x))$, we get
\[
\partial(f(x))=\partial(g(x,h(x)))=\sum_{i=1}^n (\partial_{x_i} g)(x,h(x)) \partial(x_i)+\sum_{j=1}^m (\partial_{y_j} g)(x,h(x))\partial(h_j(x))+g^{\partial}(x,h(x)).
\]

Replacing $\partial(h_j(x))$ (see equation (\ref{eqh}))
)), we get that 
\[
\partial(f(x))=\sum_i (\partial_{x_i} g)(x,h(x)) \partial(x_i)+\sum_j (\partial_{y_j} g)(x,h(x))(\sum_i \partial_{x_i} h_j(x)\partial(x_i)+h_j^{\partial}(x))+g^{\partial}(x,h(x)).
\]
\[
\partial(f(x))=\sum_i \partial_{x_i} f(x)\partial(x_i)+(\sum_j h_j^{\partial}(x))+g^{\partial}(x,h(x))).
\]
By assumption, $\sum_j h_j^{\partial}(x))+g^{\partial}(x,h(x))$ is $\cL$-definable and continuous. so we set $f^{\partial}:=\sum_{j=1}^m h_j^{\partial}(x))+g^{\partial}(x,h(x))$ \qed
\medskip
\par  We will show that it is closed by the implicit function theorem (it was already stated in \cite[Lemma 2.7]{FK} but there, one implicitly assumed that everything was defined over $\0$). 
\lem \label{implicit} Let $(K,\partial)$ be a differential topological field and $A$ an $\cL_{\partial}$-substructure. Let $U$ be an open definable subset of $K^{n}$ and $V$ an open definable subset of $K$. Let $f(\bar x, y)$, $\bar x=(x_1,\ldots,x_n)$, $y$ a single variable, be an $\cL(A)$-definable $C^1$-1-correspondence from $U\times V$ to $K$. Suppose that $\partial$ is compatible with $f$. Let $g(\bar x)$ be an $\cL(A)$-definable $C^1$-1-correspondence  from $U\to K$.
Consider the 1-correspondence $h:K^n\to K:\bar x\mapsto f(\bar x,g(\bar x))$ and assume that  $h\restriction U=0$ and $\frac{\partial}{\partial y} f(\bar x,g(\bar x))\restriction U \neq 0$.Then $\partial$ is again compatible with $g$.
Furthermore if $f$ and $g$ are $C^k$, then $g^{\partial}$ is $C^{k-1}$, $k>0$.
\elem
\par Since $\partial$ is compatible with $f$ we have: for any tuple $\bar b\in U\times V$ that
\begin{equation}\label{eq1}
\partial(f(\bar b))=\sum_{i=1}^{n} \frac{\partial}{\partial x_i} f(\bar b) \partial(b_i)+\frac{\partial}{\partial y} f(\bar b) \partial(b_{n+1})+f^{\partial}(\bar b).
\end{equation}

\par 
By the chain rule, for $\bar x\in U$, we have for each $1\leq i\leq n$,
\[
(\frac{\partial}{\partial x_i} f)(\bar x, g(\bar x)) +(\frac{\partial}{\partial y} f)(\bar x,g(\bar x)) \frac{\partial}{\partial x_i} g(\bar x)=0
\]
So taking the scalar product of $(\partial(x_1),\ldots,\partial(x_n))$ with 
\par $(\frac{\partial}{\partial x_1} f(\bar x,g(\bar x)) +\frac{\partial}{\partial y} f(\bar x,g(\bar x)) \frac{\partial}{\partial x_1} g(\bar x),\ldots, \frac{\partial}{\partial x_n} f +\partial_y f(\bar x,g(x)) \frac{\partial}{\partial x_n} g(\bar x))$, we get
\begin{equation}\label{eq2}
\sum_{i=1}^n \frac{\partial}{\partial x_i} f(\bar x,g(\bar x)) \partial(x_i)+\frac{\partial}{\partial y} f(\bar x,g(\bar x))\sum_{i=1}^n \frac{\partial}{\partial x_i} g(\bar x) \partial(x_i)=0.
\end{equation}
\par Let $\bar a\in U$. Since $f(\bar a,g(\bar a))=0$, we have $\partial(f(\bar a,g(\bar a))=0$. So we have using equation (\ref{eq1})
 \[
 \sum_{i=1}^{n} \frac{\partial}{\partial x_i} f(\bar a,g(\bar a)) \partial(a_i)+\frac{\partial}{\partial y} f(\bar a,g(\bar a)) \partial(g(\bar a))+f^{\partial}(\bar a,g(\bar a))=0.
 \] 
\par Then, using equation (\ref{eq2}), we replace  $\sum_{i=1}^{n} \frac{\partial}{\partial x_i} f(\bar a,g(\bar a)) \partial(a_i)$ by $-\frac{\partial}{\partial y} f(\bar a,g(\bar a))\sum_{i=1}^n \frac{\partial}{\partial x_i} g(\bar a) \partial(a_i)$.
\par So we get 
\[
-\frac{\partial}{\partial y} f(\bar a,g(\bar a))\sum_{i=1}^n \frac{\partial}{\partial x_i} g(\bar a) \partial(a_i)+\frac{\partial}{\partial y} f(\bar a,g(\bar a)) \partial(g(\bar a))+f^{\partial}(\bar a,g(\bar a))=0
\]

Now we assumed that $\frac{\partial}{\partial y} f(\bar a,g(\bar a)\neq 0$ and
we obtain, after dividing by $\frac{\partial}{\partial y} f(\bar a,g(\bar a)$:
\begin{equation}\label{eq_implicit}
\partial(g(\bar a))=\sum_{i=1}^n \frac{\partial}{\partial x_i} g(\bar a) \partial(a_i)- \frac{ f^{\partial}(\bar a,g(\bar a))}{\frac{\partial}{\partial y} f(\bar a,g(\bar a))}.
\end{equation}

\par Denote by $g^{\partial}(\bar x)$ the function $- \frac{ f^{\partial}(\bar x,g(\bar x))}{\frac{\partial}{\partial y} f(\bar x,g(\bar x))}$. Note that it is only defined in a open neighbourhood of $\bar a$ where 
$\frac{\partial}{\partial y} f(\bar x,g(\bar x))\neq 0$ (and by our assumption, on $U$). Also as a quotient of two continuous functions, $g^{\partial}$ is continuous and it is $\cL(A)$-definable since $g, f^{\partial}$ are (and the topology on $K$ is definable and so partial derivatives are definable). Moreover, if $\partial(A)=0$, then $g^{\partial}=0$, since it holds for $f$.\qed

\end{document}